\documentclass[a4paper,11pt, reqno]{amsart}
\usepackage[T1]{fontenc}
\usepackage[english]{babel}

\usepackage{amsmath}
\usepackage{amssymb}
\usepackage{amsthm}
\usepackage[foot]{amsaddr}
\usepackage{bbm}
\usepackage{mathrsfs}
\usepackage{epsfig}
\usepackage{enumerate}
\usepackage{mathtools}
\usepackage{graphicx}
\usepackage{tabularx}
\newcolumntype{Y}{>{\centering\arraybackslash}X}
\usepackage[usenames,dvipsnames]{xcolor}
\usepackage[colorlinks=true,linkcolor=NavyBlue,urlcolor=RoyalBlue,citecolor=PineGreen,bookmarks=true,bookmarksopen=true,bookmarksopenlevel=2,unicode=true,linktocpage, pdftitle={Renormalised Amperean area of Brownian motions and Symanzik representation of the 2D abelian Yang--Mills--Higgs field}, pdfauthor={Isao Sauzedde}]
{hyperref}
\usepackage[margin=1in]{geometry}
\usepackage{float}

\usepackage[colorinlistoftodos]{todonotes}

\newtheorem*{metatheorem*}{Metatheorem}
\newtheorem*{lemma*}{Lemma}

\newtheorem{theorem}{Theorem}[section]
\newtheorem{lemma}[theorem]{Lemma}

\newtheorem{proposition}[theorem]{Proposition}
\newtheorem{corollary}[theorem]{Corollary}

\theoremstyle{remark}

\newtheorem{remark}{Remark}

\newtheorem{conjecture}{Conjecture}

\def\d{\,{\rm d}}

\def\n{{\mathbf n}}

\def\conv{*}

\DeclareMathOperator{\grad}{grad}
\renewcommand{\div}{\operatorname{div} }
\DeclareMathOperator{\curl}{curl}

\DeclareMathOperator{\Range}{Range}
\DeclareMathOperator{\Supp}{Supp}

\def\XXint#1#2#3{{\setbox0=\hbox{$#1{#2#3}{\int}$}
     \vcenter{\hbox{$#2#3$}}\kern-.5\wd0}}

\title[Renormalised Amperean area of Brownian motions and Symanzik representation]{Renormalised Amperean area of Brownian motions and Symanzik representation of the 2D abelian Yang--Mills--Higgs field}
\author{Isao Sauzedde$^\ast$}
\address{$^\ast$ University of Warwick}
\email{isao.sauzedde@warwick.ac.uk}
\keywords{Brownian motion, winding numbers, Euclidean quantum field}
\subjclass[2020]{60J65,60J55,81T13, 60H30}

\begin{document}
\begin{abstract}
	We construct and study the renormalised Amperean area of a Brownian motion. 
	First studied by W.Werner, the Amperean area is related to Lévy area and stochastic integrals in a way akin to the relation between self-intersection measure and occupation measure. As we explain, it plays a central role in the Symanzik's polymer representation of the continuous Abelian Yang--Mills--Higgs field in 2 dimensions and allows to study this field using classical stochastic calculus and martingale theory.
	
%
\end{abstract}
\maketitle
\setcounter{tocdepth}{1}

\section{Introduction}
Let $W,W'$ be two planar Brownian motions with durations $T,T'$, and extends these oriented curves into loops by adding to each of them a straight line segment between its endpoints. We  continue to call $W,W'$ these extensions, which are parametrised by $[0,\bar{T}]$, $[0,\bar{T}']$ for some $\bar{T}>T, \bar{T}'>T'$.  Consider then the winding functions
\[ 
\n_{W }: \mathbb{R}^2\setminus \Range(W ) \to \mathbb{Z}  , \qquad
\n_{W'}: \mathbb{R}^2\setminus \Range(W') \to \mathbb{Z},  
\]
which to a point $z$ map the number of time the corresponding loop winds around $z$. 
These functions are unbounded, and not even locally integrable, in the vicinity of $\operatorname{Range} (W)$ (resp. $\operatorname{Range} (W')$), and the integrals
\[ 
\mathcal{L}_{W }\coloneqq \int_{\mathbb{R}^2} \n_{W }\d \lambda, 
\qquad \tilde{\mathcal{A}}_{W}\coloneqq\int_{\mathbb{R}^2}  \n_{W }^2  \d \lambda- \mathbb{E}[\int_{\mathbb{R}^2}  \n_{W }^2  \d \lambda] 
, 
\quad \mbox{and} \quad
\mathcal{B}_{W,W'}\coloneqq\int_{\mathbb{R}^2}  \n_{W } \n_{W'}  \d \lambda
\]  
with respect to the Lebesgue measure  $\d \lambda$ are thus ill-defined as Lebesgue integrals.

However, the first one can still be defined by several different regularisation methods which give rise to the same limit (see \cite{LAWA,Werner2}), indeed equal to the L\'evy area delimited by $W$, that is $\int_0^{\bar{T}} W^1 \circ \d W^2$. Such an equality can be understood as a formal application of Stokes' theorem.\footnote{Stokes' theorem is almost always formulated for loops \emph{without self-intersections}. In this special case, the winding function is, up to a sign that depends on the orientation of the loop, just the indicator function of the domain delimited by the loop.}

The integral ${\mathcal{A}}_{W}=\int_{\mathbb{R}^2}  \n_{W }^2  \d \lambda$ has been coined the \emph{Amperean area} of the loop. Wendelin Werner showed in \cite{Werner3}, from which this work is much inspired, that upon deleting from the range of integration the points at distance less than $\epsilon$ from $W$, one gets a finite random variable, which diverges logarithmically fast as $\epsilon$ goes to $0$.\footnote{In \cite{Werner3}, a Brownian loop is considered rather than a Brownian motion with free endpoints. Here we only consider the case of Brownian motions, but most of our results can easily be transferred to the case of fixed endpoints. The two situations (fixed or free endpoints) are qualitatively the same, and the difficulties to transfer results between them are of purely technical nature.} 

Our main results can be roughly summarised as follows.
\begin{metatheorem*}
Let $\n_W^\epsilon, \n_{W'}^\epsilon$ be the `mollications' of $\n_W,\n_{W'}$, defined below. Let the random variables 
\[\mathcal{A}^\epsilon_{W}\coloneqq\int_{\mathbb{R}^2} (\n^\epsilon_W)^2 \d \lambda, 
\qquad 
\tilde{\mathcal{A}}^\epsilon_{W}\coloneqq\mathcal{A}^\epsilon_{W}-\mathbb{E}[\mathcal{A}^\epsilon_{W}], 
\qquad 
\mathcal{B}_{W,W'}^\epsilon\coloneqq\int_{\mathbb{R}^2}  \n_{W }^\epsilon \n_{W'}^\epsilon  \d \lambda.
\]
Then, as $\epsilon\to 0$, 	$\tilde{\mathcal{A}}^\epsilon_{W}$, $\mathcal{B}^\epsilon_{W,W'}$ converge in probability and with some positive and negative exponential moments. The limits do not depend on the mollifier. The limit of the latter can be expressed as the sum of an iterated Itô integral, plus a quarter of the intersection local time between $W$ and $W'$.  
As $\epsilon \to 0$, $\mathbb{E}[\mathcal{A}^\epsilon_{W}]$ admits the asymptotic expansion 
\[
\mathbb{E}[\mathcal{A}^\epsilon_{W}]=\frac{1}{2\pi} \log(\epsilon^{-1}) +C+O(\epsilon^{\frac{1}{2}}), 
\]
where the constant $C$ depends on the mollifier. 
\end{metatheorem*} 

\medskip 
We will first explain, at the formal level, how the (yet ill-defined) integrals $\tilde{A}_W$ and $\mathcal{B}_{W,W'}$ are related to the so-called abelian Yang--Mills--Higgs (aYMH) field. More precisely, we explain how the polynomial moments of a certain large family of observables, the loop and string observables, can formally be expressed in term of exponential moments of the Amperean area between Brownian loops and bridges.  
For an excellent overview on the literature concerning the Yang--Mills--Higgs field, we refer to the introduction section in \cite{chandra}. The approach followed here is inspired mostly from the celebrated Symanzik's polymer representation \cite{Symanzik}, from Albeverio and coauthor's work in general, and from \cite{AlbeverioKusuoka} in particular. For the role played by the loop observables in the characterisation of the Yang--Mills field, and for the construction of the non-Abelian Yang--Mills field, we refer to Lévy's monograph \cite{Thierry}. 

The second and most important aspect of this work is to give mathematically rigorous definitions of these quantities $\mathcal{B}_{W,W'}$ and $\tilde{\mathcal{A}}_W$ using some regularisation and renormalisation (respectively), and to prove some of their most important properties in regard of their connection with the aYMH field, as summarised by the metatheorem above. 

The regularisation method used, relying mostly on an adequate mollification of the winding function $\n_W$, is compatible with a rigorous interpretation of the relation between these quantities and the aYMH field: 
it amounts exactly to a mollification of the abelian Yang--Mills field (the Higgs field is not mollified). The counterterm involved in the definition of the recentred Amperean area $\tilde{\mathcal{A}}_W$, i.e. the expectation of the mollified version of $\mathcal{A}_W$, can also be directly identified as a simple mass term in the aYMH model. A particularly striking fact is that if we subtract to $\mathcal{A}_W$ a specific multiple of the self-intersection local time, the counterterm needed to define both cancel each other up to a non-diverging and non-random constant. As we will see, this is not coincidental.

\medskip 

The relation between the aYMH field and the renormalised Amperean area can therefore be understood in two ways: one can have a \emph{formal} interpretation, in which we rigorously construct the renormalisated Amperean areas, (Step 1a), 
do formal computation to formally relate them with the formally defined aYHM field (1b), and finally \emph{define} the aYMH field by imposing that the moments of the loop and string observables are given by the corresponding exponential moments of the Amperean areas, relying on a result of characterisation by moments (1c). 

On the other hand, one can first define a mollified version of the aYHM field (2a), prove that the moments of its loop and string observables are given by exponential moments of mollified Amperean areas (2b), and deduce that the mollified aYHM field converges (2c), indeed toward the aYHM defined by the first approach or otherwise. 

Let us be clear that none of these two approaches are carried here: we only carry the steps (1a) and (1b). What prevents from applying the step (1c) is essentially that we lack a Nelson-type bound: we do not prove that \emph{all} the negative exponential moments of the Amperean area are finite, but only that \emph{some of them} are,{\footnote{For comparison, the method  which allows to prove the finiteness of the negative exponential moments for the self-intersection local time does not apply here. Indeed, this method relies on the positivity of the intersection local time (between different curves), when on the contrary $\mathcal{B}_{W,W'}$ is not positive. Even for a fixed and non-zero positive value of the regularisation parameter, I cannot prove these negative exponential moments are all well-defined (in the case of the recentred self-intersection local time, this is trivial from boundedness from below).} nor that the corresponding family of  moments is indeed the algebra of moments of one and only one probability measure. Remark also that the string observables play a role similar to the $k$-point kernel of a scalar field, thus should not be expected to form genuine marginals of the field (in the same way $G(x,y)``=\mathbb{E}[\Phi(x)\Phi(y)]"$ is a well-defined quantity but the Gaussian free field $\Phi$ is not actually defined pointwise).
For the second approach, we expect that the steps (2a) and (2b) are easy to carry. The point (2c) would again require to show all the negative exponential moments are finite. 

Although more work is needed to complete any of these two approaches, it is to be expected that the limiting field would match with the aYMH field obtained with other constructions, and in particular it gives genuine expressions for the polynomial moments of string and loop observables, in term of exponential moments of some stochastic integrals computed along Brownian paths. As far as the author knows, these expressions, which are sufficiently explicit that they could be numerically evaluated, were not previously known.

\medskip 

More specifically, what is achieved here is the following.  The Amperean area $\mathcal{B}_{W,W'}$ between two independent Brownian motions $W,W'$ is constructed, and an explicit formula for  $\mathcal{B}_{W,W'}$, as an iterated Stratonovich integral, is obtained. A second formula involving the intersection local time between $W$ and $W'$ is given. 
Using this formula, it is proved that $\mathcal{B}_{W,W'}$ admits some finite exponential moments (both positive and negative), and that the convergence of the regularised area toward $\mathcal{B}_{W,W'}$ holds not only in $L^2$ but with some exponential moments.

The construction of the recentred Amperean area $\tilde{\mathcal{A}}_{W}$ is then achieved, 
and it is proved it also admits some finite exponential moments (both positive and negative). The way the counterterm diverges as the regularisation parameter $\epsilon$ goes to $0$, is described up to a $o(1)$. The main order term, proportional to $\log(\epsilon)$, matches with that obtained in \cite{Werner3} although the regularisation methods here and in \cite{Werner3} differ.\footnote{With insight, it is intuitively clear that the counterterms obtained by using either our regularisation method or the cut-off method in \cite{Werner3} should differ by a $O(1)$; but a rigorous proof seems difficult to obtain. Our proof does not use the result in \cite{Werner3}.}
	

From 
results that we obtained with P.~Perruchaud in \cite{zeta}, it is to be expected that the partition function $Z^\epsilon$ of the $\epsilon$-mollified aYMH field can also be described as the expectation $Z^\epsilon=\mathbb{E}[\mathcal{Z}^\epsilon]$ of an infinite product $\mathcal{Z}^\epsilon$ which involves the mollified Amperean area between the loops of a Brownian loop soup (this will be explained in the next section). 
In consequence of our lack of exponential moment, such a formula cannot be proven here. 


If such exponential moments could be proven to be finite indeed, which might require the addition of a $\Phi^4$ term, it would certainly lead to a construction of the aYMH field much simpler than the state-of-the-art construction, 
and might also lead to a construction in the non-Abelian case. 

Other possible improvements include: 
\begin{itemize} 
	\item The addition of a $\Phi^4$, or $P(\Phi)$, interaction term, which should lead to similar formulas but with extra terms involving the renormalised power of the local times of the involved Brownian paths, see e.g. the formal expression in \cite[p.134]{AlbeverioKusuoka}. We are currently working with P. Perruchaud on a variation of this, for fixed $\epsilon>0$. The approach here is adapted to the addition of such extra terms, as they both rely on Symanzik's polymer representation.	The `true' YMH field does contain such a term. Yet the relation that we prove between the Amperean area and the intersection local time vaguely suggests that the interaction term between the Higgs field and the Yang--Mills field might contain a sort of \emph{built-in} $\Phi^4$-term.
		
\item The extension to general surfaces endowed with a Riemannian metric (and to general fibre bundles over them, would the non-abelian case be treated as well). This is in principle (but not in practice) easy, as the extra difficulties are rather well-understood and sort of orthogonal to the ones we treat here. 
\item The generalisation to $3$ dimensions. What the author expects is that the Amperean area continues to be well-defined (although it is not an \emph{area} any longer), but does not even admit second moments,\footnote{This is related to the fact that for $G_d$ the Green function in dimension $d$, $(x,y)\mapsto G_d(x)^2G_d(y)^2G_d(x-y)^2$ is locally integrable near $(0,0)\in \mathbb{R}^d\times \mathbb{R}^d$ for $d=2$ but not for $d=3$ (the Amperean area between independent Brownian motions with \emph{different} starting points should be in $L^2(\Omega)$ in dimension $d=3$, but might not have higher moments).} so that extra counterterms are needed to define the exponential moments.   
\end{itemize}

In Section \ref{sec:aYMH}, the formal relation between the Yang--Mills--Higgs field and the Amperean area of Brownian paths is established. It contains little to no rigorous mathematical result (this is abstract, symbolic, computation), but serves as the main motivation for this work. That section, which is closer from mathematical physics than from probability, is independent from the rest of the paper. Section \ref{sec:main} contains the rigorous statement of the main results in the paper. Section \ref{sec:not} introduces the notations used during the proofs. Section \ref{sec:smooth}, which contains only deterministic computation, gives a formula for the Amperean area between two smooth loops, akin to the Stokes' formula for the area delimited by one loop. This formula appears in the physics literature but with a proof that relies on ill-defined expressions during intermediate steps. 
Section \ref{sec:epsilon} and \ref{sec:zero} are devoted to the construction of the Amperean area $\mathcal{B}_{W,W'}$ between two loops: first some results for a fixed positive value of the mollification parameter $\epsilon$, then asymptotic results as $\epsilon \to 0$. In Section \ref{sec:expo}, the existence of finite exponential moments for $\mathcal{B}_{W,W'}$ is proved. In Section \ref{sec:X}, we finally build the recentred Amperean area $\tilde{\mathcal{A}}_W$. In Section~\ref{sec:expoX}, we prove it admits some exponential moments.  
In Section \ref{sec:average}, we estimate the counterterm $\mathbb{E}[\mathcal{A}_W]$, asymptotically as the mollification parameter goes to $0$.
 

\section{Relation between the abelian Yang--Mills--Higgs field and the Amperean area}
\label{sec:aYMH}

The abelian Yang--Mills--Higgs measure with self-interaction term $\lambda$, coupling constant $\alpha$, and mass $m$, on $D\subset \mathbb{R}^2$ a planar domain, is the formally defined probability measure $\mathbb{P}$ on couples $(A,\Phi)$, where $\Phi:D\to \mathbb{C}$ and $A:D\to \mathbb{R}^2$, given by 
\[ 
\d \mathbb{P}(\Phi,A)\coloneqq \frac{1}{Z}\exp(- \frac{S(\Phi,A)  }{2}   ) \mathcal{D} \Phi \mathcal{D} A,
\]
where 
\[
S(\Phi,A)\coloneqq \| \curl A \|_{L^2(D,\mathbb{R}) }^2 +S_A(\Phi), \]\[  
S_A(\Phi)\coloneqq
\| \grad \Phi+i \alpha A \Phi \|_{L^2(D,\mathbb{C}^2) }^2 + 
m \|\Phi\|_{L^2(D,\mathbb{C}) }^2+
\lambda \|\Phi\|_{L^4(D,\mathbb{C}) }^4.
\]
One could also consider, rather than a planar domain $D$, a surface $\Sigma$, with or without boundary, endowed with a Riemannian metric. Then, $A$ becomes a connection  over some complex line bundle $F$ over $\Sigma$, whilst $\Phi$ is a section of $F$. Now we only consider $\lambda=0$.  

Let us separate the terms which depend on $\Phi$ from those which do not (and recall for the time being, we are only doing formal computations):
\begin{equation}
\label{eq:defaYMH}
\d \mathbb{P}(\Phi,A)=  \frac{Z' Z_A}{Z} \d \mathbb{P}_A(\Phi) \d \mathbb{P}_{YM}(A),
\end{equation}

\[
\d \mathbb{P}_A(\Phi)\coloneqq 
\frac{1}{Z_A} \exp(- \frac{S_A(\Phi)}{2})  \mathcal{D}\Phi\ \quad  
 \d \mathbb{P}_{YM}(A)\coloneqq \frac{1}{Z'} \exp(- \frac{ \| \curl A \|_{L^2(\mathbb{R}^2,\mathbb{R}) }^2}{2}) \mathcal{D} A,
\]
where the formal constants $Z', Z_A$ are chosen such that the formal measures $\mathbb{P}_A$ and $\mathbb{P}_{YM}$ are probability measures. 

For any given smooth $A$, the formal measure $\mathbb{P}_A$ and the constant $Z_A$ can rigorously be constructed as follow. 
Let $\Delta_A$ be the Dirichlet magnetic Laplacian, associated with $A$ and with the mass $m$, and acting on complex-valued functions which vanishes on $\partial D$: 
\[ 
\Delta_A=\frac{1}{2}(\nabla+i\alpha A)^*(\nabla+i\alpha A)+m I_d.
\]
It is associated with the Dirichlet form $S_A$: $\langle \Phi, \Delta_A\Phi\rangle_{L^2(D, \mathbb{C})}= S_A(\Phi)$.
Provided $m$ is sufficiently large, which we now assume, all the eigenvalues of this operator $\Delta_A$ are positive ($m\geq 0$ suffices but is not necessary). Let $G_A$ be the corresponding Green function. 
The measure $\mathbb{P}_A$ is then defined as the centred Gaussian measure, on an adequate space of complex-valued distributions, with covariance function $G_A$. Furthermore, from this Gaussian interpretation of $\mathbb{P}_A$, the partition function $Z_A$  should formally be given as $\sqrt{\det(\Delta_A)}$. This determinant is in fact well-defined through the so-called zeta-regularisation method (see e.g. \cite{Rosenberg}), which we can take as a definition for $Z_A$. We proved with P.Perruchaud that the $\zeta$-regularised determinant of a magnetic laplacian (or its non-abelian generalisation) can be written as a product over the loop soup, which is an expression well-suited to the framework of isomorphism theorems and Symanzik's loop representation. In the case we consider here, the general formula we obtained in \cite{zeta} reduces to\footnote{ The proportionally constant in \eqref{eq:zeta} (which is equal to $Z_0=(Z_A)_{|A=0}$) is irrelevant here, as it is absorbed by the global normalisation factor $Z$ in \eqref{eq:defaYMH}.} 
\begin{equation}
\label{eq:zeta}
Z_A\propto \mathbb{E}\big[ \prod_{\ell \in \mathcal{L}} \exp( i \alpha \int A_\ell \circ \d \ell)  )  \big],
\end{equation}
where $\mathcal{L}$ is an oriented massive Brownian loop soup on $D$ with intensity $1/2$ and mass $m$, and the product is ordered in any way independent from the orientations of the loops. We use the notation $\circ \d $ for Stratonovich integrals, and $\d$ for It\^o integrals.

Expanding the Green function as a time integral of the heat kernel $p^A$ of $\Delta_A$, and using the magnetic Feynman-Kac's formula (see e.g. \cite[Proposition 4.1]{magneticFK} for a general statement), one gets 
\begin{equation} 
	\label{eq:greenFK}
G_A(x,y) = \int_0^\infty p^A_t(x,y) \d t
= \int_0^\infty p_t(x,y) \mathbb{E}^D_{t,x,y}\big[ \exp( i \alpha \int A_{W_t} \circ \d W_t )  \big] \d t ,
\end{equation}
where $W$ is,  under $\mathbb{P}^D_{t,x,y}$, a Brownian bridge from $x$ to $y$ with duration $t$, killed as it exits the planar domain $D$ (thus  $\mathbb{P}^D_{t,x,y}$ is a \emph{sub}probability measure).  

Let us now formally apply this equality when $A$ is random and distributed as $\mathbb{P}_{YM}$. Fix two points $x,y \in D$, and assume the segment $[x,y]$ is also contained in $D$.\footnote{In the following, we could replace the segment $[x,y]$ with any given smooth path from $x$ to $y$ inside $D$. However, the value of the string observable $\mathcal{S}_{[x,y]}$ does depend on which choice of path we are making.}  We then define the \emph{string observable} $\mathcal{S}_{[x,y]}$ as the (formally defined) random variable  
\[ 
\mathcal{S}_{[x,y]}\coloneqq  \Phi(x) \exp\big(i \alpha \int_{[x,y]} A(z)  \d z\big)  \overline{\Phi(y)},
\]
and $S_{[x,y]}$ its expectation. 
Remark the random variables $\mathcal{S}_{[x,y]}$ are invariant by gauge transformation, as opposed to the random variables $ \Phi(x)  \overline{\Phi(y)} $, which is why these are interesting quantities to look at.\footnote{In the non-abelian and geometrically non-trivial setting, we would consider instead 
	$\langle \Phi(x), \mathcal Hol^A_{y\to x} \Phi(y)\rangle$, which is manifestly gauge invariant. The coupling constant $\alpha$ would be absorbed as a scaling factor in the exponential map. For comparison,  $\langle  \Phi(x) , \Phi(y)\rangle $ would not even make sense in this framework, as $\Phi(x)$ and $\Phi(y)$ would live in different spaces (namely the fibers over resp. $x$ and $y$ of some fiber bundle).
}

Using \eqref{eq:defaYMH} and \eqref{eq:greenFK} (remark also $Z/Z'=\mathbb{E}_{YM}[Z_A]$, which is obtained just by looking at the expectation of $1$), we deduce, formally, 
\begin{align*}
S_{[x,y]}&=\mathbb{E}_{YM}\big[ \frac{Z'Z_A}{Z}  \exp\big(i \alpha \int_x^y A(z) \big) G_A(x,y)  \big]\\
&=\frac{\int_0^\infty p_t(x,y) \mathbb{E}_{YM}\otimes \mathbb{E}_{t,x,y}[Z_A \exp(i \alpha \int_0^{\bar{t}} A_{W_t} \circ \d W_t )]   \d t }{\mathbb{E}_{YM}[Z_A ]}.
\end{align*}

 Using formally Stokes' theorem, we get 
 \[ 
 \int_0^{\bar{t}} A_{W_t} \circ \d W_t=\int_{\mathbb{R}^2} \curl(A)(z) \n_W(z) \d z.
 \] 
Under the measure $\mathbb{P}_{YM}$ (the measure of the abelian Yang--Mills field), $\curl(A)$ is distributed as a Gaussian white noise. Thus, $W$-almost surely, the random variable $\int_{\mathbb{R}^2} \curl(A)(z) \n_W(z) \d z$ should be formally equal to a Gaussian random variable with variance $\int_{\mathbb{R}^2} \n_W(z)^2 \d z=\mathcal{A}_W$ (the non-recentred Amperean area). If we look at the \emph{quenched} measure
	 \[\d \mathbb{P}_{quenched}(\Phi,A)\coloneqq \d \mathbb{P}_A(\Phi)\d \mathbb{P}_{YM}(A)\] rather than the \emph{annealed} measure $\mathbb{P}$ which involves the partition functions $Z_A$, 
 we end up formally with 
\[
\mathbb{E}_{quenched}[\mathcal{S}_{[x,y]}]= 
\int_0^\infty p_t(x,y) \mathbb{E}_{t,x,y}[ \exp( - \alpha^2  \mathcal{A}_W) ]   \d t. \]
The expression for the annealed measure is more convoluted. 
Using the formula \eqref{eq:zeta} for the expression of $Z_A$, we end up instead with 
\begin{align*}
S_{x,y}& =\frac{\int_0^\infty p_t(x,y) \mathbb{E}_{YM}\otimes \mathbb{E}_{t,x,y}\otimes \mathbb{E}_\mathcal{L} [ \exp(i \alpha \sum_{\ell \in \mathcal{L}\cup \{W\}}  \int_0^{\bar{t}_\ell} A_{\ell_t} \circ \d \ell_t )]   \d t }{\mathbb{E}_{YM}\otimes \mathbb{E}_\mathcal{L}[  \exp(i \alpha \sum_{\ell \in \mathcal{L} }  \int_0^{\bar{t}_\ell} A_{\ell_t} \circ \d \ell_t )  ]}\\
&= \int_0^\infty  
\frac{ \mathbb{E}_{t,x,y}\otimes \mathbb{E}_\mathcal{L} [ \exp(- \alpha^2 \sum_{\ell,\ell' \in \mathcal{L}\cup \{W\}}    \mathcal{A}_{\ell,\ell'}  )  ]    }{\mathbb{E}_\mathcal{L}[  \exp(- \alpha^2 \sum_{\ell,\ell' \in \mathcal{L} }   \mathcal{A}_{\ell,\ell'}  )  ]} p_t(x,y) \d t,
\end{align*}
where $\mathcal{A}_{\ell,\ell}$ is the (non-recentred) Amperean area $\mathcal{A}_\ell$, while $\mathcal{A}_{\ell,\ell'}$ is the Amperean area that we previously denoted $\mathcal{B}_{\ell,\ell'}$.\footnote{Using different letters to distinguish the `diagonal' terms $\ell=\ell'$ from the `non-diagonal' ones will prove convenient in Section \ref{sec:X}. Only in this section will we use the notation $\mathcal{A}_{\ell,\ell'}$.} 

Higher order products are obtained in the same way. Consider $(x_i,y_i)_{i\in \{1,\dots, n\}}$ a collection of points in $D$, as well as $(\ell_j)_{j\in \{ 1,\dots ,m \}}$ a collection of smooth loops in $D$. For a smooth loop $\ell$, let also the \emph{loop observable} $\mathcal{S}_\ell\coloneqq  \exp(i \alpha \int A(\ell_t) \d \ell_t )$ (which is also invariant by gauge transformation).\footnote{In the non-abelian setting, 
	$\mathcal{S}_\ell= \operatorname{tr} (\mathcal Hol^A_\ell)$. Although $\mathcal Hol^A_\ell$ depends on the choice of a starting point on $\ell$, its trace does not.}
	
Define
\[ 
S_{quenched}\coloneqq  \mathbb{E}_{quenched}\big[ \prod_{i=1}^n  \mathcal{S}_{[x_i,y_i]} \prod_{j=1}^m \mathcal{S}_{\ell_j} \big]\quad \mbox{and} \quad 
S\coloneqq \mathbb{E}\big[ \prod_{i=1}^n  \mathcal{S}_{[x_i,y_i]} \prod_{j=1}^m \mathcal{S}_{\ell_j} \big].
\]
Let $\Sigma_n$ be the set of permutations of $\{1,\dots, n\}$. For $\sigma \in \Sigma_n$ and $\mathbf{t},\mathbf{x},\mathbf{y}$ a $n$-uplet in $[0,\infty)\times D\times D$, let 
$\mathbb{P}^\sigma_{\mathbf{t}, \mathbf{x}, \mathbf{y}}$ be a measure under which 
$W_1,\dots , W_n$ are independent Brownian bridges, with $W_i$ a Brownian bridge from $x_i$ to $y_{\sigma(i)}$ with duration $t_i$.
Let $I_\sigma$ be the set of cycles in the permutation $\sigma $, and for $i\in I$, say $i=\{i_1,\dots i_k\}$ with $i_{j+1}=i_{\sigma(j)}$, $i_1=i_{\sigma(k)}$, let $W_i$ be the loop obtained as the concatenation 
\[ 
W_i= [y_{i_1},x_{i_1}]\cdot W_{i_1}\cdot [y_{i_2},x_{i_2} ] \cdot W_{i_2}\cdot \dots \cdot  [y_{i_k},x_{i_k} ] \cdot W_{i_k},
\]
and let $W_{I_\sigma}=\{W_i: i\ \in I_\sigma\}$ (see Figure \ref{fig:ex} below for a simple example).
\begin{figure}[h]
	\includegraphics{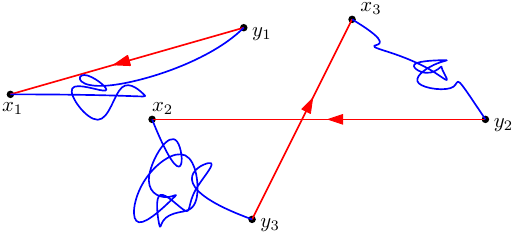}
	\caption{\label{fig:ex} An example with $k=3$ and the permutation $\sigma=(23)$. The set $W_{I_\sigma}$ contains the two depicted loops (the wiggly curves in blue represent Brownian bridges), corresponding to the cycle $(1)$ and the cycle $(2,3)$.}
\end{figure}
 
Then, with a formal computation similar to the one above (and using Isserlis' theorem), we end up with 
\begin{equation}
	S_{quenched}= \sum_{\sigma\in \Sigma_n}  \int_{[0,\infty)^n} 
	\mathbb{E}_{\mathbf{t},\mathbf{x},\mathbf{y}}^\sigma \big[ \exp\big(- \alpha^2  \sum_{\ell,\ell' \in  W_{I_\sigma}\cup \{\ell_j\}_{j\in \{1,\dots,m\}} }   \mathcal{A}_{\ell,\ell'}  \big)  \big]    
	\prod_{i=1}^n p_{t_i}(x_i,y_i) \d t_i,
\end{equation}
and
\begin{equation}
S= \sum_{\sigma\in \Sigma_n}  \int_{[0,\infty)^n} 
\frac{\mathbb{E}_{\mathbf{t},\mathbf{x},\mathbf{y}}^\sigma\otimes  \mathbb{E}_\mathcal{L} [ \exp(- \alpha^2 \sum_{\ell,\ell' \in \mathcal{L}\cup W_{I_\sigma}\cup \{\ell_j\}_{j\in \{1,\dots,m\}} }    \mathcal{A}_{\ell,\ell'}  )  ]    }{\mathbb{E}_\mathcal{L}[  \exp(- \alpha^2 \sum_{\ell,\ell' \in \mathcal{L} }   \mathcal{A}_{\ell,\ell'}  )  ]}
\prod_{i=1}^n p_{t_i}(x_i,y_i) \d t_i.
\end{equation}

Let us remark that for a smooth couple $(\Phi,A)$, the knowledge of the two families $\mathcal{S}_{x,y}$ and the  $\mathcal{S}_{\ell} )$ allows to entirely recover the gauge equivalence class of $(\Phi,A)$. Thus, it is to be expected that the family of moments $S$ that we consider fully characterises the law of $(\Phi,A)$ up to gauge equivalence, 
and this family is entirely described by exponential moments of Amperean areas between Brownian paths and smooth loops. In this paper we mostly consider $\mathcal{A}_{\ell,\ell'}$ when both $\ell$ and $\ell'$ are Brownian. The case when both loops are smooth, treated in section~\ref{sec:smooth}, is standard calculus , and the case when one loop is Brownian and the other is smooth is essentially traditional stochastic calculus and stochastic Green's formula. 

The addition of a \emph{mass} (or \emph{potential}) term $\langle \varphi, m \varphi)$ in the definition of $\mathbb{P}_{YM}$ would introduce, for each Brownian path $W$ with duration $t$ involved, an extra factor $e^{-mt}$, or more generally an extra factor $ \int_0^t m(W_s) \d s$ if $m$ varies in space. In particular, the recentering we used to define $\mathcal{A}_{W}$ do translate into a renormalisation of $\mathbb{P}_{YM}$ by addition of a negative mass term, diverging toward $-\infty$ as the regularisation parameter $\epsilon$ goes to $0$. 

\medskip 

One may wonder how much of this heuristic analysis can be made rigorous. If we fix a gauge and consider a mollification of the Yang--Mills field, i.e. a convolution with a smooth mollifier $\varphi^\epsilon$, the mollified field $A^\epsilon$ can be integrated along Brownian path $W$. With the definition of the mollified winding function $\n_W^\epsilon$ that we will soon give, formally $\n_W^\epsilon$ is the convolution of the winding function $\n_W$ with $\varphi^\epsilon$,
it does hold that 
$
\int_0^{\bar{t}} A^\epsilon_{W_t} \circ \d W_t$ is  a Gaussian random variable with variance $\int_{\mathbb{R}^2}  \n^\epsilon_W(z)^2 \d z$, which is the unnormalised value that we use to construct rigorously 
\[\tilde{\mathcal{A}}_{W}\coloneqq \mathop{\mbox{$L^2$-$\operatorname{lim}$}}_{\substack{\epsilon \to 0 \\ \epsilon\in2^{-\mathbb{N}}}}\ \Big( \int_{\mathbb{R}^2}  \n^\epsilon_W(z)^2 \d z-\mathbb{E}[\int_{\mathbb{R}^2}  \n^\epsilon_W(z)^2 \d z]\Big).\] 
At this mollified level, the formal equalities that we wrote then becomes meaningful and correct\footnote{Here the considered exponentials really are smaller than $1$, so there is no issue with the existence of these exponential moments. The positivity of the partition function follows form that it is the expectation under $\mathbb{P}_{YM}$ of a positive quantity, which in turn follows from \cite{zeta}.}: they follows from standard stochastic analysis. The annealed but mollified measure $\mathbb{P}^\epsilon$, formally $\frac{Z' Z_{A^\epsilon}}{Z} \d \mathbb{P}_{A^\epsilon}(\Phi)\d \mathbb{P}_{YM}(A)$ can then be rigorously constructed when no negative mass term is involved,  and what we prove in this paper essentially amount to the convergence of the moments of this measure (when the appropriate negative mass term is added, and except that we do not prove the finiteness of some exponential moments of the Amperean area, which corresponds to polynomial moments of the field).

\medskip

Let us summarise how 
	Symanzik's polymer representation roughly translates terms in the action defining a field into terms associated to Brownian trajectories, as long as we stay ignorant of the renormalisation issues. 
\begin{table}[H] 
	\begin{tabularx}{\linewidth}{Y|Y}
	Action $ S(\Phi,A)$ &  {Weight for Brownian curves}\\  
	\hline 
	 	& \\
	{Mass or potential term }  $\langle \Phi, m \Phi \rangle$ & { Occupation measure $\mu_W$,} $\int_D m(z) \mu_W(\d z)= \int m(W_s) \d s$  \\
 	& \\
	{Interaction with an \emph{external} magnetic potential } $\|\nabla \Phi + i \alpha A \Phi\|^2$   &  {Winding function $\n_W$,}    $\int_D \n_W(z)  \curl(A)  \d z= \int A(W_s) \circ \mathrm{d} W_s$  \\
 	& \\
	{Renormalised self-interaction potential } $\langle \Phi, \lambda \Phi\rangle_{L^4 }^2$  &  {Renormalised  self-intersection field $I_W$,} $\int_D \lambda (z) I_W(\d z)$  \\
 	& \\
	{Interaction with renormalised \emph{internal} magnetic potential } $\|\nabla \Phi+ i \alpha A \Phi\|^2 +\|\curl A\|^2$ (integration now also over $A$) &   {Squared winding function $\n_W^2$,}    $\int_D \n_W^2(z)  \alpha^2  \d z$ .
\end{tabularx}
\end{table}

Remark both the mass term and the interaction with an external magnetic potential are quadratic (in $\Phi$), but the self-interaction term or the term of interaction with an internal magnetic potential  are quartic ($\Phi$ or in $(\Phi,A)$, respectively).
It is the fourth link (which I am not aware is documented) which makes the study of the Amperean area relevant, from the point of view of quantum field theory.

\begin{remark} 
In \cite{AlbeverioKusuoka}, the authors study a related but different quantity. The starting point is also that of the holonomy of an abelian Yang--Mills connection along a Brownian path. However, where we first average over the connection, they first average over the Brownian path: 
the quantity we define rigorously is a renomalised random variable, associated with the Brownian trajectory, corresponding to $\mathbb{E}_A[ \exp( i \int A_W \circ \d W ) ]$. What they define, instead, is a random variable, associated to a random connection, corresponding to 
$\mathbb{E}_{t,x,y}[ \exp( i \int A_W \circ \d W ) ]$. Rather than a mollification approach, they use a cut-off in the Fourier space. The main advantages of our approach is first that we can give an explicit formula for the limit, secondly that both this limiting formula and the proof rely almost entirely on classical probability theory rather than higher-level technics, and finally that the renormalisation term is an explicit mass term. On the other hand, there is two main drawbacks. The first is that the existence of all the (negative) exponential moments seems more complicate to show. The main reason is that the Fourier decomposition, which we do not rely on, is of course very practical to decouple the contributions from different scales. The second drawback is that, as opposed to 
$\mathbb{E}_{t,x,y}[ \exp( i \int A_W \circ \d W ) ]$ which is a genuine observable of the physical field $A$ for which the Brownian path is but a mathematical tool, the quantity   
$\mathbb{E}_A[ \exp( i \int A_W \circ \d W ) ]=\exp(-\frac{\alpha^2}{2} \mathcal{A}_W)$ does not capture any information about $A$ itself but only about its distribution, and instead depends on the Brownian path which has no physical relevance. This would turn into an advantage in situations unrelated to the aYMH field where the Brownian particle becomes physically relevant e.g. for electron displacement in inhomogeneous magnetic fields.
\end{remark}

\section{Regularisation procedure and main results}
\label{sec:main}
Let $G(z)=\log|z|/(2\pi)$.
We fix a mollifier $\varphi:\mathbb{R}^2\to \mathbb{R}_+=[0,\infty)$ with integral $1$. We assume it is smooth (i.e. $\mathcal{C}^\infty$), compactly supported, and rotationally invariant. We then set $\varphi^\epsilon$ the rescaled function $x\mapsto \epsilon^{-2}\varphi(\epsilon^{-1} x)$. We define a vector field on $\mathbb{R}\setminus\{0\}$ by
$\theta_z= \frac{(-z^2,z^1)}{2 \pi |z|^2}$, where $z=(z^1,z^2)$. Remark for any $z\in \mathbb{R}^2\setminus \operatorname{Range}(W)$, it holds that almost surely $\n_{W}(z) = \int_0^{\bar{T}} \theta_{W_t-z} \circ \d W_t$.\footnote{This can be proved for example by taking dyadic piecewise linear approximations $W_n$ of $W$. For $z\notin \Range(W)$ and for all $n$ large enough, it holds $\n_{W}(z)=\n_{W_n}(z)= \int_0^{\bar{T}} \theta_{W_n(t)-z} \d W_n(t)$. The first equality follows from homotopy invariance of the winding number, since $W_n$ is homotopy equivalent to $W$ in $\mathbb{R}^2\setminus \{z\}$ as soon as $\|W_n-W\|_\infty<d(z, \Range(W))$. The second equality follows from the residue theorem for piecewise linear loops. 
	 We conclude by using the fact that $\int_0^{\bar{T}} \theta_{W_n(t)-z} \d W_n(t)$ converges in $L^2(\mathbb{P})$, as $n\to \infty$, toward the Stratonovich integral $\int_0^{\bar{T}} \theta_{W_t-z} \circ \d W_t$. } We define the mollified winding function $\n^\epsilon_W$ by 
\begin{equation}
	\label{eq:def:nepsilon}
\n^\epsilon_W=  \int_0^{\bar{T}} (\varphi^\epsilon\conv \theta)_{W_t-z} \circ \d W_t,
\end{equation}
where $\conv$ is the convolution on $\mathbb{R}^2$. 
If it was possible to change the order between the convolution and the stochastic integral, we would get $``\n^\epsilon_W= \varphi^\epsilon\conv \n_W"$, hence the interpretation of $\n^\epsilon$ as a mollified winding function. Yet, as we explained already, $\n_W$ is not integrable, so the convolution $``\varphi^\epsilon\conv \n_W"$ is in fact ill-defined (but $\n^\epsilon_W$ is a well-defined and measurable function, which up to modification is a smooth function---see Remark \ref{rem:RP} below). 

	Let $W,W'$ be two independent planar Brownian motions with durations $T,T'>0$ and with starting points $a,a'\in \mathbb{R}^2$. Our main results are the following.
\begin{theorem} 
	\label{th:cv}
	For all $T,T'>0$ and $a,a'\in \mathbb{R}^2$, the random variables 
	\[
	\mathcal{B}^\epsilon_{W,W'}\ \coloneqq\int_{\mathbb{R}^2} \n^\epsilon_{W }(z) \n^\epsilon_{W' }(z) \d z \qquad \mbox{and} 
	\qquad 
	\tilde{\mathcal{A}}_{W}^\epsilon \coloneqq \int_{\mathbb{R}^2} \n^\epsilon_{W }(z)^2 \d z  -\mathbb{E}[\int_{\mathbb{R}^2} \n^\epsilon_{W }(z)^2 \d z ]
	\]
	both converge in $L^2$ in the limit $\epsilon\to 0$. 
\end{theorem} 
The limits are denoted $\mathcal{B}_{W,W'}$ and $\tilde{\mathcal{A}}_{W}$. 
\begin{theorem}
	\label{th:Amperean} 
	Almost surely, 	
	\begin{equation}
	\label{eq:Amperean}
	\mathcal{B}_{W,W'}= -\int_0^{\bar{T}}\int_0^{\bar{T}'} G(W_s-W'_t) \circ \d W'_t \circ \d W_s 
	=-\int_0^{\bar{T}}\int_0^{\bar{T}'} G(W_s-W'_t) \d W'_t  \d W_s -
	I_{W,W'}/4,
	\end{equation}
	where $I_{W,W'}$ for the intersection local time between $W$ and $W'$. 
\end{theorem} 
\begin{theorem} 
	The counterterm $\mathbb{E}[\int_{\mathbb{R}^2} \n^\epsilon_{W }(z)^2 \d z  ]$ is given, asymptotically as $\epsilon \to 0$, by 
	\[ 
	\mathbb{E}\big[\int_{\mathbb{R}^2} \n^\epsilon_{W }(z)^2 \d z\big]=\frac{T}{2\pi}\log(\epsilon^{-1})+\frac{T\log{T}}{4\pi}+C_\varphi T+O(\epsilon^\frac12), 
	\]
	where $C_\varphi$ is a constant which depends on the mollifier $\varphi$ but on no other parameter. 
\end{theorem} 
\begin{theorem} 	
	For $|\beta|< \frac{8\pi}{48\sqrt{2}+\pi} \simeq 0.35$, the random variable
	$\exp (\frac{\beta}{\sqrt{TT'}} \mathcal{B}^\epsilon_{W,W'} )$ converges in $L^1$ toward $\exp (\frac{\beta}{\sqrt{TT'}} \mathcal{B}_{W,W'} )$, and the convergence is uniform in $a,a'\in \mathbb{R}^2$ and $T,T'\in(\delta,\infty)$, for all $\delta>0$. 
	The limits do not depend on the choice of the mollifier $\varphi$. 
	
	There also exists $\beta_c>0$ such that for  $|\beta|<\beta_c$, 
	$\mathbb{E}[ \exp (\frac{\beta}{T} \tilde{\mathcal{A}}_{W} )]<\infty$. 
\end{theorem} 
\begin{remark} 
	The first term in the expansion of the counterterm, 
	$	\mathbb{E}[\int_{\mathbb{R}^2} \n^\epsilon_{W }(z)^2 \d z]\sim \frac{T}{2\pi}\log(\epsilon^{-1})$, 
	is new but expected, as it matches with the result obtained in \cite{Werner3} with a different regularisation method. The terms up to order $o(1)$ are physically relevant, as they correspond to masses. 
\end{remark} 
\begin{remark}
	\label{rem:conformal}
	In the positive probability event that $W$ and $W'$ are contained inside a given simply connected domain $D$, we can replace $G(W_s-W_{t'})$ in Theorem \ref{th:Amperean} with $G^D(W_s,W'_t)$, where $G^D$ is the Green function of the domain $D$ (with Dirichlet boundary condition). This is because for any given $a\in D$, the difference $G(a-\cdot)-G^D(a,\cdot)$ is harmonic on $D$, so its Stratonovich integral along a loop contained inside $D$ is equal to $0$. 
	\end{remark} 
\begin{remark}
	\label{rem:`miraculous`cancellation} 
	Although we do not prove it here, for a smooth enough function $f$, it is be possible to similarly construct the random variables 
	$\mathcal{B}_{W,W'}(f)=`` \int f(z) \n_{W}(z) \n_{W'}(z) \d z"$ and it should also be possible to define  $\tilde{\mathcal{A}}_{W}(f) =``\int f(z) (\n_{W}(z)^2-\mathbb{E}[\n_W(z)^2]) \d z"$. Unlike the case $f=1$, the counter-term that needs to be subtracted to define $\tilde{\mathcal{A}}_{W}(f)$ is then random, given at the main order by 
	\begin{equation}\label{eq:counter}  \frac{\log(\epsilon^{-1})}{2\pi} \int f(W_t) \d t.
	\end{equation} 
	Remark this is half of the counter-term needed to define $I_W(f)$ the self-intersection local time of $W$ with respect to $f$: formally, 
	\begin{equation}
		\label{eq:defI}
	I_W(f)=\int_{\mathbb{R}^2}\int_{0\leq s \leq t \leq T} f(z) \delta_{W_s=W_t}\delta_{W_s=z} \d s \d t \d z,
	\end{equation} 
 see e.g. \cite[Theorem 1.1]{Dynkin}. It thus seems indeed that $\tilde{\mathcal{A}}_{W}(f)-I_W(f)/2$ can be defined without renormalisation term (i.e. the counterterms cancel each other up to $o(1)$). The term $I_W(f)/2$ can be seen as an It\^o-to-Stratonovich correction term for $\tilde{\mathcal{A}}_{W}(f)$ (compare with Theorem \ref{th:Amperean} above and Conjecture~\ref{conj:formula} below.\footnote{While the intersection local time $I_{W,W'}$ is defined from the rectangle $[0,T]\times [ 0,T']$, the self-intersection local time $I_W$ is usually normalised so that it formally matches the expression \eqref{eq:defI}, where the time integrals is performed over the triangle $0\!<\!s\!<\!t\!<\!T$. This explains why the former appears with a prefactor $1/4$ when the latter appears with a prefactor $1/2$.
}
\end{remark} 	

\begin{remark} 
		From geometric consideration, it is natural to consider, in Remark \ref{rem:`miraculous`cancellation} above, the case $f={\det(g)}$ for some underlying Riemannian metric $g$, and $W$ a Brownian motion with respect to the same metric $g$. Then, the above integral $\int f(W_t)\d t$ is the $g$-quadratic variation of $W$, so that the counterterm~\eqref{eq:counter} matches with the former expression $T\log(\epsilon^{-1})/2\pi$. In the case there is a global conformal transformation $\phi: D\to D'$ which maps the metric $g$ to the euclidean metric,  Remark \ref{rem:conformal} above suggests the definition $\tilde{\mathcal{A}}_W(f)\coloneqq \tilde{\mathcal{A}}_{\phi(W)}$. However the convergence of the mollified  $\tilde{\mathcal{A}}^\epsilon_W(f)$ toward $\tilde{\mathcal{A}}_W(f)$ would not be immediate as the map $\phi$ would distort the mollifier (although such a convergence would certainly be true, at least up to addition of a local time term corresponding to the $O(1)$ terms in the asymptotic expansion of the counterterm). When such conformal transformations exist locally but not globally, it should still be possible to define $\tilde{\mathcal{A}}_W(f)$ in this way, by appropriately chopping the Brownian trajectory into small pieces fitting the domain of definition of these conformal transformations. 
\end{remark}

\begin{remark} From the formula for $\mathcal{B}_{W,W'}$, it is possible to infer that the renormalised Amperean area $\tilde{\mathcal{A}}_W$ can also be expressed by stochastic integrals, given by the following conjecture.
\begin{conjecture}
	\label{conj:formula}
The renormalised Amperean area $\tilde{\mathcal{A}}_W$ can be expressed as 
	\begin{align*} 
	\mathcal{A}_W&=-2\int_0^{\bar{T}} \int_0^s G(W_s-W_t) \d \overleftarrow{W}_t \d W_s  
	-I_W/2\\
	&=-2 \int_0^{\bar{T}} \int_0^s G(W_s-W_t) \d W_t \d W_s  +I_W/2,
	\end{align*}
where $\d \overleftarrow{W}_t $ stands for a Backward integral and $I_W$ is the recentred self-intersection local time. 
\end{conjecture}	
Note the counterterm needed to define $\mathcal{A}_W$ would be absorbed as the counter-term needed to define $I_W$. It is not clear to the author that the corresponding Stratonovich integral 
\[
 \int_0^{\bar{T}} \int_0^s G(W_s-W_t) \circ \d W_t \d W_s  
\]
is well-defined (although we can for example define it as the average between the backward and forward integrals). 

\end{remark}

To conclude this section, let us remark that by analogy with the case of the Lévy area, there is possibly a different way to construct the Amperean area, which would not rely on mollification but on the use of symmetric cut-off on the image set $\mathbb{Z}$ where the winding functions take their values. Although I believe such a construction would be theoretically interesting, I do not know how it could be rigorously related to the aYMH field.  One advantage of such a normalisation is that it is more intrinsic than mollifications. It would allow to compute 
the Amperean areas with the data of the curves and the Lebesgue measure on the plane, without the detail of the underlying Riemmanian (indeed Euclidean) metrics. Conformal covariance, for example, would follow directly.

\begin{conjecture}
	For relative integers $k,j$, let $A_k$ and $B_{k,j}$ be the measure of the sets 
	\[ 
	\mathcal{A}_k(W)\coloneqq \{z:  \n_W(z)=k\}, \qquad \mathcal{B}_{k,j}\coloneqq 	\mathcal{A}_k(W) \cap 	\mathcal{A}_j(W').
	\] 
	Then , the sums $\sum_{k\in \mathbb{Z}} k^2 |A_{k}-  \mathbb{E}[ A_k]|$ and 
	$ \sum_{(k, j) \in \mathbb{Z}^2 }  |kj| B_{k,j}$ are both infinite, but 
	the random variables $\mathcal{A}_{W}$ and $\mathcal{B}_{W,W'}$ are given by 
	\[
	\mathcal{B}_{W,W'}= \lim_{k_0,j_0\to \infty } \sum_{k=-k_0}^{k_0} \sum_{j=-j_0}^{j_0} kj B_{k,j}   , 
	\quad
	\mathcal{A}_{W} = 
	\lim_{k_0\to \infty } \sum_{k=-k_0}^{k_0}  k^2 (A_{k}-  \mathbb{E}[ A_k]),
	\]
	where the convergence hold both in the almost sure sense and in the $L^p$ sense for all $p\in[1,\infty)$ (and in fact with convergence of some exponential moments).
\end{conjecture}


\section{Notations}
\label{sec:not}
The Euclidean norm in $\mathbb{R}^2$ is written $|\cdot |$. The ball centred at $0$ with radius $r$ is written $B_r$, and the ball centred at $z$ with radius $r$ is written $B_r(z)$.

We fix a mollifier $\varphi:\mathbb{R}^2\to \mathbb{R}_+=[0,\infty)$ with integral $1$. We assume it is smooth, compactly supported, and rotationally invariant. We then set $\varphi^\epsilon$ the rescaled function $x\mapsto \epsilon^{-2}\varphi(\epsilon^{-1} x)$, the integral of which is also equal to $1$. We let $K_\varphi>0$ be the minimal value such that the open support of $\varphi$ is contained inside $B_{K_\varphi}$. 

We often write the evaluation of $f$ at $z\in \mathbb{R}^2$ by $f_z=f(z)$. Coordinates are usually written as superscripts, e.g. $x=(x^1,x^2)\in \mathbb{R}^2$. 
The derivative of a function $f$ indexed by a subset of $\mathbb{R}$ is always written $\dot{f}$, while $f'$ usually just refers to a different function unrelated to $f$.


The gradient of a function $f\in \mathcal{C}^1(D\subseteq \mathbb{R}^2,\mathbb{R})$ is the vector field 
$\grad f= (\partial_1 f, \partial_2 f)\in \mathcal{C}^0(D,\mathbb{R}^2)$ defined on $D$. The curl of a vector field  $V\in \mathcal{C}^1(D\subseteq \mathbb{R}^2,\mathbb{R}^2)$ is the function 
$\curl V=\partial_1 f^2 -\partial_2 f^1\in   \mathcal{C}^0(D,\mathbb{R}^2)$, and its divergence is the function $\div V=\partial_1 f^1 +\partial_2 f^2\in   \mathcal{C}^0(D,\mathbb{R}^2)$. For a vector field $V=(V_1,V_2)$, we write $\star V$ for the vector field rotated by a quarter of turn, $\star V=(-V^2,V^1)$ (one can think of $\star$ either as the multiplication by the complex number $i$, or as a Hodge dual). 
These operators commute with the translation operators and with the operators of convolution (written~$\conv$) with a given function. It holds $\star \star V=-V$, $\curl(\star V)=\div(V)$, $\div(\star V)=-\curl(V)$, $\div\circ \grad=\Delta$, $\curl\circ \grad=0$, and all these identities also holds in the distributional sense. 

We define the Green kernel $G:  \mathbb{R}^2\setminus \{0\}\to \mathbb{R}$, and the vector field $E=\grad G $, $\theta= \star  E$, by
\[ 
G_z=\frac{\log|z|}{2\pi}, \qquad 
E_z=\frac{(z^1,z^2)}{2\pi |z|^2}, \quad 
\theta_z=\star E_z= \frac{(-z^2,z^1)}{2 \pi |z|^2}. 
\]
Recall $\Delta G=\delta_0$ is the Dirac mass at $0$, so that \[\div(E)=\delta_0, \quad \curl(\theta)=\delta_0.\] Furthermore, it follows from the orthogonality in $L^2(\mathbb{R}^2,\mathbb{R}^2)$ between the subspaces \[\{ V: \exists f \text{ s.t. }V=\grad(f) \}
\text{ \quad and \quad } \{ V: \exists f\text{ s.t. }V=\star  \grad(f) \}\] that 
\[\curl(E)=\div(\theta)=0,\qquad \text{not only pointwise but also in the distributional sense.}\]
We write $G^\epsilon, E^\epsilon, \theta^\epsilon$ for the convolutions of $G,E,\theta\in L^1+L^\infty$ with the mollifier $\varphi^\epsilon\in L^1\cap L^\infty$.

We write $p_t$ for the whole-plane heat kernel, $p_t(x,y)=\exp(-|x-y|^2/(2t))/(2\pi t)$. For $r\geq 0$, we use the shortcut notation $p_t(r)$ for the common value $p_t(x,y)$ of any $x,y: |x-y|=r$, i.e. $p_t(r)=\exp(-r^2/(2t))/(2\pi t)$.

For $T>0$ and $a \in \mathbb{R}^2$, we write $\mathbb{P}_{T,a}$ for a probability measure under which $W$ is a Brownian motion started from $a$ and with duration $T$. For $T'>0$ and $a' \in \mathbb{R}^2$,  we write $\mathbb{P}_{T,T'a,a'}$ for a probability measure under which a second Brownian motion $W'$ is defined, independent from $W$, started from $a'$ and with duration $T'$.
These subscripts are omitted from the notation except when we want to stress the dependency on these variables.

We continue to write $W$ (resp. $W'$) for the loop obtained by concatenation of $W$ with the straight line segment between its endpoints. This loop is parameterised by $[0,\bar{T}]$ for some $\bar{T}>T$, in such a way that this parametrisation matches with the initial process on $[0,T]$, and has constant speed on $[T,\bar{T}]$. The specific value chosen for $\bar{T}$ is irrelevant in the following, because the integrals we consider are invariant by increasing reparametrisation. 

For an oriented segment $[a,b]$ in the plane and a smooth vector field $V$, we write either $\int_a^b V$ or $\int_a^b V(u) \d u$ for the integral of $V$ along $[a,b]$: if $\gamma:[0,1]\to [a,b]$ is a smooth bijection with $\gamma_0=a$, then $\int_a^b V=\int_a^b V(u)\d u= \int_0^1 \langle V_{\gamma_s} ,\dot{\gamma}_s\rangle \d s$.
Stratonovich and It\^o integrals, either in one or two dimensions, are written respectively with $\circ \d W_t$ and $\d W_t$. In particular, for a vector field~$V$,  
\[ 
\int_0^{\bar{T}} V_{W_t} \circ \d W_t= \int_0^{\bar{T}} V^1_{W_t} \circ \d W^1_t+
\int_0^{\bar{T}} V^2_{W_t} \circ \d W^2_t=\int_0^{T} V_{W_t} \circ \d W_t + \int_{W_T}^{W_0} V  ,
\]
and the same equalities hold for It\^o integrals. Be careful that $\int_T^{\bar{T}} V_{W_t} \d W_t$ may look like a stochastic integral when it is actually a Riemann integral, as $W_{[T,\bar{T}]}$ is a parametrisation of the straight line segment from $W_T$ to $W_0$. 

Recall the mollified index function of $W$ is defined as the Stratonovich integral 
\begin{equation}
\label{def:nepsilon}
\n^\epsilon_W(z)\coloneqq \int_0^{\bar{T}} \theta^{\epsilon}_{W_t-z} \circ \d W_t, \qquad \theta^\epsilon=\theta \conv \varphi^\epsilon.
\end{equation}
For $\epsilon=0$, and $z\in \mathbb{R}^2\setminus \Range(W)$, it does hold $\n_W(z)= \int_0^{\bar{T}} \theta_{W_t-z} \circ \d W_t$.

When $W$ is replaced with a smooth loop $\gamma$, it can be checked that this definition \eqref{def:nepsilon} gives indeed   $\n^\epsilon_\gamma(z)=\varphi^\epsilon \conv \n_\gamma(z)$, where $\conv$ is the convolution operator, and also for $\epsilon=0$, $\n_\gamma(z)=
\int \theta_{\gamma_t-z} \d \gamma_t$. 

Note that, since the Stratonovich integral satisfies $\int_0^T \grad f (W) \circ \mathrm{d} W=f(W_T)-f(W_0)$ for any smooth $f$, hence $\int_0^{\bar{T}} \grad f (W) \circ \mathrm{d} W=f(W_{\bar{T}})-f(W_0)=0$,  
we could have replace the vector field $\theta^{\epsilon}$ in the definition of $\n^\epsilon_W(z)$ with any other vector field $\tilde{\theta}^{\epsilon}$ so long as $\curl(\tilde{\theta}^{\epsilon}    )=\varphi^{\epsilon}$.

 Since $\theta$ is harmonic on $\mathbb{R}^2\setminus\{0\}$ and $\varphi^\epsilon$ is supported on $B_{\epsilon K_\varphi}$, the vector field $\theta^\epsilon$ is harmonic on $\mathbb{R}^2\setminus B_{\epsilon K_\varphi}$,  and it follows that $\n_W^\epsilon(z)$ is compactly supported on $B_{\|W\|_\infty+\epsilon K_\varphi }$ (a result that is to be expected from the interpretation of $\n^\epsilon_W=``\n_W\conv \phi^\epsilon"$). By derivation below stochastic integral (see Remark \ref{rem:RP} below), for any given $\epsilon>0$, the smoothness of $\theta^\epsilon$ ensures that $\n_W^\epsilon$ is almost surely a smooth function. Thus, 
for $\epsilon,\epsilon'>0$, \begin{equation} 
	\label{eq:defB}
	\mathcal{B}^{\epsilon,\epsilon'}_{W,W'}=\int_{\mathbb{R}^2} \n_W^\epsilon(z) \n_{W'}^{\epsilon'}(z)\d z \quad \mbox{and} \quad \mathcal{A}^{\epsilon}_{W}=\int_{\mathbb{R}^2} \n_W^\epsilon(z)^2 \d z
	\end{equation} 
 are both well-defined.
\begin{remark} 
	Nothing prevents from mollifying $\n_W$ and $\n_{W'}$ with two different mollifiers. Here we consider possibly $\epsilon\neq \epsilon'$ not with the intent to be as general as possible but simply because some expressions are simpler to keep track of when we distinguish the two mollifiers. 
\end{remark} 
\begin{remark}
	\label{rem:RP}
	There is in fact two very similar ways to interpret the definition  \eqref{def:nepsilon} of $\n^\epsilon_W$, which are \emph{versions} of each other (i.e. for any given $z$, these two definitions of $\n^\epsilon_W(z)$ are almost surely equal, but they may not be almost surely equal as functions of $z$). 
	
	The first, which we shall call $\n^\epsilon_{W,\mathrm{stoc}}$ for the sake of this remark, is as a Stratonovich integral, by which we mean as the sum of the corresponding It\^o integral and Stratonovich-to-It\^o correction term. Remark that $\omega$-almost surely, this is defined simultaneously for all $z$, as the It\^o map is more than just a collection of random variables.%
	\footnote{ In fact, this is already a \emph{version} of the initial definition of Stratonovich. I do not know whether it is the case that the initial definition of Stratonovich makes sense simultaneously for all $z$, and whether it is then jointly measurable in $(z,\omega)$.} 	
	Furthermore, it follows from \cite[Theorem 2.2, (1)]{StochasticFubini} that the function $(\omega,z)\mapsto \n^\epsilon_{W(\omega),\mathrm{stoc}}(z)$ is jointly measurable. 
	
	The second interpretation, which is in fact the one we use and which we denote $\n^\epsilon_W$, is as a \emph{rough path} Stratonovich integral. That is, we endow the path $W$ with the data of the almost surely continuous versions of $t\mapsto \int_0^t W^i_s\circ \d W^j_s$ ($i,j\in \{1,2\}$), and we use the enriched path $(X(\omega),\int_0^\cdot W^i_s\circ \d W^j_s(\omega) )$ to construct, pathwise, other integrals along $X(\omega)$ using rough path theory.
	It follows from \cite[Theorem 4.4, Equation (4.16)]{FrizHairer} that the map $z\mapsto  \n^\epsilon_W(z)$ is not only almost surely well-defined simultaneously for all $z$, but also continuous in $z$, and it follows in particular that the map $(\omega, z)\mapsto  \n^\epsilon_{W(\omega)}(z)$ is also jointly measurable. 
	
	We have the following relations between these two random functions $\n^\epsilon_{W}$ and $\n^\epsilon_{W,\mathrm{stoc}}$ : they are indeed version of each other, which follows for example from \cite[Corollary~5.2]{FrizHairer}. It follows in particular that 
	$\n^\epsilon_W$ is the continuous version of $\n^\epsilon_{W,\mathrm{stoc}}$ (uniquely defined up to almost sure equality of functions). Furthermore, since both versions are jointly measurable, it follows from Fubini theorem that almost surely, for all integrable function $f:\mathbb{R}^2\to \mathbb{R}$, 
	\[ 	
	\int f(z) \n_W^\epsilon(z)\d z=\int f(z) \n^\epsilon_{W,\mathrm{stoc}}(z)\d z \qquad \text{and} 
	\int f(z) \n_W^\epsilon(z)^2\d z=\int f(z) \n_{W,\mathrm{stoc}}^\epsilon(z)^2\d z.
	\]
	We will use this in the proof of Lemma~\ref{le:green1}, as it will allow us to rely on the stochastic Fubini theorem. 
	
	If one discards this subtlety of versions and measurability, and treats $ \n_W^\epsilon(z)$ as a Stratonovich integral, no 
	 knowledge of rough path theory is necessary to understand the paper. Being careful about this issue, one must keep in mind that every time we write a stochastic integral which depend on the parameter $z$, it is to be understood either in the rough path sense or as the version of the stochastic integral that is continuous in $z$ (which are equivalent).
	
	As for the fact that $z\mapsto \n_\epsilon(z)$ is not only continuous but in fact smooth in $z$, it follows from a lemma of derivation below rough path integrals (Lemma~\ref{le:derivBelowRP}) which we prove using elementary rough path technics. For the pendant of this result in the framework of stochastic integration (and outside the framework of rough paths), see \cite[Theorem 2.2]{Hutton}.  
\end{remark} 


\section{Amperean area of smooth loops}
\label{sec:smooth}
In this section, we derive the formula equivalent to \eqref{eq:Amperean} for the case of smooth loops. This formula is not used later in the paper, but it provides us with a hint at the formula~\eqref{eq:Amperean}, and with a general strategy for the proof in the Brownian setting (although some of the steps have to be circumvented). 

\begin{theorem}
	\label{th:StokesLisse}
	Let $\gamma\in \mathcal{C}^1([0,T], \mathbb{R}^2)$ and $\gamma' \in \mathcal{C}^1([0,T'], \mathbb{R}^2)$ be two loops.  Then,
	\[ 
	\int_{\mathbb{R}^2} \n_{\gamma}(z)  \n_{\gamma'}(z) \d z =-  \int_0^{T'} \int_0^T G(\gamma_t-\gamma'_{s})  \langle \dot{\gamma}_t,  \dot{\gamma}'_{s} \rangle \d t \d {s}.
	\] 
\end{theorem}
\begin{proof} 
	By Banchoff-Pohl inequality \cite{Banchoff}, it holds that the index functions $\n_\gamma$ and $\n_{\gamma'}$ are both in $L^2(\mathbb{R}^2)$. Since $\varphi$ is a mollifier, it follows that 
	\[ 
	\int_{\mathbb{R}^2} \n_{\gamma}(z)  \n_{\gamma'}(z)  \d z=\lim_{\epsilon \to 0} 
	\int_{(\mathbb{R}^2)^2} \n_{\gamma}(z)  \n_{\gamma'}(w) \varphi^\epsilon(w-z) \d z \d w.
	\]
	For an arbitrary $z\in \mathbb{R}^2$ and since $\curl \theta^{\epsilon}=\varphi^\epsilon$, 
	 Green's formula applied to the loop $\gamma'$ and the smooth vector field $w\mapsto \theta^{\epsilon}_{w-z}$ gives 
	\[
	\int_{\mathbb{R}^2} \n_{\gamma'}(w) \varphi^\epsilon(w-z) \d w= \sum_{i\in\{1,2\}} \int_0^{T'}    (\theta^{\epsilon}_{\gamma'_{s}-z})^i (\dot{\gamma}'_{s})^i \d s.
	\]
	Using Fubini-Tonelli theorem (recall $\n_{\gamma}\in L^1_c(\mathbb{R}^2)$, and $\theta^\epsilon$ is bounded), we deduce 
	\begin{align*}
	\int_{(\mathbb{R}^2)^2} \n_{\gamma}(z)  \n_{\gamma'}(w) \varphi^\epsilon(w-z) \d z \d w
&=		\int_{\mathbb{R}^2} \n_{\gamma}(z) \sum_{i\in\{1,2\}} \int_0^{T'}    (\theta^{\epsilon}_{\gamma'_{s}-z})^i (\dot{\gamma}'_{s})^i \d s \d z\\
&=	 \sum_{i\in\{1,2\}} \int_0^{T'}  (\dot{\gamma}'_{s})^i  \int_{\mathbb{R}^2} \n_{\gamma}(z)  (\theta^{\epsilon}_{\gamma'_{s}-z})^i  \d z \d s , 
	\end{align*}
	hence	
	\begin{align*}
	\int_{\mathbb{R}^2} \n_{\gamma}(z)  \n_{\gamma'}(z)  \d z
	&=\phantom{-} \lim_{\epsilon \to 0}  \sum_{i\in\{1,2\}}
	\int_0^{T'}  (\dot{\gamma}'_{s})^i 	\int_{\mathbb{R}^2} \n_{\gamma}(z) (\theta^{\epsilon}_{\gamma'_{s}-z})^i  \d z \d s.
	\end{align*}
	Let $G_1$ and $G_2$ be the vector field $G_1(w)= ( G_w , 0 )$ and 
	$G_2(w)= (0, G_w  )$ (recall $G$ is the Green function), and let $G_i^{\epsilon}=\varphi^\epsilon \conv G_i$ (which is a well-defined convolution). Then, for $i\in \{1,2\}$, the curl  of $G_i^{\epsilon}$ is equal to the function $(\theta^{\epsilon})^i$, thus $\curl ( G_i^{\epsilon} (\ \cdot \ -\gamma'_s))(z) =(\theta^{\epsilon})^i(z-\gamma'_s )=-(\theta^{\epsilon})^i(\gamma'_s-z )$.
	For an arbitrary $s\in [0,T']$, and $i\in\{1,2\}$, Green's formula applied to the loop $\gamma$ and the vector field 
	$z\mapsto -G^{\epsilon}_i(z-\gamma'_{s})$ gives 
	\[ 
	\int_{\mathbb{R}^2} \n_{\gamma}(z) (\theta^{\epsilon}_{\gamma'_s-z})^i \d z= -\sum_{j\in\{1,2\}} \int_0^{T}     (G_{i}^{\epsilon}(\gamma_t-\gamma'_s))^j     (\dot{\gamma}_{t})^j \d t
	=-\int_0^{T}     G^{\epsilon}(\gamma_t-\gamma'_s)     (\dot{\gamma}_{t})^i \d t.
	\]
	We end up with 
	\begin{align*}
		\int_{\mathbb{R}^2} \n_{\gamma}(z)  \n_{\gamma'}(z)  \d z
		&=-   \sum_{i\in\{1,2\}} \lim_{\epsilon \to 0} 
		\int_0^{T'}  (\dot{\gamma}'_{s})^i (\dot{\gamma}_{t})^i G^{\epsilon}(\gamma_t-\gamma'_{s})     \d t \d s\\
		&= -  \lim_{\epsilon \to 0} 
		\int_0^{T'}  \langle \dot{\gamma}'_{s} , \dot{\gamma}_{t}\rangle (\varphi^\epsilon \conv G)( \gamma_t-\gamma'_s  )       \d t \d s.
	\end{align*}
	Let $p\in(1,\infty)$, which can be chosen arbitrarily. Since $G\in L^{p}(B_{\|\gamma\|_\infty+\|\gamma'\|_\infty+K_\varphi})$, the convergence of $\varphi^\epsilon \conv G$ toward $G$ holds in $L^{p}(B_{\|\gamma\|_\infty+\|\gamma'\|_\infty})$, which ensures that $\varphi^\epsilon \conv G$ is dominated by a function in $L^{1}(B_{\|\gamma\|_\infty+\|\gamma'\|_\infty})$. We can thus apply the dominated convergence theorem and conclude indeed that 
	\[ 		\int_{\mathbb{R}^2} \n_{\gamma}(z)  \n_{\gamma'}(z)  \d z 
	=	-\int_0^{T'}  \langle \dot{\gamma}'_{s} , \dot{\gamma}_{t}\rangle  G( \gamma_t-\gamma'_s  )       \d t \d s.
	\]
\end{proof} 
\begin{remark}
Unlike the Stokes' theorem which can be extended to loops in $\mathcal{C}^\alpha$ for any $\alpha>1/2$, Theorem \ref{th:StokesLisse} cannot be extended to the case of loops in $\mathcal{C}^\alpha$ for any $\alpha<1$. Indeed, one can find a loop $\gamma\in \bigcap_{\epsilon>0} \mathcal{C}^{1-\epsilon}$ such that $\n_\gamma\notin L^2(\mathbb{R}^2)$, and one can find two such curves $\gamma,\gamma'$ such that neither the positive nor the negative part of $\n_\gamma\n_{\gamma'}$ are integrable.  
\end{remark}

\section{The case of independent Brownian motions: mollified results}
\label{sec:epsilon}Recall $\mathcal{B}^{\epsilon,\epsilon'}_{W,W'}$ is defined in \eqref{eq:defB}, thanks to the function $\n_W^\epsilon$ given by \eqref{def:nepsilon}. 
The goal of this section is to prove the following formula (Proposition \ref{prop:mollif} and Lemma \ref{le:ItoStrato}), which holds for fixed~$\epsilon>0$: almost surely, 
\begin{align*}
\mathcal{B}^{\epsilon,\epsilon'}_{W,W'}  
&=-\sum_{i=1}^2	\int_0^{\bar{T}'}  	\int_0^{\bar{T}}    (\varphi^\epsilon\conv \varphi^{\epsilon'} \conv G)( W_t-W'_s  )      \circ \d W^i_{t}\circ  \d {W}'^i_{s}\\
&=	- \sum_{i=1}^2\int_0^{\bar{T}'} \hspace{-0.3cm}\big( \int_0^{\bar{T}}(\varphi^\epsilon\conv\varphi^{\epsilon'}\conv G)(W_t-W'_s)
\mathrm{d} W^i_t \big)    \mathrm{d} {{W'}}\phantom{\hspace{-0.1cm})}^i_s
- \frac{1}{4} \int_0^{T'} \hspace{-0.2cm}  \int_0^{T} (\varphi^\epsilon\conv\varphi^{\epsilon'})(W_t-W'_s)
\mathrm{d} t \d s.
\end{align*}
We will make several use of the following stochastic Fubini theorem. 
\begin{theorem}[Stochastic Fubini theorem,  {\cite[Theorem 2.2]{StochasticFubini}}]
	Let $S$ be a continuous semimartingale defined on a probability space $\Omega$, with semimartingale decomposition $S=M+A$. Let $(X,\mu)$ be a $\sigma$-finite measured space and $f: X\times \mathbb{R}_+\times  \Omega\to \mathbb{R}$ be progressively measurable. Assume the integrals
	\[ 
	\int_X \big(\int_0^T f(x,t,\omega)^2 \d \langle M\rangle_t\big)^{\frac{1}{2}} \d \mu_x \quad \text{and} \quad 
	 \int_X \int_0^T |f(x,t,\omega)| |\mathrm{d} A|_t  \d \mu_x
	\]
are almost surely finite. Then, almost surely, 
	\[ 
	\int_X \int_0^T f(x,t,\omega) \d S_t  \d \mu_x=  \int_0^T  \int_X f(x,t,\omega)   \d \mu_x \d S_t.
	\]
\end{theorem}
Note the facts all the integrals appearing in the last identity are well-defined random variables is part of the result.

\begin{lemma}
	\label{le:green1}
	Let $f\in\mathcal{C}^1(\mathbb{R}^2, \mathbb{R})$ be such that $\int_{\mathbb{R}^2\setminus B_1} \frac{|f(z)|}{ |z| } \d z<\infty$ and $\int_{\mathbb{R}^2\setminus B_1} \frac{|\grad f(z)|}{ |z| } \d z<\infty$. 
	Then, for all $\epsilon>0$, almost surely,
	\[ 
	 \int_{\mathbb{R}^2} f(z) \n^{\epsilon}_{W}(z)  \d z =
	\sum_i  \int_0^{\bar{T}} \Big( \int_{\mathbb{R}^2} f(z) \theta^{\epsilon,i}_{W_t-z}  \d z \Big) \circ \mathrm{d} W^i_t.
	\]
\end{lemma} 
\begin{proof}
	Using the definition of  $\n^\epsilon_{W}$, we see this is a simple change in the order of integration, justified by the stochastic Fubini theorem provided it can be applied here. Using the It\^o-to-Stratonovich correction formula
	to separate the martingale and the finite variation parts of  ${\int_0^T\big( \int_{\mathbb{R}^2} f(z) (\theta^{\epsilon,z}_{W_t})^i  \d z \big) \circ \mathrm{d} W^i_t}$, we see that the conditions to apply the stochastic Fubini theorem reduce in our case to the almost-sure finiteness of the integrals
	\[
	\int_{\mathbb{R}^2} |f(z)|\big( \int_0^T   |\theta^{\epsilon}_{W_t-z}|^2 \d t\big)^{\frac{1}{2}} \d z
	\quad\text{and}\quad  
	\int_{\mathbb{R}^2} \int_0^T |\div (f \theta^{\epsilon}_{W_t-\ \cdot\ })(z)|  \d t \d z. 
	\]
		
	For fixed $\epsilon>0$, $W$-almost surely,  
	$(t,z)\mapsto |\theta^{\epsilon}_{W_t-z}|$ is globally bounded. Indeed, using a decomposition $\theta=\theta_1+\theta_\infty$ with $\theta_1\in L^1,\theta_\infty\in L^\infty$, we get $\|\theta^\epsilon\|_{\infty} \leq 
	\|\theta_1\|_{L^1}\|\phi^\epsilon\|_{\infty} +
	\|\theta_\infty\|_{{\infty}}<\infty$.	
	Furthermore, for $|z|\geq |W_t|+\epsilon K_\varphi $, \[|\theta^{\epsilon}_{W_t-z}|\leq (|z|-|W_t|-\epsilon K_\varphi)^{-1}/(2\pi).\] In particular, for $|z|\geq 2(|W_t|+\epsilon K_\varphi) $, $|\theta^{\epsilon}_{W_t-z}|\leq |z|^{-1}/\pi$.
	 Thus, there exists a random constant $C$ such that for all $z$,  $|\theta^{\epsilon}_{W_t-z}|\leq C /\max(1,|z|)$. Thus, 
	\[
	\int_{\mathbb{R}^2} |f(z)| \big( \int_0^T  |\theta^{\epsilon}_{W_t-z}|^2 \d t\big)^{\frac{1}{2}} \d z
	\leq C 
	\int_{\mathbb{R}^2}   \frac{|f(z)|}{\max(1,z)}\d z<\infty.	
	\]
	Similarly, there exists a random $C'$ such that for all $z$, $|\grad \theta^{\epsilon}_{W_t-z}| \leq C' /\max(1,|z|^2)$, and thus 
	\[ 
	\int_{\mathbb{R}^2} \int_0^T |\div (f \theta^{\epsilon}_{W_t-\ \cdot\ } )(z)|  \d t \d z
	\leq 
	\int_{\mathbb{R}^2}  \frac{C' |f(z)| }{\max(1,z^2)}+ \frac{C |\grad f(z) |}{\max(1,z)}  \d z<\infty,	
	\]
	which concludes the proof.
\end{proof}
	We now prove $\theta$ is a Green kernel for the curl operator defined on the space of divergence-free vector fields:
\begin{lemma}
	\label{le:inverseCurl}
	Let $A\in\mathcal{C}^1_c(\mathbb{R}^2, \mathbb{R}^2)$ be such that $\div(A)=0$. 
	Then,  for $z \in \mathbb{R}^2$, \[A_z= \int_{\mathbb{R}^2} \curl A_{z-w} \theta_{w-z}  \d w .\]
\end{lemma}
\begin{proof}
	Let $\hat{A}^j_z$ be the $j^{\text{th}}$ coordinate of the right-hand side. 
	The fact that it is well-defined is clear, since the integrand is both compactly supported and bounded in absolute value by  the locally integrable function $w\mapsto \frac{\ \| \nabla A\|_\infty }{2 \pi |z-w|}$. Note it also follows that 
	\begin{equation}
		\label{eq:lap0}
	\hat{A}^j_z=  \lim_{\delta \to 0} \frac{1}{2\pi } \int_{\mathbb{R}^2\setminus B_\delta(z)} (\partial_1 A^2 -\partial_2 A^1)_{z-w} \theta^{j}_{w-z}  \d w.
	\end{equation}
	Let $\hat{A}^\delta_{jkl}(z)\coloneqq \int_{\mathbb{R}^2\setminus B_\delta(z)} \partial_k A^j_{z-w}  E^{l}_{w-z} \d w$. 
	Let also $e(t)=(\cos(t),\sin(t))$. Recall $E^{l}_{w-z}=\partial_l G_{w-z}$ and note $E^{l}_{\delta e(t)}=e_l(t)/(2 \pi \delta)$. The integration by part formula gives 
	\begin{align*}
	\hat{A}^\delta_{jkl}(z)
	&= 
	\int_{\mathbb{R}^2\setminus B_\delta(z)} A^j_{z-w} \partial_k E^{l}_{w-z} \d w
	+ \delta \int_0^{2\pi} A^j_{z-\delta e(t)} E^l_{\delta e(t)  } e_k(t) \d t\\
	&=_{\delta\to 0} 
	\int_{\mathbb{R}^2\setminus B_\delta(z)} A^j_{z-w} \partial_{kl} G_{w-z}  \d w
	+ \frac{A^j_z}{2\pi} \int_0^{2\pi}  e_k(t) e_l(t) \d t+o(1)\\
	&=_{\delta\to 0} 
	\int_{\mathbb{R}^2\setminus B_\delta(z)} A^j_{z-w} \partial_{kl} G_{w-z}  \d w
	+ \frac{A^j_z}{2} \mathbbm{1}_{l=k} +o(1).
	\end{align*}
Specifically, the $o(1)$ term is smaller than $\delta \|\grad A^j\|_{\infty}$.
	Since $\Delta G^z =0$ on $\mathbb{R}^2\setminus \{z\}$, 
	we deduce that
	\begin{equation}
		\label{eq:lap1} 
	\hat{A}_{j11}^\delta(z)
		+\hat{A}_{j22}^\delta(z)
		=_{\delta\to 0} 
		\int_{\mathbb{R}^2\setminus B_\delta(z)} A^j_{z-w} \Delta G_{w-z}  \d w
		+  A^j_z +o(1)\underset{\delta\to 0}\longrightarrow A^j_z.
	\end{equation}
	Furthermore, 
	\begin{equation}		
	\label{eq:lap2}
	\hat{A}_{j12}^\delta(z)
	-\hat{A}_{j21}^\delta(z)
	=_{\delta\to 0} 
	\int_{\mathbb{R}^2\setminus B_\delta(z)} A^j_{z-w} (\partial_{12}-\partial_{21}) G_{w-z}  \d w
	+0 +o(1)\underset{\delta\to 0}\longrightarrow 0. 
	\end{equation}
	Furthermore, by definition of $	\hat{A}_{jkl}^\delta$, and since $\div(A)=0$, 
	\begin{equation}
	\label{eq:lap3}
	\hat{A}_{11l}^\delta+\hat{A}_{22l}^\delta=  
	\int_{\mathbb{R}^2\setminus B_\delta(z)} \div A_{z-w} \partial_l G_{w-z}  \d w= 0,
	\end{equation}
	Using \eqref{eq:lap0}, then \eqref{eq:lap1} with $j=1$, then \eqref{eq:lap2} with $j=2$, then \eqref{eq:lap3} with $l=1$, we get 
	\[ 
	\hat{A}^1=\lim_{\delta\to 0} (\hat{A}_{122 }^\delta-\hat{A}_{212}^\delta)
	=\lim_{\delta\to 0} (A^1 - \hat{A}_{111}^\delta-\hat{A}_{212}^\delta)
	=\lim_{\delta\to 0} (A^1 - \hat{A}_{111}^\delta-\hat{A}_{221}^\delta)
=A^1,\]
and similarly 
\[
	\hat{A}^2=
	\lim_{\delta\to 0} (\hat{A}_{211 }^\delta-\hat{A}_{121}^\delta)
	=\lim_{\delta\to 0} (A^2 - \hat{A}_{222}^\delta-\hat{A}_{121}^\delta)
	=\lim_{\delta\to 0} (A^2 - \hat{A}_{222}^\delta-\hat{A}_{112}^\delta)
	=A^2,
	\]
	which concludes the proof.
\end{proof}

\begin{corollary}[mollified stochastic Green's formula]
	\label{coro:green2}
	Let $A\in\mathcal{C}^2(\mathbb{R}^2, \mathbb{R}^2)$. Let $A^\epsilon=A\conv \varphi^\epsilon$. Then, almost surely,
	\[ 
	\int_{\mathbb{R}^2} \curl A_z  \n^\epsilon_{W}(z)  \d z =
	 \int_0^{\bar{T}}   A^{\epsilon}_{W_t}    \circ \mathrm{d} W_t.
	\]
\end{corollary} 
\begin{proof}
	First, remark that neither side depend on the value of $A$ outside $B_{\|W\|_\infty+\epsilon K}$, so we can assume without loss of generality that $A$ is compactly supported. Furthermore, both sides are invariant by the addition to $A$ of a gradient term: for the left-hand side, this is  because $\curl \circ \grad=0$; for the right-hand side, this is because the Stratonovich integral of a gradient along a closed loop is equal to $0$. 
	Thus, we can assume that $A$ is compactly supported and divergence-free without loss of generality (recall a compactly supported vector field can be decomposed into the sum of a compactly supported gradient part and a compactly supported divergence-free part: this is essentially Helmoltz or Hodge decomposition theorem).
		
	 Then, the function $f=\curl A $ satisfies the assumption of Lemma \ref{le:green1}, and we get 
		\[ 
	\int_{\mathbb{R}^2} \curl A_z \n^\epsilon_{W}(z)  \d z =
	\sum_i  \int_0^T\big( \int_{\mathbb{R}^2} \curl A_z \theta^{\epsilon,i}_{W_t-z}  \d z \big) \circ \mathrm{d} W^i_t.
	\]
	For all $z, w\in \mathbb{R}^2$, it holds $\theta^{\epsilon}_{w-z}=  - (\theta_{\ \cdot\ -w} \conv \varphi^\epsilon)(z)$. Using first Fubini's theoren and then Lemma \ref{le:inverseCurl}, we deduce that for all $t\in[0,T]$ 
	\[
	 \int_{\mathbb{R}^2} \curl A_z \theta^{\epsilon,i}_{W_t-z}  \d z
	 =
	\int_{\mathbb{R}^2} \curl A^\epsilon_z \theta^{i}_{W_t-z}  \d z
	= A^{\epsilon,i}_z,\]
	which concludes the proof.
\end{proof}

\begin{proposition}[Mollified stochastic Amperean area formula]\label{prop:mollif} 
	For all $\epsilon,\epsilon'>0$, almost surely, 
	\[ 
	\mathcal{B}^{\epsilon,\epsilon'}_{W,W'}  
	= 
	-\sum_i \int_0^{\bar{T}}\big(\int_0^{\bar{T}'} \varphi^\epsilon\conv \varphi^{\epsilon'}\conv G (W_t-W'_s) \circ \mathrm{d} W'^i_s\big)\circ \mathrm{d} W_t^i.
	\] 
\end{proposition}
\begin{proof} 
	Let again $G^{1}_w=(G_w,0)$, $G^{2}_w=(0,G_w)$, so that $\curl G^{i}=\theta^{i}$. 
	Define 
	\[
	\mathcal{A}^{W,\epsilon,j}_z\coloneqq - \int_0^{\bar{T}} (\varphi^\epsilon \conv G^{j})_{ W_t-z}\circ \mathrm{d} W_t^j.
	\]

	We apply iteratively Lemma~\ref{le:derivBelowRP} (derivation below rough path integrals), to the smooth function $(x,z) \mapsto (\varphi^\epsilon \conv G^{j})_{ x-z}$  and its derivatives. We deduce that for any multi-index $I$, 
	\[ 
	- \partial_I \big(  \int_0^{\bar{T}} (\varphi^\epsilon \conv G^{j})_{ W_t-\ \cdot\ }\mathrm{d} W_t^j\big)(z) = - \int_0^{\bar{T}} 	(\partial_I(\varphi^\epsilon \conv G^{j}))_{ W_t-z} \mathrm{d} W_t^j.
	\] 
Applying the Stratonovich-to-Itô correction formula and the theorem of derivation below an integral (the classical one for Lebesgue integrals), we deduce 
	\[ 
	\partial_I 	\mathcal{A}^{W,\epsilon,j}_z= - \int_0^{\bar{T}} 	(\partial_I(\varphi^\epsilon \conv G^{j}))_{ W_t-z}\circ \mathrm{d} W_t^j.
	\] 
	In particular, 
	\begin{align*}
	\curl \mathcal{A}^{W,\epsilon}_z
	&= 
	\int_0^{\bar{T}} \varphi^\epsilon \conv \partial_1  G^{2}_{W_t-z} \circ \mathrm{d}W_t^2
	-\int_0^{\bar{T}} \varphi^\epsilon \conv \partial_2  G^{1}_{W_t-z} \circ \mathrm{d} W_t^1\\
	&	= \int_0^{\bar{T}} (\varphi^\epsilon \conv \theta_{\ \cdot\ -z}  )_{W_t} \circ \mathrm{d} W_t= \n^\epsilon_{W}(z),
	\end{align*}
	and since  the  vector field $\mathcal{A}^{\epsilon}_W$ is almost surely smooth, we can apply Corollary 
	\ref{coro:green2} (to the Brownian motion $W'$) to deduce that 
	\[ 
	\int_{\mathbb{R}^2} \n^\epsilon_{W}(z)  \n^{\epsilon'}_{W'}(z)  \d z =
	\sum_i  \int_0^{\bar{T}'} \varphi^{\epsilon'} \conv \mathcal{A}^{W,\epsilon,i}_{W'_s}    \circ \mathrm{d} W'^i_s
	=-\sum_i  \int_0^{\bar{T}'}
	\int_0^{\bar{T}} \varphi^{\epsilon'} \conv \varphi^\epsilon \conv G(W_t-W'_s)\circ \mathrm{d} W_t^i
	\circ \mathrm{d} W'^i_s.
	\]
\end{proof}

\begin{lemma}
	\label{le:ItoStrato}
	We have the It\^o-to-Stratonovich correction formula 
$\mathcal{B}^{\epsilon,\epsilon'}_{W,W'} 
		=  \mathcal{B}^{\epsilon,\epsilon',\mathrm{It\hat{o}}}_{W,W'}-\frac{1}{4} I^{\epsilon,\epsilon'}_{W,W'}$,
		where 
		\[ 
		\mathcal{B}^{\epsilon,\epsilon',\mathrm{It\hat{o}}}_{W,W'}\coloneqq
		-\sum_{i=1}^2\int_0^{\bar{T}'} \hspace{-0.3cm} \big( \int_0^{\bar{T}}\hspace{-0.1cm} (\varphi^\epsilon\conv\varphi^{\epsilon'}\conv G)_{W_s-W'_t}
		\mathrm{d} W^i_s \big)    \mathrm{d} {{W'}}\phantom{\hspace{-0.1cm})}^i_t \ \text{and} \  
		I^{\epsilon,\epsilon'}_{W,W'}\coloneqq \int_0^{T'}  \hspace{-0.3cm}\int_0^{T} (\varphi^\epsilon\conv\varphi^{\epsilon'})_{W_s-W'_t}
		\mathrm{d} s \d t.\]
\end{lemma}
\begin{proof}	
	Let $f$ be a smooth function from $\mathbb{R}^2$ to $\mathbb{R}$, and let $r_{s,t}=W_s-W'_t$.
	Then, in differential notations, 
	\begin{align*}
		f_{W_s-W'_t}
		\circ \mathrm{d} W^i_s  \circ   \mathrm{d} {{W'}}^i_t
		=& \big(f_{W_s-W'_t}
		\d W^i_s+ \frac{\mathbbm{1}_{s\leq T} }{2} \partial_i f_{W_s-W'_t} \d s    \big) \d {{W'}}^i_t
		-\frac{\mathbbm{1}_{t\leq T'}}{2} \partial_i f_{W_s-W'_t}
		\circ \mathrm{d} W^i_s \d t\\
		=& f_{W_s-W'_t}
		\d W^i_s \d {{W'}}^i_t    + \frac{\mathbbm{1}_{s\leq T} }{2} \big(\partial_i f_{W_s-W'_t}  \d s    \circ \mathrm{d} {{W'}}^i_t 
		+\frac{\mathbbm{1}_{t\leq T'}}{2} \partial_{ii} f_{W_s-W'_t}  \d s\d t
		\big)\\
&\hspace{6.7cm}		-\frac{\mathbbm{1}_{t\leq T'}}{2} \partial_i f_{W_s-W'_t}
		\circ \mathrm{d} W^i_s \d t.
	\end{align*}
Remark the sign in front of the $\partial_{ii}$ term: a minus sign comes from the fact we look at the derivative of  $x\mapsto f_{W_s-x}$, a second minus sign comes from the fact we are converting It\^o differential back into Stratonovich differential (and not the other way around). 

	Since the Stratonovich integral of a gradient vector field along a loop is equal to $0$, for all $t$,
	$
	\sum_i \int_0^{\bar{T}}\partial_i f_{W_s-W'_t}
	\circ \mathrm{d} W^i_s=0$.
	
	Similarly, 
	$
	\sum_i \int_0^{\bar{T}'} \int_0^{{T}} \partial_i f_{W_s-W'_t} \d s 
	\circ \mathrm{d} W'\phantom{\hspace{-0.1cm})}^i_t=
	\sum_i \int_0^{\bar{T}'} \partial_i \big(\int_0^{{T}}  f (W_s-\ \cdot\ ) \d s\big)_{W'_t} 
	\circ \mathrm{d} W'\phantom{\hspace{-0.1cm})}^i_t=
	0$, and we are left with
	\begin{multline*}
		\sum_{i=1}^2\int_0^{\bar{T}'}\hspace{-0.3cm} \int_0^{\bar{T}}f_{W_s-W'_t}
		\circ \mathrm{d} W^i_s \circ   \mathrm{d} {{W'}}\phantom{\hspace{-0.1cm})}^i_t
		= \sum_{i=1}^2\int_0^{\bar{T}'} \hspace{-0.3cm}  \int_0^{\bar{T}}f_{W_s-W'_t}
		\mathrm{d} W^i_s  \mathrm{d} {{W'}}\phantom{\hspace{-0.1cm})}^i_t
		+ \frac{1}{4}\int_0^{T'} \hspace{-0.3cm} \int_0^{T}\Delta f_{W_s-W'_t}
		\mathrm{d} s \d t.	
	\end{multline*}
	We obtain the lemma by taking $f=\frac{1}{2\pi}\varphi^\epsilon\conv\varphi^{\epsilon'}\conv G$, for which $\Delta f =\varphi^\epsilon\conv\varphi^{\epsilon'}$. 
\end{proof}

\section{The case of independent Brownian motions: the \texorpdfstring{$\epsilon\to 0$ }{ 𝝐→ 0} limit}
\label{sec:zero}
We define $I:\mathbb{R}^2 \to \mathbb{R}_+$ the jointly continuous intersection local time of $W$ (up to time $T$) with $W'$ (up to time $T'$), as defined by \cite[Theorem 1]{Geman} (where it is called $\alpha$). Formally, 
\[ I_z = ``\int_0^T \int_0^{T'} \delta_{z}(W_s-W'_t) \d t\d s". \]
It satisfies the property that for any bounded Borel function $g:\mathbb{R}^2 \to \mathbb{R}$ and any $s,t\in [0,1]$, almost surely, 
\begin{equation}
	\label{eq:loctimeprop}
\int_{\mathbb{R}^2} g_z I_z \d z=\int_0^T\int_0^{T'}  g_{W_s-W'_t} \d t \d s.
\end{equation}
In particular, the quantity $I^{\epsilon,\epsilon'}_{W,W'}$ defined in Lemma \ref{le:ItoStrato} is given by $(\varphi^\epsilon\conv\varphi^{\epsilon'}\conv I) (0)$ 
We will write $I_{W,W'}$ for the value $I_0$. 
We also define 
	\[ \mathcal{B}_{W,W'}^{\mathrm{It\hat{o}}} =- \sum_{i=1}^2 \int_0^{\bar{T'}} \Big( \int_0^{\bar{T}} G(W_s-W'_t)
	\mathrm{d} W^i_s \Big)  \mathrm{d} {{W'}}\phantom{\hspace{-0.1cm})}^i_t, \qquad 
	\mathcal{B}_{W,W'}= \mathcal{B}_{W,W'}^{\mathrm{It\hat{o}}}-\frac{I_{W,W'}}{4}.\]

\begin{proposition}
	\label{prop:cvY}
	As $\epsilon,\epsilon'\to 0$, $\mathcal{B}^{\epsilon,\epsilon'}_{W,W'}$ converges in probability. The limit is given by $\mathcal{B}_{W,W'}$. 
\end{proposition}
\begin{proof}
	In proposition~\ref{prop:mollif} we obtained the formula for $\mathcal{B}^{\epsilon,\epsilon'}_{W,W'}$ as an iterated Stratonovich integral. In Lemma~\ref{le:ItoStrato}, we proved this iterated Stratonovich integral can be decomposed as the sum of the iterated It\^o integral  $\mathcal{B}_{W,W'}^{\epsilon,\epsilon',\mathrm{It\hat{o}}}$ and the quadratic variation term $I_{W,W'}^{\epsilon,\epsilon'}$. 
	
	By \eqref{eq:loctimeprop},  the fact $\varphi^\epsilon \conv \varphi^{\epsilon'}$ is a mollifier with integral one, and the continuity of $I$ at $0$, we deduce that almost surely (hence also in probability),
	\begin{equation} 
		\label{eq:cvLoctime}
	I_{W,W'}^{\epsilon,\epsilon'}=	\int_0^T\int_0^{T'} (\varphi^\epsilon \conv \varphi^{\epsilon'})_{W_u-W'_v}  \d v \d u
		=\int_{\mathbb{R}^2}  (\varphi^\epsilon \conv \varphi^{\epsilon'})_z I_z \d z 
		\underset{\epsilon,{\epsilon'}\to 0}\longrightarrow I_{W,W'}.
	\end{equation} 
	
	The conclusion follows from the fact $\mathcal{B}_{W,W'}^{\epsilon,\epsilon',\mathrm{It\hat{o}}}$ converges toward $\mathcal{B}_{W,W'}$ in $L^2(\mathbb{P})$, thus in probability. To prove this, it suffices to prove that for $i\in \{1,2\}$,
	\begin{equation}
		\Big\|   
		\int_0^{\bar{T}'} \Big( \int_0^{\bar{T}} ((\varphi^\epsilon\conv\varphi^{\epsilon'} -\delta )\conv G) (W_s-W'_t)
		\mathrm{d} W^i_s \Big)  \mathrm{d} {{W'}}\phantom{\hspace{-0.1cm})}^i_t\Big\|_{L^2(\mathbb{P})}\underset{\epsilon,\epsilon' \to 0}\longrightarrow 0. \label{eq:conv} 
	\end{equation}
	This is fairly technical and postponed to Appendix \ref{sec:appendixTechnic} (Lemma \ref{le:mollifok} and \ref{le:ap2}). 
\end{proof}

In order to later improve the type of convergence, we now show the moments of $\mathcal{B}^{\epsilon,\epsilon',\mathrm{It\hat{o}}}_{W,W'}$ and $I^{\epsilon,\epsilon'}_{W,W'}$ are bounded by their respect values at $\epsilon=\epsilon'=0$. 
\begin{lemma} 
	\label{le:swapswap}
	Let $g:z \in \mathbb{R}^2 \mapsto \mathcal{B}^{\mathrm{It\hat{o}}}_{W+z,W'}$. 
	For all $\epsilon,\epsilon'>0$, it holds $\mathcal{B}^{\epsilon,\epsilon',\mathrm{It\hat{o}}}_{W,W'}= (\varphi^\epsilon\conv\varphi^{\epsilon'}\conv g)(0)$.
\end{lemma} 
\begin{proof}
	We start with the expression of $\mathcal{B}^{\epsilon,\epsilon',\mathrm{It\hat{o}}}_{W,W'}$ provided by Lemma \ref{le:ItoStrato}. 
	Formally, it only remains to change the order of integration again in order to conclude (i.e. to take the convolution with $\phi^\epsilon\conv \phi^{\epsilon'}$ outside from the stochastic integrals, rather than inside), for which we need to apply the stochastic Fubini lemma twice. Verifying that the assumption holds is straightforward: to prove that for almost all $t$,
	\[\int_0^{T'} \int_{\mathbb{R}^2} \varphi^\epsilon\conv\varphi^{\epsilon'}(z)\ G(W_s-W'_t+z) \d z \d W^i_s	
	= \int_{\mathbb{R}^2}\int_0^{T'} \varphi^\epsilon\conv\varphi^{\epsilon'}(z) G(W_s-W'_t+z) \d W^i_s \d z, 
	\]
	we need to check that almost surely,
	\[ 
	A= \int_{\mathbb{R}^2}\Big(\int_0^{T'} \varphi^\epsilon\conv\varphi^{\epsilon'}(z)^2 G(W_s-W'_t+z)^2 \d s\Big)^{1/2} \d z<\infty.
	\]
	To this end, we prove that for all $R$ finite, $A\mathbbm{1}_{E_R}$ is bounded in $L^2(\Omega)$, where $E_R$ is the event $\{\Range(W)\cup \Range(W')\subset B_R\}$. Let $a,a'$ be the starting points of $W,W'$. Then, 
	\begin{align*}
	\mathbb{E}[A \mathbbm{1}_{E_R} ]^2
	&\leq \| \varphi^\epsilon\conv\varphi^{\epsilon'}\|_\infty^2 
	\mathbb{E}\Big[  \mathbbm{1}_{E_R} \int_{B_{K_\varphi}} \int_0^{T'}  G(W_s-W'_t+z)^2 \d s \d z \Big]\\
	&\leq \| \varphi^\epsilon\conv\varphi^{\epsilon'}\|_\infty^2 
	\int_{B_R^2\times B_{K_\varphi}} \int_0^{T'}   p_s(a,x)p_t(a',y)  G(x-y+z)^2 \d s  \d x \d y\d z \\
	&\leq \| \varphi^\epsilon\conv\varphi^{\epsilon'}\|_\infty^2 \frac{|B_{K_\varphi}|}{2\pi t}
	\int_{B_R} (\int_0^{T'}   p_s(a,x) \d s) (\int_{B_{R+K_\varphi}} G(x-w)^2\d w  ) \d x <\infty. 
	\end{align*}
	We can thus apply the stochastic Fubini theorem. Since we can also use the classical Fubini theorem to control the integral from $T'$ to $\bar{T}'$, we deduce that 
	\[\int_0^T \int_0^{\bar{T}'}
	 \varphi^\epsilon\conv\varphi^{\epsilon'} \conv G(W_s-W'_t) \d W^i_s	\d W'^i_t
	=	\int_0^T 
	\int_{\mathbb{R}^2}\int_0^{\bar{T}'} \varphi^\epsilon\conv\varphi^{\epsilon'}(z) G(W_s-W'_t+z) \d W^i_s\d z	\d W'^i_t.
	\]
	To swap the two uttermost integrals, we apply the stochastic Fubini theorem again. We need to check that almost surely, 
	\[ 
	B=\int_{\mathbb{R}^2}\Big(\int_0^{T'} \Big(\int_0^{T} \varphi^\epsilon\conv\varphi^{\epsilon'}(z) G(W_s-W'_t+z) \d W^i_s\Big)^2 \d t   \Big)^{1/2} \d z<\infty. 
	\]
	From computation similar to the previous one, we deduce 
	\begin{align*}
	\mathbb{E}[B\mathbbm{1}_{E_R}]^2
	&\leq \| \varphi^\epsilon\conv\varphi^{\epsilon'}\|_\infty^2 \int p_s(a,x)p_t(a',y) G(x-y+z) \d s \d t \d x \d y \d z\\
	&\leq \| \varphi^\epsilon\conv\varphi^{\epsilon'}\|_\infty^2  \int \frac{\Gamma(0, |x-a|^2/2 )\Gamma(0, |y-a'|^2/2 )}{4\pi^2}  G(x-y+z) \d x \d y \d z,
	\end{align*}
	where the range of integration is a compact set, namely $B_R\times B_R\times B_{(\epsilon+\epsilon') K_\varphi}$, and where $\Gamma$ is the incomplete Gamma function. Since both $\Gamma(0,|\cdot|)$ and $G$ are locally integrable (they both have logarithmic divergences at the origin), we conclude by performing first the integration with respect to the variable $z$ and then the other two integrations.
	Again, the parts from $T$ to $\bar{T}$ and from $T'$ to $\bar{T}'$ are controlled in a similar but simpler way, we leave it as an exercise.
\end{proof}
This allows us to deduce the following bounds (in which, from translation invariance, we can fix any $b\in \mathbb{R}^2$ and replace $\sup_{b,b'}$ with $\sup_{b'}$).
\begin{lemma}
	\label{le:bound}
	For all $a,a'\in \mathbb{R}^2$, $T,T',\epsilon,\epsilon'>0$ and $k\in \mathbb{N}$,
	\[ 
	\mathbb{E}_{T,T',a,a'}[  |\mathcal{B}^{\epsilon,\epsilon',\mathrm{It\hat{o}}}_{W,W'}|^k]
	\leq \sup_{b,b'}
	\mathbb{E}_{T,T',b,b'}[  |\mathcal{B}^{\mathrm{It\hat{o}}}_{W,W'}|^k] \quad \mbox{and} \quad 
	\mathbb{E}_{T,T',a,a'}[  |I^{\epsilon,\epsilon'}_{W,W'}|^k]
	\leq \sup_{b,b'}
	\mathbb{E}_{T,T',b,b'}[  |I_{W,W'}|^k].
	\]
\end{lemma}
\begin{proof}
Let $\psi=\varphi^\epsilon\conv \varphi^{\epsilon'}$ and $g$ as in Lemma \ref{le:swapswap}. 
Using the equality given by Lemma \ref{le:swapswap} and the fact $\psi$ is non-negative with integral $1$, we have 
\begin{align*}
\mathbb{E}_{T,T',a,a'}[  |\mathcal{B}^{\epsilon,\epsilon',\mathrm{It\hat{o}}}_{W,W'}|^k]
&=	\mathbb{E}_{T,T',a,a'}[ | \int_{(\mathbb{R}^2)k} \psi(z_1)\dots \psi(z_k) g(z_1)\dots g(z_k)\d z_1\dots \d z_k|]\\
&\leq	\int_{(\mathbb{R}^2)k} \psi(z_1)\dots \psi(z_k) \mathbb{E}_{a,a'}[|g(z_1)\dots g(z_k)|]\d z_1\dots \d z_k\\
&\leq \sup_z \mathbb{E}_{T,T',a,a'}[|g(z)|^k]=\sup_{b}\mathbb{E}_{T,T',a,b}[|\mathcal{B}^{\mathrm{It\hat{o}}}_{W,W'}|^k], \qquad (b=a'-z).
\end{align*} 	

The same proof also shows the second inequality,  simply replacing the function $g$ with the function $I$.

Remark we do not claim yet these moments are finite, but this will be proved later.
\end{proof}

\begin{corollary}
	\label{coro:unifBound} \ 
	\begin{itemize}
	\item 	Let $0 \leq \beta^{\mathrm{It\hat{o}} }$. 
	 Assume $\sup_{a,a'} \mathbb{E}_{a,a'}[ \exp( \beta |\mathcal{B}^{\mathrm{It\hat{o}} }_{W,W'}|  )   ]<\infty$ for all $\beta\in[0,\beta^{\mathrm{It\hat{o}} })$. Then, for all $\beta'\in(-\beta^{\mathrm{It\hat{o}} },\beta^{\mathrm{It\hat{o}} })$,  
	$\exp( \beta' \mathcal{B}^{\epsilon,\epsilon',\mathrm{It\hat{o}}}_{W,W'}  ) $ converges in $L^1(\mathbb{P}_{a,a'})$ toward 
	$\exp( \beta' \mathcal{B}^{\mathrm{It\hat{o}} }_{W,W'}  ) $ as $\epsilon,\epsilon'\to 0$. 
	
	\item Let $  \beta^{I }\geq 0$. Assume $\sup_{a,a'} \mathbb{E}_{a,a'}[ \exp( \beta |I_{W,W'}|  )   ]<\infty$ for all $\beta<\beta^{I }$. Then, for all $\beta'\in  (-\beta_+^I, \beta_+^{I})$, 
	$\exp( \beta' I^{\epsilon,\epsilon'}_{W,W'}  ) $ converges in $L^1(\mathbb{P}_{a,a'})$ toward 
	$\exp( \beta' I_{W,W'}  ) $  as $\epsilon,\epsilon'\to 0$. 
	
\end{itemize}
\end{corollary}
\begin{proof}
	Let $\beta'\in(-\beta^{\mathrm{It\hat{o}} },\beta^{\mathrm{It\hat{o}} })$, and let $p>1$ be such that $\beta\coloneqq p |\beta'|$ is smaller than $\beta^{\mathrm{It\hat{o}} }$.
	Then, the assumption ensures that $\exp( p \beta' |\mathcal{B}^{\mathrm{It\hat{o}}}_{W,W'}|  )$ is finite. 
	By Lemma \ref{le:bound}, it follows that the $L^p$ norm of $\exp( \beta' \mathcal{B}^{\epsilon,\epsilon',\mathrm{It\hat{o}}}_{W,W'}  )$ is bounded uniformly in $\epsilon,\epsilon'$, thus the random variables $\exp( \beta' \mathcal{B}^{\epsilon,\epsilon',\mathrm{It\hat{o}}}_{W,W'}  )$ are uniformly integrable. Since furthermore 
	$\exp( \beta' \mathcal{B}^{\epsilon,\epsilon',\mathrm{It\hat{o}}}_{W,W'}  )$ converges in probability toward $\exp( \beta' \mathcal{B}^{\mathrm{It\hat{o}}}_{W,W'}  )$, 
	Vitali convergence theorem allows to conclude that 
	$\exp( \beta' \mathcal{B}^{\epsilon,\epsilon',\mathrm{It\hat{o}}}_{W,W'}  )$ converges in $L^1$ toward $\exp( \beta' \mathcal{B}^{\mathrm{It\hat{o}}}_{W,W'}  )$. 
	The result for the intersection local time is similar.
\end{proof}

\section{Exponential moments} 
\label{sec:expo}
In this section, we will prove that the Amperean area admits some finite exponential moments, i.e. there exists $\beta_c>0$ such that for all $\beta\in(-\beta_c,\beta_c)$, 
\begin{equation}
\label{eq:strech}
	\mathbb{E}\big[\exp\big( \beta \int_0^{\bar{T}}\int_0^{\bar{T}'}  G(W_s-W'_t) \circ \d W_s \circ\d W'_t  \big)\big]<\infty.
  \end{equation}
We will track down how $\beta_c$ depends on $T,T'$. From scale invariance, if both $T$ and $T'$ are multiplied by the same factor $\lambda$, then $\beta_c$ must be divided by $\lambda$, but we will show the finer result that we can choose $\beta_c$ proportional to $1/\sqrt{T T'}$.\footnote{Recall we expect the partition function here to be given, up to counterterms, by $\mathbb{E}[ \exp(-\alpha^2 \sum_{\ell,\ell' \in \mathcal{L} }  \mathcal{B}_{\ell,\ell'} )]$, where $\mathcal{L}$ is a Brownian loop soup. If we wish to bound the contribution from the infinitely many small loops $\ell$, say by applying Hõlder's inequality for infinite products, for which the exponents $p_k$ goes to infinity, we must be able, when either $\ell$ or $\ell'$ has a small size, to control a  $\mathbb{E}[ \exp(-\beta \mathcal{B}_{\ell,\ell'} )]$ for a \emph{large} $\beta$.}
  
We will use the Burkholder inequality to bound the polynomial moments of $\mathcal{B}_{W,W'}$, with the asymptotically optimal constant given by Zakai \cite[Theorem 1]{Zakai}, as we need to keep track of the factors $k^k$ when looking at the $k^{th}$ moment. 
\begin{theorem}
	\label{th:Zakai}
	Let $(\Omega,\mathcal{F},\mathbb{P})$ be a probability space endowed with a filtration $(\mathcal{F}_s)_{s\geq 0}$, $W$ be a standard Brownian motion adapted to $(\mathcal{F}_s)_{s\geq 0}$, and $Y$ a process which is $\mathcal{F}_s$-progressively measurable and real-valued. Then,  for $p>1$,
	\begin{equation}
		\label{eq:burkholder}
	\mathbb{E}\Big[ \Big| \int_0^T Y_s \d W_s \Big|^{2p} \Big]  \leq (2p-1)^{p} \Big( \int_0^T \mathbb{E}[|Y_s|^{2p} ]^{\frac{1}{p}} \d s\Big)^{p}\leq  (2p)^{p} T^{p-1} \int_0^T \mathbb{E}[|Y_s|^{2p} ]  \d s.
	\end{equation}
\end{theorem} 
Remark the direct use of the Burkholder inequality, even with the optimal constant, would be far from sufficient to conclude, as it would only gives finiteness of stretch-exponential moments,  
	\[\mathbb{E}\big[\exp\big( \beta \big|\int_0^{\bar{T}}\int_0^{\bar{T}'}  G(W_s-W'_t) \circ \d W_s \circ\d W'_t \big|^\eta \big)\big]<\infty,\]
for $\eta<\frac{2}{3}$ (and $\eta=\frac{2}{3}$ with some extra work): we invite the reader to try by themself. The main issue is that we would end up with an integral, in the vicinity of $0$, of $|\log|z||^k$, which grows extremely fast as $k \to \infty$. The problem is essentially that we would have apply the \emph{second} inequality in Theorem \ref{th:Zakai}, which is very rough here, as the $L^k$ norm of  $\log|z|$ is much bigger than its $L^2$ norm. 

We will also use the Hardy--Littlewood rearrangement inequality:
\begin{lemma}
	Let $f,g: \mathbb{R}^n\to \mathbb{R}_+$ be measurable functions, and assume the super-level sets 
	$\{ x: f(x)\geq a \mbox{ or } g(x)\geq a\}$ have finite measure for all $a>0$. Then, 
	\[ 
	\int_{\mathbb{R}^n} f(x) g(x) \d x \leq 	\int_{\mathbb{R}^n} f^*(x) g^*(x) \d x, 
	\]
	where $f^*$ (resp. $g^*$) is the symmetric decreasing rearrangement of $f$ and $g$.
\end{lemma}
For our purpose we do not need the whole definition of the symmetric decreasing rearrangement, but only the following facts: 
\begin{itemize} 
	\item The symmetric decreasing rearrangement is invariant under translation: if there exists $a$ such that $f(x+a)=g(x)$ for all $x$, then $f^*=g^*$.
	\item If there exists a non-decreasing function $\tilde{f}:\mathbb{R}_+\to \mathbb{R}_+$ such that $f(x)=\tilde{f}(|x|)$, then $f^*=f$. 
	\item For a measurable set $A$, $(\mathbbm{1}_A)^*=\mathbbm{1}_B$ where $B$ is the ball centred at $0$ with the same volume as $A$. 
\end{itemize}	
It follows from the first two points that for any $a$, $p_t(a,\cdot )^*=p_t(0,\cdot)$.

\medskip 
For a non-negative integer $k$, we define the tetrahedra 
\[\mathbb{T}^k_T=\{ (t_1,\dots, t_k)\in [0,T]^k: t_0\leq t_1\leq \dots \leq t_k \leq T \}. \]
\begin{lemma} 
	\label{le:integralT}
	Let $s>-1$. There exists $C_s<\infty$ such that for all $k\in \mathbb{N}\setminus\{0 \}$ and $T>0$, 
	\[ 
	\int_{ \mathbb{T}^k_T} \prod_{i=1}^k (t_{i}-t_{i-1})^s \d t_1\dots \d t_k
	= T^{(1+s)k} \frac{\Gamma(1+s)^k}{\Gamma(1+k(1+s)) }\leq  C_s^k T^{(1+s)k}  k^{-k(1+s)}. 
	\]
	For $s=-1/2$, there is a universal constant $C<\infty$ such that for all $k\in \mathbb{N}\setminus\{0 \}$ and $T>0$, 
	\[ \int_{ \mathbb{T}^k_T} \prod_{i=1}^k (t_{i}-t_{i-1})^{-\frac{1}{2}} \d t_1\dots \d t_k \leq C (2\pi e)^\frac{k}{2}   T^{\frac{k}{2} }  k^{-\frac{k}{2}}.
	\]
\end{lemma} 
\begin{proof} 
	Let $f(k,T)$ be the right-hand side. The change of variables $t'_i=T^{-1} t_i$ gives $f(k,T)= T^{(1+s) k} f(k,1)$.
	Then, using the change of variables $u_i=t_i-t_{i-1}$ ($u_1=t_1$ instead) and isolating the last variable $u=u_{k+1}$, we get
	\[ 
	f(k+1,1)=\int_0^1\hspace{-0.2cm} u^s  f(k,1-u)  \d u=   f(k,1) \int_0^1 u^s (1-u)^{(1+s)k }  \d u
	=  f(k,1)  \frac{ \Gamma(1+s) \Gamma(1+k(1+s)) }{ \Gamma(1+  (k+1)(1+s))   },
	\]
	from which we deduce 
	\[ 
	f(k,1)= \prod_{j=0}^{k-1} \frac{ \Gamma(1+s) \Gamma(1+j(1+s)) }{ \Gamma(1+  (j+1)(1+s))   }
	= \frac{\Gamma(1+s)^k}{\Gamma(1+k(1+s)) },
	\]
	from which we deduce the equality in the lemma. The two inequalities in the lemma are then obtained by using Stirling's approximation and $\Gamma(1/2)=\sqrt{\pi}$. 
\end{proof}

We can finally start our estimation proper. 
\begin{lemma}
	\label{le:finally}
	There exists a universal constant $C'$ such that for all $T>0$, for all integer $k\geq 2$, for $i\in \{1,2\}$,
	\[ 
	\sup_{a,b\in \mathbb{R}^2, a\neq b}
	\mathbb{E}_{T,a} \Big[ \Big|  	\int_0^{\bar{T}} G(W_s-b) \d W^i_s\Big|^{2k}\Big]  \leq C' (18/\pi^2)^{k}  T^{k} k^{k}.
	\]
\end{lemma}
\begin{proof}
	By translation invariance, we can assume $b=0$ (and then  $a\neq 0$).  
	Be careful $s\mapsto G(W_s)$ is notoriously a local martingale but not a true martingale. 
	First remark for all $s\in[0,T]$, it holds \cite[(5.2)]{MansuyYor}
	\[ G(W_s)-G(W_0)= \int_0^s \frac{W^1_u\d W^1_u+W^2_u\d W^2_u}{2\pi |W_u|^2}.
	\]
	In particular, the quadratic covariation $[G(W), W^i ]_s$ is given by $\int_0^t \frac{W^i_s}{2\pi |W_s|^2} \d s$.
	By It\^o's lemma applied to $f:t\in [0,\bar{T}]\mapsto W^i_tG(W_t)$, and since $W^i_{\bar{T}}G(W_{\bar{T}})=W^i_{0}G(W_{0})=aG(a)$, 
	\[
	0=f(\bar{T})-f(0)= \int_0^{\bar{T}} W^i_s \frac{W^1_s\d W^1_s+W^2_s\d W^2_s}{2\pi|W_s|^2}
	+\int_0^{\bar{T}} G(W_s) \d W^i_s +\frac{1}{2}\int_0^T \frac{W^i_s}{2\pi |W_s|^2} \d s, \quad \text{i.e.}
	\]
	\begin{equation}
	\label{eq:temp:cutin3}
\int_0^{\bar{T}} G(W_s) \d W^i_s = -\Big( \int_0^{T} +  \int_T^{\bar{T}}\Big) W^i_s \frac{W^1_s\d W^1_s+W^2_s\d W^2_s}{2\pi|W_s|^2}
-\frac{1}{2}\int_0^T \frac{W^i_s}{2\pi |W_s|^2} \d s.
	\end{equation}
	We now estimate the $L^k(\mathbb{P}_{T,a})$-norm of each of the three terms on the right-hand side of \eqref{eq:temp:cutin3}. 
	\medskip 
	
	By Burkholder inequality \eqref{eq:burkholder}, and since ${| W^i_s| |W^j_s|/|W_s|^2\leq 1}$ for all $s$, 
	\begin{align}
	\mathbb{E}\Big[  \Big( \int_0^T W^i_s \frac{W^1_s\d W^1_s+W^2_s\d W^2_s}{2\pi|W_s|^2} \Big)^{2k} \Big] 
&	\leq 2^{2k-1} \sum_{j\in \{1,2\}}
		\mathbb{E}\Big[  \Big( \int_0^T  \frac{W^i_s W^j_s\d W^j_s}{2\pi|W_s|^2} \Big)^{2k} \Big] 
\nonumber\\
&	\leq 
2^{2k-1} \sum_{j\in \{1,2\}} \int_0^T
\mathbb{E}\Big[  \Big(   \frac{W^i_s W^j_s}{2\pi|W_s|^2} \Big)^{2k} \Big] \d t
\nonumber\\
&\leq 
	2^{2k}
	(2k)^k T^{k-1} \int_0^T  \big(  \frac{1}{2\pi} \big)^{2k} \d t 
	=2^k \pi^{-2k} k^{k}  T^{k}.
	\label{eq:temp:expo1}
	\end{align}
	\medskip 
	
Recall $W_{[T,\bar{T}]}$ is a straight line segment between $W_T$ and $W_0$. Since 	
	\[\int_T^{\bar{T}} W^i_s \frac{W^1_s\d W^1_s+W^2_s\d W^2_s}{|W_s|^2}\leq 2  
	\int_T^{\bar{T}} |\d W_s|=2 |W_T-W_0|, 
	\]
	the second term in \eqref{eq:temp:cutin3} is bounded as follows: 
	\begin{equation}		\label{eq:temp:expo2}
	\mathbb{E}\Big[ \big(\int_T^{\bar{T}} \hspace{-0.2cm}W^i_s \frac{W^1_s\d W^1_s \hspace{-0.05cm}+ \hspace{-0.05cm}W^2_s\d W^2_s}{2\pi |W_s|^2}\big)^{2k}\Big]
	\leq 
	\pi^{-2k} \mathbb{E}[|W_T \hspace{-0.05cm}- \hspace{-0.05cm}W_0|^{2k} ]=2^{k}\pi^{-2k}  k!T^{k} \leq 2^k \pi^{-2k} k^{k}  T^{k}.
	\end{equation}
	\medskip 	
	
	To estimate the last term in  \eqref{eq:temp:cutin3}, we first use Hardy-Littlewood rearrangement inequality to obtain, for all $c\in \mathbb{R}^2$ and all $t>0$, 
	\[ 
	\int_{\mathbb{R}^2} \frac{p_t(c,x)}{2\pi |x|}\d x\leq 	\int_{\mathbb{R}^2} \frac{p_t(0,x)}{2 \pi |x|}\d x= \frac{1}{2\sqrt{2\pi t}}.
	\]
	Applying iteratively, we deduce that for all $x_0\in \mathbb{R}^2$ and for any $0=t_0<t_1<\dots< t_n<+\infty$,
	\[ 
	\int_{(\mathbb{R}^2)^{2k}} \prod_{j=1}^{2k} \frac{p_{t_j-t_{j-1}}(x_{j-1},x_j) |x_j^i|}{2\pi |x_j|^2}\d x_1\dots \d x_i\leq 
	\frac{1}{2^{2k}(2\pi)^{k}} \prod_{i=1}^{2k} \frac{1}{\sqrt{t_i-t_{i-1} }}.
	\]
	Since $2k$ is an integer, we deduce 
	\begin{align*}
	\mathbb{E}\Big[ \Big( \int_0^T \frac{|W^i_t|}{2 \pi  |W_t|^2} \d t \Big)^{2k} \Big]	
	&= (2k)! \mathbb{E}\Big[ \int_{0<t_1<\dots<t_{2k}<T} 	  \prod_{j=1}^{2k} \frac{|W^i_{t_j}|}{2 \pi |W_{t_j}|^2} \d t_1\dots \d t_j \Big]\\	
	&=(2k)! \int_{\mathbb{T}^{2k}_T} 	\int_{(\mathbb{R}^2)^{2k}} \prod_{j=1}^{2k} \frac{p_{t_j-t_{j-1}}(x_{j-1},x_j) x_j^i}{2 \pi  |x_j|^2}\d x_1\dots \d x_i\d t_1\dots \d t_j\\
	&\leq  ({2k})! 	\frac{1}{2^{2k}(2\pi)^{k}}  \int_{\mathbb{T}^{2k}_T} \prod_{i=1}^{2k} \frac{1}{\sqrt{t_i-t_{i-1} }} \d t_1\dots \d t_j.
	\end{align*}
	Using Lemma \ref{le:integralT}, we obtain
	\begin{align}
	\mathbb{E}\Big[ \Big( \frac{1 }{2}\int_0^T \frac{W^i_s}{2\pi |W_s|^2} \d s \Big)^{2k} \Big]
	&\leq k! 	\frac{1}{4^{2k}(2\pi)^{k}}  \int_{\mathbb{T}^{2k}_T}     
	\prod_{i=1}^{2k} (t_i-t_{i-1})^{-\frac12} \d t_1\dots \d t_{2k}\nonumber \\
	&\leq C  \frac{ e^{k}}{4^{2k}}  T^{k}   (2k)^{k}\ll 2^k \pi^{-2k} k^{k}  T^{k},		\label{eq:temp:expo3}	
	\end{align}
	for a universal constant $C<\infty$. 
	
	Using the H\"older inequality $|a+b+c|^{2k}\leq 3^{2k-1}(|a|^{2k}+|b|^{2k}+|c|^{2k})$ inside \eqref{eq:temp:cutin3}, and using the bounds \eqref{eq:temp:expo1}, \eqref{eq:temp:expo2} and \eqref{eq:temp:expo3}, we get 
	\begin{align*} 
	\mathbb{E}\Big[ \Big|  	\int_0^{\bar{T}} G(W_s) \d W^i_s\Big|^{2k}\Big] 
	&\leq 3^{2k-1} \Big( 
	\mathbb{E}\Big[\Big| 
	\int_0^T W^i_s \frac{W^1_s\d W^1_s+W^2_s\d W^2_s}{2\pi |W_s|^2}\Big|^{2k}\Big]\\
	& +
	\mathbb{E}\Big[\Big| 
	\int_T^{\bar{T}} W^i_s \frac{W^1_s\d W^1_s+W^2_s\d W^2_s}{2\pi|W_s|^2}\Big|^{2k}\Big]
	+\mathbb{E}\big[\Big| \frac{1}{2}\int_0^T \frac{W^i_s}{2\pi|W_s|^2} \d s\Big|^{2k}\Big]\Big)\\
	&\leq 	 
	3^{2k-1}(  2 \cdot 2^k \pi^{-2k} k^{k}  T^{k}  + C  \frac{\ e^{k}}{4^{2k}}  T^{k} (2k)^{k}  )\\
	&\leq 	 
	C'18^{k} \pi^{-2k} k^{k}  T^{k},
	\end{align*}
	for a universal constant $C'<\infty$.	
\end{proof}
\begin{remark}
	If we were to consider the integral up to $T$ rather than $\bar{T}$, we would end up with an extra term 
	$\mathbb{E}[ |W_T \log|W_T| |^k]$, which grows more quickly than $(Ck)^{k/2}$ for any constant $C$, as $k\to \infty$. In other words, we would have an extra ultraviolet divergence. With such a modification, Lemma \ref{le:finally} would become false, but we could recover it making a second modification: by pinning the endpoint of the Brownian motion (i.e. by looking at Brownian bridges instead), or by killing it as it exits a given compact set, or by conditionning it to stay in a given compact set. In the first case (pinning), we would get a locally uniform bound (in the endpoints of the Brownian bridge), but not a uniform one. This ultraviolet divergence is related to some of the difficulties we will encounter in subsection~\ref{sec:average2}.
\end{remark}

\begin{corollary} 
	\label{coro:finally} 
	There exists a constant $C$ such that for all $T,T'>0$, for all $a,b\in \mathbb{R}^2$, for all integer $k \geq 2$, 
	\[ \mathbb{E}_{T,T',a,b}\big[  (\mathcal{B}^\mathrm{It\hat{o}}_{W,W'})^{2k} \big]\leq C  (6\sqrt{2}/\pi)^{2k}    T^k T'^k  k^{2k},\]
	where $W,W'$ are two independent Brownian motions with durations $T,T'$ and starting points~$a,b$.
\end{corollary} 
\begin{proof}
	We apply Burkholder's inequality \eqref{eq:burkholder}, with $Y_s= \int_0^{\bar{T}'} G(W_s-W'_t) \d W'_t$, with the Brownian motion $W$, and with the filtration given by $\mathcal{F}_s=\sigma( (W_u)_{u\leq s}, (W'_v)_{v\in [0,T']} )$. We obtain 
	\begin{align*}
	\mathbb{E}\Big[ \Big|\int_0^{T }  \int_0^{\bar{T}' }  G(W_s-W'_t)  \d W'_t \d W_s\Big|^{2k} \Big]
	\leq& 
	(2k)^{k} T^{k-1 } \int_0^{T}  \mathbb{E}\Big[  \Big| \int_0^{\bar{T}'} G(W_s-W'_t) \d W'_t\Big|^{2k} \Big] \d s.
	\end{align*}	
	For Lebesgue-almost all $s \in [0,\bar{T}]$, $W_s\neq b$. We get 
	\[ 
	\mathbb{E}\Big[  \Big|  \int_0^{\bar{T}'} G(W_s-W'_t) \d W'_t\Big|^{2k} \Big]
	\leq \sup_{b,c \in \mathbb{R}^2 , b\neq c } 	  \mathbb{E}_b\Big[  \Big|   \int_0^{\bar{T}'} G(c-W'_t) \d W'_t\Big|^{2k} \Big],
	\]
	With Lemma \ref{le:finally}, we obtain 
	\[ 
	\mathbb{E}\Big[ \Big|\int_0^{T }  \int_0^{\bar{T}' }  G(W_s-W'_t)  \d W'_t \d W_s\Big|^{2k} \Big] \leq  C' 2^k (18/\pi^2)^{k}    T^k T'^k  k^{2k}= C' (6/\pi)^{2k}    T^k T'^k  k^{2k} . 
	\]
	Now we need to show a similar bound for 
	\[ 
	\mathbb{E}\Big[ \Big|\int_T^{\bar{T }}  \int_0^{\bar{T}' }  G(W_s-W'_t)  \d W'_t \d W_s\Big|^{2k} \Big]=  	\mathbb{E}\Big[ \Big|\int_{[W_T, W_0]}   \int_0^{\bar{T}' }  G(z-W'_t)  \d W'_t \d z\Big|^{2k} \Big].
	\]
	Using first the density of $W_T$, then applying $|\int_a^w f(z)\d z|^{2k}\leq |w-a|^{2k} \sup f(z)^{2k}$, then using  both the fact $\int_{\mathbb{R}^2} p_T(a,w) |w-a|^{2k} \d w= k! 2^kt^k$ and Lemma \ref{le:finally}, we get 
	\begin{align*}
	\mathbb{E}&\Big[ \Big|\int_{[W_T, W_0]}   \int_0^{\bar{T}' }  G(z-W'_t)  \d W'_t \d z\Big|^{2k} \Big]
	=
	\int_{\mathbb{R}^2} \frac{p_T(a,w)}{(2\pi)^{2k}} \mathbb{E}_b\Big[  \Big|\int_{[w, a]}   \int_0^{\bar{T}' }  \log|z-W'_t|  \d W'_t \d z\Big|^{2k}   \Big] \d w\\
	&\leq \int_{\mathbb{R}^2} \frac{p_T(a,w)}{(2\pi)^{2k}} |w-a|^{2k} \d w \sup_z  \mathbb{E}_b\Big[  \Big|   \int_0^{\bar{T}' }  \log|z-W'_t|  \d W'_t \Big|^{2k}   \Big] \\		
	&\leq \frac{k! 2^k T^k }{(2\pi)^{2k}} C 16^k T'^k k^k= k^{2k} \frac{8^{k+o(1)}}{\pi^{2k}e^k }
	\ll   (6/\pi)^{2k} T^k T'^k k^{2k}.
	\end{align*}
Using the Hölder's inequality $(a+b)^k \leq 2^{k-1} (a^k+b^k)$, we deduce 
\[
\mathbb{E}_{T,T',a,b}\big[  (\mathcal{B}^{\mathrm{It\hat{o}}}_{W,W'})^{2k} \big]
\leq 2^{k-1} C'' (6/\pi)^{2k}    T^k T'^k  k^{2k} 
\]
\end{proof}

\begin{lemma} 
	\label{le:expoLocalTime}
	Let $W,W'$ be independent standard planar Brownian motions with respective durations $T,T'$ and with starting points $a,b$ under $\mathbb{P}_{a,b,T,T'}$. Let $I_{W.W'}$ be their intersection local time. Let $\beta\in \mathbb{R}$ with $|\beta|<2$. Then, 
	\[ 
	\sup_{a,b,T<T'} \mathbb{E}_{a,b,T,T'}[ \exp( \frac{\beta}{\sqrt{TT'}} I_{W,W'}  )  ]<\infty.
	\]	
\end{lemma} 
\begin{proof}
	In \cite[Lemma 2]{LeGall}, it is proven that for $T=T'=1$ and $a=b$, 
	\[ 
	\mathbb{E}[ I_{W,W'}^p ]\leq C 2^{-p} p!
	\]	 
	 for a universal constant $p$, which gives the finiteness of exponential moment up to order $2$ (remark the exponential moment  are also known to be finite for all $\beta<0$). 
	The general bound 
	\[ 
	\mathbb{E}_{a,b,T,T'}[ I_{W,W'}^p ]\leq C 2^{-p} p! (TT')^{p/2}
	\]	 
	 for the same constant $C$, is a mere adaptation of the proof, which we sketch without entering into details. The extension to general $T,T'$ is direct: the proof of \cite[Lemma 2]{LeGall} extends directly to $T=T'$, possibly different from $1$, and furthermore  
	 uses the Cauchy--Schwarz inequality (\cite[p.175]{LeGall}) in a way that effectively decouples the $T$ and $T'$ variables. The extension to $a\neq b$ is more subtle, but follows from Hardy-Littlewood inequality applied at an early stage of the proof, showing in fact that  $\mathbb{E}_{a,b,T,T'}[ I_{W,W'}^p ]$ is maximal for $a=b$. 
	 
	 Remark that using monotonicity of the intersection local time and a stopping time argument, which might look natural, does not in fact work, as the Brownian motion $W$ started from the first time it hits $W'$ is not independent from $W'$.
\end{proof}

We deduce the existence of exponential moments for the Amperean area between two independent Brownian motions.
\begin{theorem}
	\label{th:mainY}
	Assume $|\beta|<\frac{8\pi}{8(6\sqrt{2})+\pi}\simeq 0.35$. Then, 
	\[ \sup_{a,a'\in \mathbb{R}^2,\ T,T'>0} \mathbb{E}_{T,T',a,a'}\big[  \exp\big( \frac{\beta}{\sqrt{TT'}} \mathcal{B}_{W,W'}\big) \big]<\infty,\]
	Furthermore, 
	$\exp\big( \frac{\beta}{\sqrt{TT'}} \mathcal{B}_{W,W'}^{\epsilon,\epsilon' }\big) $ converges toward 
	$\exp\big( \frac{\beta}{\sqrt{TT'}} \mathcal{B}_{W,W'} \big)$ in $L^1(\mathbb{P}_{T,T',a,a'})$  as $\epsilon,\epsilon'\to 0$.
	
	The same holds when $\mathcal{B}_{W,W'}$ and $\mathcal{B}^{\epsilon,\epsilon'}_{W,W'}$  are replaced with $\mathcal{B}_{W,W'}^{\mathrm{It\hat{o}}}$ and $\mathcal{B}^{\epsilon,\epsilon',\mathrm{It\hat{o}}}_{W,W'}$, under the weaker condition $|\beta|<\frac{\pi}{6\sqrt{2}}\simeq 0.37$.
\end{theorem}
\begin{proof} 
	The bound on $\exp\big( \frac{\beta}{\sqrt{TT'}} \mathcal{B}_{W,W'}^{\mathrm{It\hat{o}}  } \big)$ follows directly from Corollary \ref{coro:finally}. The $L^1$ convergence of 
	$\exp\big( \frac{\beta}{\sqrt{TT'}} \mathcal{B}_{W,W'}^{\epsilon,\epsilon',\mathrm{It\hat{o}} }\big) $ can then be deduced from Corollary \ref{coro:unifBound}. 
	
	For the bound on 
	$\exp\big( \frac{\beta}{\sqrt{TT'}} \mathcal{B}_{W,W'} \big)$, the condition on $\beta$ ensures that there exists $p,q\in(1,\infty)$ such that $1/p+1/q=1$, $q\beta<8$, $p\beta< \frac{\pi}{6\sqrt{2}}$. This is due to
	\[ 
\big(	1/p+1/q=1,\ q \beta=8,\ p\beta=\frac{\pi}{6\sqrt{2}}\big) \implies \beta=\frac{8\pi}{8(6\sqrt{2})+\pi}.
	\]
	
	Combining the It\^o-to-Stratonovich correction formula given by Proposition \ref{prop:cvY} with $L^p-L^q$ inequality, we get 
	\[ 
	\mathbb{E}[\exp\big( \frac{\beta}{\sqrt{TT'}} \mathcal{B}_{W,W'} \big) ]\leq 
		\mathbb{E}[ \exp\big( \frac{p \beta}{\sqrt{TT'}} \mathcal{B}^{\mathrm{It\hat{o}}  }_{W,W'} \big)]^{\frac{1}{p}} \mathbb{E}[ \exp ( - \frac{q\beta}{4 {\sqrt{TT'}} } I_{W,W'}) ]^{\frac{1}{q}}, 
	\]
	which is finite and bounded over $a,a',T,T'$ by Lemma \ref{le:expoLocalTime} and Corollary \ref{coro:finally}. The $L^1$ convergence similarly follows, using  Corollary \ref{coro:unifBound} again.
\end{proof}

\section{Convergence for the recentred variable \texorpdfstring{$X$}{X}}
\label{sec:X}
In this section, we study the random variables
\[ 
\mathcal{A}^\epsilon_W\coloneqq \int_{\mathbb{R}^2}  {\n}^\epsilon_W(z)^2 \d z \qquad \mbox{and} \qquad \tilde{\mathcal{A}}^\epsilon_W \coloneqq \mathcal{A}^\epsilon_W- \mathbb{E}[\mathcal{A}^\epsilon_W].
\]
For the whole section, we assume $W$ has duration $1$ and starts from $0$ (without loss of generality, because of scale and translation invariance) and 
$\epsilon=2^{-n/2}$ for some $n\in \mathbb{N}$, and we consider the limit as $n\to \infty$.
We will show that, along this sequence, $\tilde{\mathcal{A}}^\epsilon_W$ converges in $L^2$. The choice of this subsequence of $\epsilon$ is only for the clarity of the presentation. 
The general case with $\epsilon$ varying continuously can be similarly obtained, simply by writing $\epsilon$ as $\epsilon=\sqrt{s} 2^{-n/2}$ with $s\in[1,2)$ and $n\in \mathbb{N}$, and by showing that the estimations we prove here hold uniformly over $s\in[1,2]$, which is very straightforward but would divert our attention from the relevant ideas.

To prove this convergence, we will decompose $
\mathcal{A}^\epsilon_W$ in a way reminiscent of Paul Lévy's construction of the Lévy area (without recentering term) or Varadhan's construction of the renormalised self-intersection local time (with a recentering term)---except because of the \emph{boundary terms}, i.e. the terms coming from fact we have to turn open paths into loops by adding line segments, there is more cross-terms involved, compared to the Lévy area or the self-intersection local time. 

For $\gamma:[0,T]\to \mathbb{R}^2$, 
which is not assumed to be a loop, recall the winding function $\n_\gamma$ is defined as the winding function of the loop obtained by concatenation of $\gamma$ with a straight-line segment between its endpoints. 
 For $\gamma':[0,T']\to \mathbb{R}^2$ a second path such that $\gamma'_0=\gamma_T$, the winding function of the concatenation $\gamma\cdot \gamma'$ satisfies (almost everywhere, provided the range of $\gamma\cdot \gamma'$ has a vanishing Lebesgue measure on $\mathbb{R}^2$)
 \[ 
 \n_{\gamma\cdot \gamma'}=\n_\gamma+\n_\gamma'+ \n_{T_{\gamma_0,\gamma_T=\gamma'_0,\gamma'_{T'}}},
 \]
 where $T_{a,b,c}$ is the oriented triangle between the points $a,b,c$: considered as a loop, $T_{a,b,c}$ starts from $a$, goes to $b$ and then to $c$, before going back to $a$. 
 It is this triangle which is responsible for the extra terms mentionned above in the iterative decomposition of $\mathcal{A}^\epsilon_W$: these are the terms $T_{k,j}$ and $Z^\epsilon_{k,j}$ we will soon define.

We define $\mathbb{M}$ the set of couples $(k,j)$, where $k$ is a non-negative integer, $j$ is a positive integer and $j\leq 2^k$. 
For such a couple, we set $W_{(k,j)}$ the restriction of $W$ to $[2^{-k}(j-1), 2^{-k}j]$.
 We will write $(k,j)_-$ for $(k+1,2j-1)$, and $(k,j)_+$ for $(k+1,2j)$. Think of $W_{(k,j)_-}$ and $W_{(k,j)_+}$ as ``the first and second half of $W_{(k,j)}"$: indeed, for all $(k,j)\in \mathbb{M}$, $W_{(k,j)}$ is equal to the concatenation of 
 $W_{(k,j)_- }$ and $W_{(k,j)_+}$ (see the time segment $[2^{-k}(j-1), 2^{-k}j]$ parameterised by some $(k,j)\in \mathbb{M}$ on Figure~\ref{fig:kjs}).
 \begin{figure}[h]
 	\includegraphics{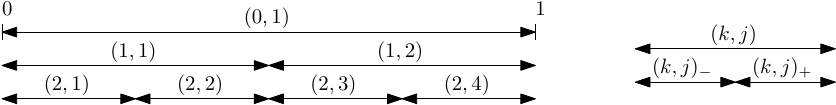}
 	\caption{\label{fig:kjs} Some of the time segments parameterised by the couples $(k,j)\in \mathbb{M}$. } 		
 \end{figure}

 Remark for all $(k,j)\in \mathbb{M}$, the Brownian motion $W_{(k,j)}$ is independent from $(W_{(k,l)})_{l<j}$ conditionally on the starting point of $W_{(k,j)}$. In particular, the equivalence class of $W_{(k,j)}$ up to translation is independent from $(W_{(k,l)})_{l<j}$.

Define $\delta_{k,j}$ be the winding function of the triangle delimited by $W_{2^{-k}(j-1)}$,  $W_{2^{-k-1}(2j-1)}$ and $ W_{2^{-k}j }$ (i.e. the starting point, midpoint and endpoint of $W_{(k,j)}$). It takes value either in $\{0,1\}$ or in $\{0,-1\}$. Let also $\delta_{k,j}^\epsilon=\delta_{k,j}\conv \varphi^\epsilon$. 
We use it to define the following random variables.
\begin{align*} 
\mathcal{A}^\epsilon_{k,j}& \coloneqq \int_{\mathbb{R}^2} {\n}^\epsilon_{W_{(k,j)}}(z)^2 \d z,  \quad	
\mathcal{B}^\epsilon_{k,j} \coloneqq  \int_{\mathbb{R}^2} {\n}^\epsilon_{W_{(k,j)_-}}(z)
{\n}^\epsilon_{W_{(k,j)_+}}(z) \d z , \\
\mathcal T^\epsilon_{k,j} 	& \coloneqq  \int_{\mathbb{R}^2} \delta^\epsilon_{k,j}(z)^2 \d z, \quad
\mathcal{C}^{\pm,\epsilon}_{k,j} 	\coloneqq  \int_{\mathbb{R}^2}\delta^\epsilon_{k,j}(z)  {\n}^\epsilon_{W_{(k,j)_\pm}}(z) \d z , \quad
  \mathcal{C}^\epsilon_{k,j} 	 \coloneqq \mathcal{C}^{-,\epsilon}_{k,j}+\mathcal{C}^{+,\epsilon}_{k,j} +\mathcal  T^\epsilon_{k,j}.
\end{align*} 
Remark in particular, $\mathcal{A}^\epsilon_{0,1}=\mathcal{A}^\epsilon_W$, $\mathcal{B}^\epsilon_{k,j}= \mathcal{B}^{\epsilon,\epsilon}_{W_{(k,j)_-}, W_{(k,j)_+}}$.
We will show that for all $(k,j)\in \mathbb{M}$, the random variables  $\mathcal{B}^\epsilon_{k,j}$, $\mathcal T^\epsilon_{k,j}$ and $\mathcal{C}^{\pm, \epsilon}_{k,j}$ (hence also $\mathcal{C}^\epsilon_{k,j}$) converge, in the $L^2$ sense,  
as $\epsilon\to 0$: we then write  
$\mathcal{B}_{k,j}$, $\mathcal T_{k,j}$, $\mathcal{C}^\pm_{k,j}$, and $\mathcal{C}_{k,j}$ for the limits. These limits can all be interpret, informally, as areas (these interpretations are rigorous when $W$ is replaced with a smooth loop).
For any such variable $\mathcal D_{k,j}$, $\mathcal{D}\in \{\mathcal{A},\mathcal{B},\mathcal{C}, \mathcal{C}^{\pm}, \mathcal{T},\mathcal{A}^\epsilon,\mathcal{B}^\epsilon,\mathcal{C}^\epsilon,  \mathcal{C}^{\pm,\epsilon}, \mathcal{T}^\epsilon \}$, we also set $\tilde{\mathcal D}_{k,j}$ the recentred variable $\mathcal D_{k,j} - \mathbb{E}[\mathcal D_{k,j}]$. 


The point of these definitions is that the relation of concatenation $W_{(k,j)}=W_{(k,j)_-} \cdot W_{(k,j)_+}$, gives (by Chasles relation for stochastic integrals) the pointwise equality 
\[ 
{\n}^\epsilon_{W_{(k,j)}}(z)= {\n}^\epsilon_{W_{(k,j)_-}}(z)
+{\n}^\epsilon_{W_{ (k,j)_+ }}(z)+ \delta^\epsilon_{k,j}.
\]
Taking the square and integrating over $z$, we deduce the recurrence relation 
\[
\mathcal{A}^\epsilon_{k,j}= \mathcal{A}^\epsilon_{(k,j)_-}+  \mathcal{A}^\epsilon_{(k,j)_+}+ \mathcal T_{k,j}^\epsilon
+ 2 \mathcal{B}^\epsilon_{k,j} + 2\mathcal{C}^{+,\epsilon}_{k,j} + 2\mathcal{C}^{-,\epsilon}_{k,j}. 
\]
Iterating over $k$, we deduce that for any integer $m$, 
\begin{align}
\mathcal{A}^\epsilon_{0,1}&= \sum_{j=1}^{2^m}  \mathcal{A}^\epsilon_{m,j}
+ \sum_{k=0}^{m-1}  \sum_{j=1}^{2^k}\mathcal  T_{k,j}^\epsilon
+ 2 \sum_{k=0}^{m-1}  \sum_{j=1}^{2^k} \mathcal{B}^\epsilon_{k,j}
+ 2 \sum_{k=0}^{m-1}  \sum_{j=1}^{2^k}\mathcal{C}^{-,\epsilon}_{k,j}
+ 2 \sum_{k=0}^{m-1}  \sum_{j=1}^{2^k}\mathcal{C}^{+,\epsilon}_{k,j}
\nonumber\\
&= \sum_{j=1}^{2^m}  \mathcal{A}^\epsilon_{m,j}
 - \sum_{k=0}^{m-1}  \sum_{j=1}^{2^k}\mathcal  T_{k,j}^\epsilon
 + 2 \sum_{k=0}^{m-1}  \sum_{j=1}^{2^k} \mathcal{B}^\epsilon_{k,j}
 + 2 \sum_{k=0}^{m-1}  \sum_{j=1}^{2^k} \mathcal{C}^\epsilon_{k,j}.
	\label{eq:decompoX}
\end{align}
In the following, we apply this decomposition with $m=n$ (recall $n$ is the integer such that $\epsilon=2^{-n/2}$) in order to deduce the $L^2$ convergence of $\tilde{\mathcal{A}}^\epsilon_{(0,1)}$, as $n\to \infty$. A key point and the reason to make this choice for $m$ is that the terms $\mathcal{A}^\epsilon_{m,j}=\mathcal{A}^{2^{-n/2}}_{n,j}$ that appear in \eqref{eq:decompoX} are exactly scale invariant in distribution: $\mathcal{A}^{2^{-n/2}}_{n,j}\overset{(d)}{=} 2^{-n}\mathcal{A}^{1}_{0,1} $. 
This will allow us to prove, simply using the law of large numbers,  that the first sum in \eqref{eq:decompoX} converges toward a finite deterministic constant, as $\epsilon\to 0$. The three other sums, however, each contribute to the limit, even after recentering. As for the counterterm, the $\mathcal{B}^\epsilon_{k,j}$ do not contribute to it, but the other terms do. 

\begin{lemma}
	\label{le:inL2}
	For all $\epsilon>0$, the random variables $\mathcal D^\epsilon_{0,1}$, $\mathcal{D}\in\{\mathcal{A},\mathcal{B},\mathcal{C},\mathcal{C}^\pm,\mathcal{T}\}$, lies in $L^2(\mathbb{P})$. 
\end{lemma}
\begin{proof}
	Using various Cauchy--Shwarz inequalities, it suffices to prove the result for $\mathcal  D=\mathcal T$ and  
	$\mathcal D=\mathcal{A}$.
	
	$\diamond$ Since $\|\varphi^\epsilon\|_{L^1}=1$, Young's convolution identity gives 
	\begin{equation}
	\label{eq:Tboundunifepsilon}
	|\mathcal T^\epsilon_{0,1}|\leq \int_{\mathbb{R}^2} \delta_{0,1}(z)^2\d z\leq \pi \|W\|_\infty^2\in L^2(\mathbb{P}).
	\end{equation}

	$\diamond$ As for $\mathcal{A}^\epsilon_{0,1}$,
	remark first \[C\coloneqq \sup_z \mathbb{E}[ \n^\epsilon_W(z)^8 ]=\sup_z  \mathbb{E}[ (\int_0^{\bar{T}} \theta^\epsilon_{W_t-z} \circ \d W_t)^8]<\infty.\] This can be proved with the same strategy we already used several time: split the integral at the 
	time $t=T$, so that the integral from $0$ to $T$ is a Stratonovich integral along a Brownian motion and the integral from $T$ to $\bar{T}$ is a Riemann integral along a line segment. 
	To control the Stratonovich integral, use the Stratonovich-to-It\^o correction formula, and use Burkolder's inequality to estimate the It\^o integral. Conclude by using the fact that $\theta^\epsilon$ is bounded. To control the Riemann integral from $T$ to $\bar{T}$, use the fact $\|W\|_\infty\in L^8(\mathbb{P})$.

	Recall furthermore that for $|z|> \|W\|_\infty+\epsilon$,  $\n^\epsilon_W(z)^2=0$. 
	Thus, using Cauchy--Schwarz inequality, we get 
	\begin{align*}
	\mathbb{E}[ (\mathcal{A}^\epsilon_{0,1})^2]
	&=\int_{(\mathbb{R}^2)^2} \mathbb{E}[\mathbbm{1}_{|z|\leq \|W\|_\infty+\epsilon } \mathbbm{1}_{|w|\leq \|W\|_\infty+\epsilon }  \n^\epsilon_W(z)^2 \n^\epsilon_W(w)^2 ]\d z \d w\\
	&\leq 	\int_{(\mathbb{R}^2)^2} \mathbb{P}(  |z|\leq \|W\|_\infty+\epsilon)^{\frac{1}{4}}  \mathbb{P}(|w|\leq \|W\|_\infty+\epsilon   )^{\frac{1}{4}}  \mathbb{E}[\n^\epsilon_W(z)^8]^\frac{1}{4} \mathbb{E}[\n^\epsilon_W(w)^8]^\frac{1}{4}\d z \d w\\
	&= \Big(
	\int_{\mathbb{R}^2} \mathbb{P}(  |z|\leq \|W\|_\infty+\epsilon)^{\frac{1}{4}}   \mathbb{E}[\n^\epsilon_W(z)^8]^\frac{1}{4} \d z
	\Big)^2\\
	&\leq \sqrt{C} \Big(	\int_{\mathbb{R}^2} \mathbb{P}(  |z|\leq \|W\|_\infty+\epsilon)^{\frac{1}{4}}   \d z \Big)^2<\infty.
	\end{align*} 
\end{proof}
\begin{remark} 
	An adaptation of the same proof, with a quantitative Burkholder inequality, already shows that all the polynomial moments of  $\mathcal{A}^\epsilon_{0,1}$, and even the strech-exponential moments $\mathbb{E}[ \exp( \beta \sqrt{\mathcal{A}^\epsilon_{0,1}} )]$ with $|\beta|>0$ small enough, are finite. 
\end{remark}

\begin{lemma} 
	\label{le:cvsumX}
	Recall $\epsilon=2^{-n/2}$. There is a constant $C$ such that for all $n$ , 
	\begin{equation} 
	\label{eq:Acv}
	\mathbb{E}	\Big[  \Big( \sum_{j=1}^{2^{n}} \tilde{\mathcal{A}}^{\epsilon}_{n,j} \Big)^2 \Big] \leq C \epsilon^2.
	\end{equation}
\end{lemma} 
\begin{proof} 
	From Markov property and translation invariance of both the Brownian motion and the Lebesgue measure on the plane,  for any given $n\in \mathbb{N}$ and $\epsilon'>0$, the random variables $(\tilde{\mathcal{A}}^{\epsilon'}_{n,j})_j$ are independent and identically distributed. Furthermore, from the scaling properties the Brownian motion and of the family $(\varphi^{\epsilon})_{\epsilon}$, we deduce that for any $j$, and for $\epsilon=2^{-n/2}$,
	\[	\tilde{\mathcal{A}}^\epsilon_{n,j}\overset{(d)}{=} 2^{-n}  \tilde{\mathcal{A}}^1_{0,1}= 2^{-n}  ( \mathcal{A}^1_{0,1}-\mathbb{E}[\mathcal{A}^1_{0,1}]).
	\] 
	Furthermore, $C\coloneqq \mathbb{E}[(\tilde{\mathcal{A}}^1_{0,1})^2]<\infty$ by Lemma \ref{le:inL2}.
	It follows from these facts that for all $j\neq k$, 
	$\mathbb{E}[\tilde{\mathcal{A}}^\epsilon_{n,j} \tilde{\mathcal{A}}^\epsilon_{n,k} ]=0$ and  $\mathbb{E}[(\tilde{\mathcal{A}}^\epsilon_{n,j})^2]= 2^{-2n} C$.
	Thus, 
	\[ 
	\mathbb{E}	\Big[  \Big( \sum_{j=1}^{2^{n}} \tilde{\mathcal{A}}^{\epsilon}_{n,j} \Big)^2 \Big] \leq C 2^{-n}=C\epsilon^2.\]
\end{proof}

In order to estimate the contribution of the three other terms in \eqref{eq:decompoX}, we will follow the same strategy, using in each case the following general result. 
\begin{lemma} 
	\label{le:cvgen}
	Let $(X^{(\epsilon)}_{k,j})_{(k,j)\in \mathbb{M},\epsilon\in 2^{-\mathbb{N}}} $ be a collection of random variables which satisfies the properties:
	\begin{enumerate}
		\item For any fixed $k$ and $\epsilon>0$, the random variables $(X^{(\epsilon)}_{k,j})_{ j\in \{1,\dots, 2^k\}}$ are globally independent.
		\item For all $(k,j)\in \mathbb{M}$, $(X^{(\epsilon)}_{k,j})_{\epsilon}$ is equal in distribution to $2^{-k} (X^{(2^{k/2}\epsilon) }_{0,1})_{\epsilon}$. 
		\item For all $\epsilon$, $\mathbb{E}[X_{0,1}^{(\epsilon)}]=0$.
		\item As $\epsilon\to 0$, $X^{(\epsilon)}_{0,1}$ converges in $L^2(\mathbb{P})$ toward a random variable $X_{0,1}$.  
	\end{enumerate}
Then, the following holds. 
\begin{enumerate}
	\item[(a)] For all $(k,j)\in \mathbb{M}$, $X^{(\epsilon)}_{k,j}$ converges in $L^2(\mathbb{P})$. Let $X_{k,j}$ be the limit. 
	\item[(b)] For any fixed $k$, the random variables $(X_{k,j})_{ j\in \{1,\dots, 2^k\}}$ are globally independent.	
	\item[(c)] $X_{k,j}$ is equal in distribution to $2^{-k} X_{0,1}$.
	\item[(d)] For all $(k,j)\in \mathbb{M}$, $\mathbb{E}[X_{k,j}]=0$ 
	\item[(e)] The sequence $n\mapsto X_n\coloneqq \sum_{k=0}^n \sum_{j=1}^{2^k} X_{k,j}$ converges in $L^2$. Let $X$ be the limit. 
	\item[(f)] The sequence 
	$n\mapsto Y_n\coloneqq \sum_{k=0}^n \sum_{j=1}^{2^k} X^{(2^{-n/2})}_{k,j}$ also converges  toward $X$ in $L^2(\mathbb{P})$.
\end{enumerate}
\end{lemma}
\begin{proof}
	The first four points are immediate: \emph{(a)} follows from \emph{(2)} and \emph{(4)}; \emph{(b)} follows from \emph{(a)} and \emph{(1)}; \emph{(c)} follows from \emph{(a)} and \emph{(2)}; \emph{(d)} follows from \emph{(a)} and \emph{(3)}.
	
	For \emph{(e)}, the independence and centring properties \emph{(b)} and \emph{(d)} ensure that for $j_1\neq j_2$, 
	$\mathbb{E}[ X_{k,j_1}X_{k,j_2}]=0$. Using also the scaling property \emph{(c)}, we deduce  $\mathbb{E}[(\sum_{j=1}^{2^k} X_{k,j} )^2]
	=\sum_{j=1}^{2^k} \mathbb{E}[ (X_{k,j})^2]=2^{-k} \mathbb{E}[(X_{0,1})^2]$. Since 
	$ 
	\sum_{k=0}^\infty 2^{-k/2} \mathbb{E}[(X_{0,1})^2]^{\frac{1}{2}}<\infty$, 
	we deduce that the sequence $S_n$ is a Cauchy sequence in $L^2(\mathbb{P})$, hence converges in $L^2(\mathbb{P})$. 
	
	We now prove \emph{(f)}. For $m\in \mathbb{N}$, let $X_n^{(m)}\!=\!\sum_{k=0}^n \sum_{j=1}^{2^k} X^{(2^{-m/2})}_{k,j}$. In particular, ${Y_n=	X_n^{(n)}}$.
	Let $C\coloneqq \sup\{ \mathbb{E}[ (X^{(\epsilon)}_{0,1})^2]^{\frac{1}{2}}: \epsilon \in 2^{-\mathbb{N}}\}$. This supremum is finite: indeed, as $X^{(\epsilon)}_{0,1}$ is convergent in $L^2(\mathbb{P})$, there exists $\epsilon_0$ such that for all $\epsilon\leq\epsilon_0$,  
	$\mathbb{E}[ (X^{(\epsilon)}_{0,1})^2]<1$. Since there is finitely many $n\geq 0$ such that $2^{-n}\geq  \epsilon_0$, we deduce the finiteness of $C$.

	Let $m\leq n$. Using the triangle inequality in $L^2(\mathbb{P})$, then the independence,scaling and centring assumptions \emph{(1), (2),(3)}, we obtain  
	\begin{align*}
	\|X_n^{(n)}-X_m^{(n)}\|_{L^2(\mathbb{P})}
	&=		
	\mathbb{E}\Big[\Big(\sum_{k=m}^n \Big|\sum_{j=1}^{2^k} X_{k,j}^{(2^{-n/2})}\Big|\Big)^2\Big]^{\frac12} \leq 
	\sum_{k=m}^n
	\mathbb{E}\Big[\Big(\sum_{j=1}^{2^k} X_{k,j}^{(2^{-n/2})}\Big)^2\Big]^{\frac12}\\
	& =
	\sum_{k=m}^n\sum_{j=1}^{2^k}
	\mathbb{E}\big[ (X_{k,j}^{(2^{-n/2})})^2\big]^{\frac12}
	=
	\sum_{k=m}^n 2^{-k/2}
	\mathbb{E}\big[ (X_{0,1}^{(2^{(k-n)/2})})^2\big]^{\frac12}\leq 
	C 2^{-m/2+1}.
	\end{align*}
	Let $\epsilon>0$. Let $n_0$ be such that $C2^{-m/2+1}<\epsilon/3$, so that 
	for all $n_0\leq m\leq n$, 
	\[\|X_n^{(n)}-X_m^{(n)}\|_{L^2(\mathbb{P})}<\epsilon/3.\] 
	Let $n_1\geq n_0$ be such that for all $n_2\geq n_1$, $\|X_{n_2}-X\|_{L^2(\mathbb{P})}\leq \epsilon/3$, which exists by \emph{(e)}. 
	
Using \emph{(a)} and taking finite linear combinations, we deduce
	$X^{(n)}_{n_1}$ converges in $L^2(\mathbb{P})$ toward $X_{n_1}$ as $n\to \infty$. Thus, there exists 
	$n_2\geq n_1$ such that for all $n\geq n_2$, \[\|X^{(n)}_{n_2}-X_{n_2}\|_{L^2(\mathbb{P})}\leq \epsilon/3.\]
	Thus, for all $n\geq n_2$, 
	\[ 
	\|X_n^{(n)}-X\|_{L^2(\mathbb{P})}\leq \|X_n^{(n)}-X_{n_2}^{(n)}\|_{L^2(\mathbb{P})}+
	\|X^{(n)}_{n_2}-X_{n_2}\|_{L^2(\mathbb{P})}
	+\|X_{n_2}-X\|_{L^2(\mathbb{P})}\leq \epsilon,
	\] 
	which concludes the proof.
\end{proof}

We now proceed to show that the families of random variables
 $(\mathcal{B}_{k,j}^\epsilon)$, $(\tilde{\mathcal{C}}_{k,j}^\epsilon)$ and 
$(\tilde{\mathcal{T}}_{k,j}^\epsilon)$ each satisfies the four conditions in Lemma 
\ref{le:cvgen}. Conditions \emph{(1)} and \emph{(2)} are easily checked, and follow from the same arguments as in Lemma \ref{le:cvsumX}. Furthermore, the families    
$(\tilde{\mathcal{C}}_{k,j}^\epsilon)$ and 
$(\tilde{\mathcal{T}}_{k,j}^\epsilon)$ clearly satisfy the centring condition \emph{(3)} , since these have been defined by recentring. We now prove that $(\mathcal{B}_{k,j}^\epsilon)$ is already centred. 
\begin{lemma} 
	\label{le:sympassym}
	For all $\epsilon>0$, it holds $\mathbb{E}[ \mathcal{B}^\epsilon_{0,1}]=0$.	
\end{lemma}
\begin{remark}
	Although for any fixed $z$, we will show ${\n}^\epsilon_{W_{(0,1)_-}}(z){\n}^\epsilon_{W_{(0,1)_+}}(z)$ is symmetric in distribution, this does \emph{not} imply $ \mathcal{B}^\epsilon_{0,1}$ is symmetric in distribution (this symmetry in distribution does not hold for the joint distributions over different points $z$). 
\end{remark} 
\begin{proof}[Proof of Lemma \ref{le:sympassym}]
	For a given $\epsilon>0$, recall $\mathcal{A}^\epsilon_{(0,1)_{-}}$ and $\mathcal{A}^\epsilon_{(0,1)_{+}}$ both lie in $L^2(\mathbb{P})$ (Lemma \ref{le:inL2}), hence in $L^1(\mathbb{P})$. The Cauchy--Schwarz inequality in $L^2(\Omega\times \mathbb{R}^2)$ gives 
	\begin{align*} \mathbb{E}[ 
		\int_{\mathbb{R}^2} |{\n}^\epsilon_{W_{(0,1)_-}}(z)
		{\n}^\epsilon_{W_{(0,1)_+}}(z)| \d z] 
		&\leq 
		\mathbb{E}[ 
		\int_{\mathbb{R}^2} {\n}^\epsilon_{W_{(0,1)_-}}(z)^2
		\d z]^\frac12 
		\mathbb{E}[ 
		\int_{\mathbb{R}^2} 
		{\n}^\epsilon_{W_{(0,1)_+}}(z)^2 \d z]^\frac12\\
		&= \mathbb{E}[\mathcal{A}^\epsilon_{(0,1)_{+}} ] 
		\mathbb{E}[\mathcal{A}^\epsilon_{(0,1)_{-}} ]<\infty.
	\end{align*}
	Since the random variable  $\int_{\mathbb{R}^2} |{\n}^\epsilon_{W_{(0,1)_-}}(z)
	{\n}^\epsilon_{W_{(0,1)_+}}(z)| \d z$ is 
	in $L^1(\Omega)$, we can now apply Fubini-Tonelli's theorem:
	\[ 
	\mathbb{E} [ \mathcal{B}^\epsilon_{0,1}]=
		\mathbb{E} [ \int_{\mathbb{R}^2} f(z) \d z ]=
	\int_{\mathbb{R}^2} \mathbb{E}[ f(z)] \d z ,
	\quad \text{where}
	\quad f(z)\coloneqq {\n}^\epsilon_{W_{(0,1)_-}}(z){\n}^\epsilon_{W_{(0,1)_+}}(z).
	\]	
	Thus, in order to prove the lemma, it suffices to show that for all $z\in \mathbb{R}^2$, $f(z)\overset{(d)}=-f(z)$.
	 Recall $W_{(0,1)^-}=W_{1,1}$ is the restriction of $W$ to the time interval $[0,1/2]$, while $W_{(0,1)^+}=W_{1,2}$ is the restriction of $W$ to the time interval $[1/2,1]$.
	Fix $z\in \mathbb{R}^2$, and let $W^\dagger_{1/2}$ be the reflection of the Brownian motion $W_{(1,2)}$ with respect to the axis $(W_{1/2},z)$. By reflection principle,  $(W_{(1,1)},W_{(1,2)})\overset{(d)}=(W_{(1,1)},W_{(1,2)}^\dagger)$. Remark ${\n}^\epsilon_{W_{(1,2)}^\dagger }(z)= -{\n}^\epsilon_{W_{(1,2)}}(z)$, which is due to the fact that the point $z$ lies in the axis of reflection and that the symmetry under reflections of the mollifier $\phi^\epsilon$.
	We deduce 
	\[f(z)={\n}^\epsilon_{W_{(1,1)} }(z)	{\n}^\epsilon_{W_{(1,2)} }(z)\overset{(d)}= 
	{\n}^\epsilon_{W_{(1,1)} }(z)	{\n}^\epsilon_{W^\dagger_{(1,2)} }(z)
	=-f(z),
	\]	
	hence $\mathbb{E}[f(z)]=0$, which concludes the proof. 
\end{proof}

We now prove that the families of random variables $(\mathcal{B}_{k,j}^\epsilon)$, $(\tilde{\mathcal{C}}_{k,j}^\epsilon)$ and 
$(\tilde{\mathcal{T}}_{k,j}^\epsilon)$ satisfy the condition \emph{(d)} in Lemma \ref{le:cvgen}, i.e. that $\mathcal{B}_{0,1}^\epsilon$, $\tilde{\mathcal{C}}_{0,1}^\epsilon$ and 
$\tilde{\mathcal{T}}_{0,1}^\epsilon$ converge in $L^2(\mathbb{P})$ as $\epsilon\to 0$. 
In the process, we identify the limits, and give a bound on the convergence speed which we will not use in this section but in  Section~\ref{sec:average}, when we will estimate the counterterm $\mathbb{E}[\mathcal{A}^\epsilon_W]$. 

Let us remind that the convergence of $\mathcal{B}_{0,1}^\epsilon$ is essentially the content of Sections \ref{sec:epsilon} and \ref{sec:zero}, and follows from Proposition \ref{prop:cvY}.\footnote{The Brownian motions $W_{(0,1)_-}$ and $W_{(0,1)_+}$ are not independent from each other so we cannot directly apply Proposition  \ref{prop:cvY}: we must consider instead the Brownian motions $W_{(0,1)_-}':t\in[0,1/2]\mapsto W_{1/2-t}-W_{1/2}$ and $W_{(0,1)_+}':t\in[0,1/2]\mapsto W_{t+1/2}-W_{1/2} $. These Brownian motions are independent, so we can apply 
	Proposition  \ref{prop:cvY}. We then use the symmetry of the Amperean area under simultaneous translation of both curves and its antisymmetry under reversal of orientation. 
} The limit is the Amperean area between $W_{(0,1)_-}=W_{[0,1/2]}$ and $W_{(0,1)_+}=W_{[1/2,1]}$.
\begin{lemma} 
	\label{le:Tcv}
	As $\epsilon \to 0$, $\mathcal{T}_{0,1}^\epsilon$ converges toward $\mathcal{T}_{0,1}=\int_{\mathbb{R}^2} \delta_{k,j}(z)^2 \d z $, both almost surely and in $L^2(\mathbb{P})$. Furthermore, there exists a constant $C$ such that for all $\epsilon\in(0,1]$,  $\mathbb{E}[(\mathcal{T}_{0,1}^\epsilon -\mathcal{T}_{0,1})^2]\leq C \epsilon^2$.  	
\end{lemma} 
\begin{proof}
	The almost sure convergence follows directly from dominated convergence theorem. 
	
	Let $T^\epsilon$  be the $\epsilon K_\varphi$-neighbourhood of the ``contour'' of the triangle $T_{0,1}$ with vertices $0, W_{1/2}, W_1$, i.e. the union of the three segments $[0,W_{1/2}]$, $[W_{1/2,W_1}]$, $[W_0,W_1]$. 
	The set  $T^\epsilon$ is equal to the union of three balls with radius $\epsilon K_\varphi$ and three rectangles with width $2\epsilon K_\varphi$ and respective lengths $|W_1|$, $|W_{1/2}|$, and $|W_1-W_{1/2}|\leq |W_1|+|W_{1/2}|$. It follows that, for a constant $C$ which depends on $\varphi$ but not on the realisation of $W$ nor on $\epsilon\in(0,2]$, the area of $T^\epsilon$ is less than $C\epsilon (\|W\|_\infty +2)$ (we allow $\epsilon\in[1,2]$ for a later usage in Lemma \ref{le:Zcv}).  
	
	Since  $|\delta^\epsilon_{0,1}-\delta_{0,1}|$ is equal to $0$ except on $T^\epsilon$ where it is bounded by $1$,
	$|\mathcal{T}_{0,1}-\mathcal{T}^\epsilon_{0,1}|$ is smaller than the area of $T^\epsilon$. Thus, 
	\[ 
	\mathbb{E}[   |\mathcal{T}_{0,1}-\mathcal{T}^\epsilon_{0,1}|^2  ]\leq C^2 \epsilon^2 \mathbb{E}[ (\|W\|_\infty +2)^2 ]=C'\epsilon^2,
	\]
	which concludes the proof.
\end{proof}

\begin{lemma} 
	\label{le:Zcv}
	Let 
	\[ 
	\mathcal{C}^-_{0,1}\coloneqq \int_0^{\frac{1}{2}} (\delta_{0,1} \conv \theta)_{W_t}  \d W_t+\int_{W_{1/2}}^{W_0}  (\delta_{0,1} \conv \theta)_u \d u,\quad
	\mathcal{C}^+_{0,1}\coloneqq \int_{\frac{1}{2}}^1 (\delta_{0,1} \conv \theta)_{W_t}  \d W_t+\int_{W_1}^{W_{1/2}}  (\delta_{0,1} \conv \theta)_u \d u,
	\]
	\[ 
	\mathcal{C}_{0,1}\coloneqq \mathcal{C}^-_{0,1}+\mathcal{C}^+_{0,1}+\mathcal{T}_{0,1}= \int_{0}^1 (\delta_{0,1} \conv \theta)_{W_t} \d W_t+\int_{W_1}^{W_{0}}  (\delta_{0,1} \conv \theta)_u \d u.
	\]	
	Then, there exists a constant $C$ such that for all $\epsilon\in(0,1]$, 
	\begin{equation}
	\label{eq:Zboundtemp}
	\mathbb{E}[ (\mathcal{C}_{0,1}^{\pm,\epsilon}-\mathcal{C}_{0,1}^{\pm})^2]\leq C \epsilon \quad \text{and} \quad 
		\mathbb{E}[ (\mathcal{C}_{0,1}^{\epsilon}-\mathcal{C}_{0,1})^2]\leq C \epsilon.
	\end{equation}
	In particular, the random variables $\mathcal{C}^{\pm,\epsilon}_{0,1}$ and $\mathcal{C}^{\epsilon}_{0,1}$  converge toward $\mathcal{C}^{\pm}_{0,1}$ and $\mathcal{C}_{0,1}$ in $L^2(\mathbb{P})$ as $\epsilon\to 0 $. 
\end{lemma} 
\begin{remark}
	Since $\div( \theta)=0$, $\div(\delta_{0,1}\conv \theta)=0$, and thus the integral defining $\mathcal{C}^\pm_{0,1}$ and $\mathcal{C}^{\pm,\epsilon}_{0,1}$ can be understood in the Ito or Stratonovich sense. 
	
	The random function $\delta_{0,1}$ is not adapted to the natural filtration of $W$, but this integral are well-defined by using the enlarged filtration $\mathcal{F}_t=\sigma ( (W_s)_{s\leq t}, W_{1/2}, W_1)$, with respect to which $W_{[0,1/2]}$ and $W_{[1/2,1]}$ are Brownian bridge hence semimartingales. See e.g. \cite{Nualart} for a much more general situation and equivalent definitions of stochastic integrals of non-adapted processes.
	
	With respect to this enlarged filtration, the process $W$ (say, restricted to $[0,1/2]$) is not a Brownian motion, nor even a local martingale, but a Brownian bridge from $0$ to the $\mathcal{F}_0$-measurable point $W_{1/2}$. It thus has a non-zero finite-variation part. Be careful that, as a consequence, the process $
	t\in[0,1/2]\mapsto  \int_{0}^t (\delta_{0,1} \conv \theta)_{W_s} \d W_s$
	is  \emph{not} a local martingale with respect to this extended filtration. In particular, $\mathbb{E}[\mathcal{C}_{0,1}^-]\neq 0$ although $\mathbb{E}[\int_{W_{1/2}}^{W_{0}}  (\delta_{0,1} \conv \theta)_u \d u]=0$.
\end{remark}
\begin{proof}[Proof of Lemma 	\ref{le:Zcv}]
	We first prove the bound  on $\mathcal{C}^{-,\epsilon}_{0,1}$ in two steps. We set 
	\[
	F_\epsilon(x)\coloneqq \varphi^\epsilon \conv \varphi^\epsilon \conv \delta_{0,1} \conv \theta(x) -\delta_{0,1} \conv \theta(x).
	\]
	Remark $F_\epsilon$ is a random vector field, measurable with respect to $\sigma(W_{1/2}, W_1)$. 
	
	\emph{step 1:} 
	We first prove that there exists a constant $C$ such that for all $\epsilon\in( 0,1]$, \[	\mathbb{E}[ \| F_\epsilon\|_{\infty}^4]\leq C \epsilon^2.\]
	Recall the definition of $T^\epsilon$ from Lemma \ref{le:Tcv}.
	Since $\Supp(\varphi\conv \varphi)\subset B_{2K_\varphi}$, for $v\in \Supp(\varphi\conv \varphi)$, 
	\begin{equation}
		\label{eq:temp:dumbbound}
		|\delta_{0,1}(\ \cdot \ -\epsilon v   )- \delta_{0,1}  |\leq  \mathbbm{1}_{T^{2\epsilon} }.
	\end{equation}  
	Let $r_\epsilon$ be such that the area of $B_{r_\epsilon}$ is equal to the area of $T^{2\epsilon}$. 
	Thus, the symmetric decreasing rearrangement of $w\mapsto \mathbbm{1}_{T^{2\epsilon}}(z-w)$ is $\mathbbm{1}_{B_{r_\epsilon}}$. 
As we proved in Lemma \ref{le:Tcv}, for $\epsilon\in(0,1]$, the area of $T^{2\epsilon}$ is less than $2 C\epsilon (\|W\|_\infty+2)$, where $C$ only depends on $\varphi$.  It follows that $r_\epsilon\leq C'  \sqrt{ \epsilon(\|W\|_\infty +1) }$ for a constant $C'$ which only depends on $\varphi$. 
	Using first $|\int f g|\leq \|f\|_{L^1}\|g\|_{L^\infty}$, then \eqref{eq:temp:dumbbound},  and then Hardy--Littlewood rearrangement inequality, we obtain
	\begin{align*}	 |F_\epsilon(x)|
	& =\Big| \int_{(\mathbb{R}^2)^2} (\varphi\conv \varphi)(v) \theta(w) (\delta_{0,1}(z-\epsilon v -w  )- \delta_{0,1}(z -w  ) )  \d w \d v\Big| \\
	&\leq \sup_{v \in \Supp(\varphi \conv \varphi) }
	\int_{\mathbb{R}^2} |\theta(w)| |\delta_{0,1}(z-\epsilon v -w  )- \delta_{0,1}(z -w  ) |  \d w \\
	&\leq 2 \sup_{v \in \Supp(\varphi \conv \varphi) }
	\int_{\mathbb{R}^2} |\theta(w)| \mathbbm{1}_{T^{2\epsilon}}(z-w)   \d w \\
	&\leq 2\int_{\mathbb{R}^2}  |\theta(w)| \mathbbm{1}_{B_{r_\epsilon}}(w)  \d w= 2 r_\epsilon\leq C'\sqrt{ \epsilon(\|W\|_\infty+1)}.
	\end{align*} 
	It follows that 	
	\[	\mathbb{E}[ \|F_\epsilon\|_{\infty}^4]
	\leq C'^4 \epsilon^2 \mathbb{E}[ (\|W \|_\infty+1)^2  ]	.\]
	
	\emph{step 2:}
	Using again Lemma \ref{le:derivBelowRP} to swap integrals, we obtain 
	\[ 
	\mathcal{C}^{-,\epsilon}_{0,1}= \int_0^{\frac{1}{2}} (\varphi^\epsilon \conv \varphi^\epsilon \conv \delta_{0,1} \conv \theta)_{W_t} \circ \d W_t+\int_{W_{1/2}}^{0}  (\varphi^\epsilon \conv\varphi^\epsilon \conv\delta_{0,1} \conv \theta)_u \d u,
	\]	
	thus 
	\[ 
	\mathcal{C}^{-,\epsilon}_{0,1}-	\mathcal{C}^{-}_{0,1}=  \int_0^{\frac{1}{2}} F_\epsilon(W_t) \d W_t +\int_{W_{1/2}}^{0} F_\epsilon(u) \d u, 
	\]
	and it suffices to prove there exists a constant $C$ such that for all $\epsilon\in (0,1]$, 
	\[ 
	\mathbb{E}\Big[ \Big( \int_0^{1/2} F_\epsilon(W_t)\d W_t \Big)^2 \Big] \leq C \epsilon  \quad \text{and} \quad 
	\mathbb{E}\Big[ \Big( \int_{W_{1/2 }}^0 F_{\epsilon}(u) \d u\Big)^2 \Big]\leq C \epsilon  .  
	\]
	For the first term, it is tempting to use Itô isometry directly, but this is not possible since $W$ is not a Brownian motion with respect to the filtration for which  $F_\epsilon(W_t)$ is adapted. Let thus, for $t\in[0,1/2]$, $B_t\coloneqq W_t-\int_0^t \frac{W_{1/2} -W_s}{1/2-s} \d s$, so that conditionnaly on $( W_{1/2}, W_1)$, $B$ is a Brownian motion. Let $\mathbb{E}'[\ \cdot\ ]=\mathbb{E}[\ \cdot\ | (W_{1/2}, W_1)]$, and remark $ \mathbb{E}'[f(F_\epsilon)]=f(F_\epsilon)$ for any deterministic function $f$, as $F_\epsilon$ is $\sigma(W_{1/2}, W_1)$-measurable. Using $\d W_t= \d B_t+ \frac{W_{1/2} -W_t}{1/2-t}\d t$, $(a+b)^2\leq 2(a^2+b^2)$, and then Itô isometry, we deduce 
	\begin{align*} 
	&	\mathbb{E}'\Big[ \Big( \int_0^{1/2} F_\epsilon(W_t)\d W_t \Big)^2\Big]\\
	& \leq 2	\mathbb{E}'\Big[ \Big( \int_0^{1/2} F_\epsilon(W_t)\d B_t \Big)^2\Big]
	+2\mathbb{E}'\Big[ \Big( \int_0^{1/2} |F_\epsilon(W_t)| \frac{|W_{1/2} -W_t|}{1/2-t}\d t \Big)^2 \Big]\\
	&= 
	2	 \int_0^{1/2} \mathbb{E}'[ |F_\epsilon(W_t)|^2     ] \d t
	+2\int_0^{1/2} \int_0^{1/2}  \mathbb{E}'[   |F_\epsilon(W_t)||F_\epsilon(W_s)| \frac{|W_{1/2} -W_t||W_{1/2} -W_s|}{(1/2-t)(1/2-s)}]\d t \d s \\
	&\leq	\mathbb{E}'[ \|F_\epsilon\|_\infty^2] +
2	 \int_0^{1/2}  \int_0^{1/2}  \frac{ \mathbb{E}'[\|F_\epsilon\|_\infty^2 |W_{1/2} -W_t|  |W_{1/2} -W_s| ]  }      { (1/2-t)(1/2-s) } \d t \d s\\	
&= \|F_\epsilon\|_\infty^2 +
2\|F_\epsilon\|_\infty^2 \int_0^{1/2} \int_0^{1/2}  \frac{  \mathbb{E}'[ |W_{1/2} -W_t|  |W_{1/2} -W_s| ]  }      { (1/2-t)(1/2-s) } \d t \d s\\
&\leq \|F_\epsilon\|_\infty^2 +
2\|F_\epsilon\|_\infty^2 \int_0^{1/2} \int_0^{1/2}  \frac{  \mathbb{E}'[ |W_{1/2} -W_t|^2]^{\frac12}  \mathbb{E}'[ |W_{1/2} -W_s|^2]^{\frac12} }      { (1/2-t)(1/2-s) } \d t \d s\\
	&= \|F_\epsilon\|_\infty^2 +
	2\|F_\epsilon\|_\infty^2 \Big( \int_0^{1/2}  \frac{  \mathbb{E}'[ |W_{u}|^2 ]^{\frac{1}{2}}}{ u } \d u  \Big)^2. 	
	\end{align*}
	For $u\in [0,1/2]$, conditionally on $W_{1/2}$, the random vector $W_{u}$ is Gaussian, centred at $2u W_{1/2}$ and with variance $2u(1/2-u)\leq u$ . Thus, 
	\[ 
	\mathbb{E}'[ |W_{u}|^2 ]\leq 4u^2 |W_{1/2}|^2 + u\leq u (1+ |W_{1/2}|^2), 
	\]
	from which we deduce
	\[\mathbb{E}'\Big[ \Big( \int_0^{1/2} \hspace{-0.3cm} F_\epsilon(W_t)\d W_t \Big)^2\Big]
	\leq 
	\|F_\epsilon\|_\infty^2 (1 + 2   ( \int_0^{1/2}  \frac{1 }{ \sqrt{u}} \d u  )^2(1+ |W_{1/2}|^2) )\leq 9 \|F_\epsilon\|_\infty^2(1+|W_{1/2}|^2 ).
	\]
	We take the expectation (i.e. we remove the conditioning on $(W_{1/2}, W_1)$) and apply Cauchy--Schwarz inequality. 
	Since $|W_{1/2}|$ has a finite fourth moment, we deduce that for a some constant $C'$, 
	\[ 
	\mathbb{E}\Big[ \Big( \int_0^{1/2} F_\epsilon(W_t)\d W_t \Big)^2\Big]
	\leq 9 \mathbb{E}[  (1+|W_{1/2}|^2)^2 ]^\frac12  \mathbb{E}[   \|F_\epsilon\|_\infty^4   ]^{\frac{1}{2}}\leq C' \epsilon,
	\]
using step 1.
	%
	%
	%
	Furthermore, 
	\begin{align*} 
	& \mathbb{E}\Big[ \Big(  
	\int_{W_{1/2}}^{0} F_\epsilon(u)  \d u\Big)^2 \Big]  \\
	&\leq \mathbb{E}[|W_{1/2}|^2 \| F_\epsilon \|_{\infty}^2] 
	\leq 
	\mathbb{E}[|W_{1/2}|^4]^{\frac{1}{2}} \mathbb{E}[\| F_\epsilon\|_{\infty}^4]^\frac{1}{2} 
	\leq C'' \epsilon,
	\end{align*}
	using step 1 again. This proves the bound 
\eqref{eq:Zboundtemp} for $\mathcal{C}^{-,\epsilon}_{0,1}$.  
	
	The case of $\mathcal{C}^{+,\epsilon}_{0,1}$ is similar. We then deduce that 
	$\mathcal{C}^\epsilon_{0,1}$ converges toward $ \mathcal{C}^{-}_{0,1}+\mathcal{C}^{+}_{0,1}+\mathcal{T}_{0,1}$, with the same asymptotic bound of order $\epsilon$, using Lemma \ref{le:Tcv}. The fact that the given value for $\mathcal{C}_{0,1}$ coincide with  $ \mathcal{C}^{-}_{0,1}+\mathcal{C}^{+}_{0,1}+\mathcal{T}_{0,1}$ follows from Green's formula (or equivalently here, from residue theorem) applied in the triangle $(0,W_{1/2},W_1)$: for $f$ bounded, measurable, and compactly supported, 
	\[ 
	\int_{W_0}^{W_{1/2}}  (f \conv \theta)_u \d u
	+\int_{W_{1/2}}^{W_1}  (f \conv \theta)_u \d u
	+\int_{W_1}^{W_0}  (f \conv \theta)_u \d u= 
	\int_{\mathbb{R}^2} f(z) \delta_{0,1}(z) \d z.
	\]	
	We use this here with $f= \delta_{0,1}$.
\end{proof} 
\begin{remark} 
	The term $\mathcal{C}_{k,j}$ should be interpreted, through a formal Green's formula, as equal to the ill-defined integral $\int_{\mathbb{R}^2} \delta_{k,j}(z) \n_{W_{(k,j)}}(z) \d z$. 
\end{remark} 

Since we have now proved that the families 
 $(\mathcal{B}_{k,j}^\epsilon)$, $(\tilde{\mathcal{C}}_{k,j}^\epsilon)$ and 
$(\tilde{\mathcal{T}}_{k,j}^\epsilon)$ each satisfy the assumptions of Lemma \ref{le:cvgen}, we deduce that the sums 
\[ 
\sum_{k=0}^{n-1} \sum_{j=1}^{2^k} {\mathcal{B}}^\epsilon_{k,j} , \qquad 
\sum_{k=0}^{n-1} \sum_{j=1}^{2^k} \tilde{\mathcal{C}}^\epsilon_{k,j}, \qquad 
\sum_{k=0}^{n-1} \sum_{j=1}^{2^k} \tilde{\mathcal{T}}^\epsilon_{k,j}, 
\]
with $\epsilon=2^{-n/2}$, are convergent as $n\to \infty$, and that the limits $\mathcal{B}$, $\tilde{\mathcal{C}}$ and $\tilde{\mathcal{T}}$  are respectively equal to the sum of the individual limits: for $\mathcal{D}\in \{ \mathcal{B}, \tilde{\mathcal{C}}, \tilde{\mathcal{T}}\}$,  
\begin{equation}
	\label{eq:sumsum}
{\mathcal{D}}=\sum_{k=0}^{\infty}  \sum_{j=1}^{2^k} \lim_{\epsilon\to 0}{\mathcal{D}}^\epsilon_{k,j}=\sum_{k=0}^{\infty}  \sum_{j=1}^{2^k} {\mathcal{D}}_{k,j}.
\end{equation}

We can now conclude: 
\begin{theorem} 
	\label{th:cvX} 
	As $\epsilon \in 2^{-\mathbb{N}}$ goes to $0$, $\tilde{\mathcal{A}}_W^\epsilon$ converges in $L^2(\mathbb{P})$.	
	The limit is given by 
	\[ 
	\sum_{k=0}^\infty \sum_{j=1}^{2^n} (\mathcal{B}_{k,j}+ \mathcal{C}_{k,j}-\mathbb{E}[\mathcal{C}_{k,j}]-\mathcal{T}_{k,j}+\mathbb{E}[\mathcal{T}_{k,j}]).\]
\end{theorem} 
\begin{proof} 
	This follows from the decomposition \eqref{eq:decompoX}, Lemma \ref{le:cvsumX}, and 
	\eqref{eq:sumsum}. 
\end{proof}

As a preparation for the next section in which we will no longer use the details of the definitions of the variables $\mathcal{D}_{k,j}$, let us conclude this section by showing that $\mathcal{C}_{0,1}^\pm$ admits some finite exponential moments. 
\begin{lemma} 
\label{le:Cexpmom}
	There exist $\beta_c>0$ such that for all $\beta\in(-\beta_c, \beta_c)$,
	\[ 
	\mathbb{E}[\exp(\beta \mathcal{C}_{0,1}^\pm)]<\infty. 
	\]
\end{lemma} 
\begin{proof}
	We only prove it for $\mathcal{C}_{0,1}^-$, it is similar for $\mathcal{C}_{0,1}^+$. 
	We use the notations $B$ and $\mathbb{E}'$ from the proof of Lemma~\ref{le:Zcv} ($\mathbb{E}'$ is the expectation conditional on $(W_{1/2},W_1)$, and conditionnally on $W_{1/2}$, the process $W_{[0,1/2]}$ is a Brownian bridge driven by the Brownian motion 
	$B$). Set $\rho$ such that $\pi\rho^2$ is equal to the area of the convex set delimited by the three points $0,W_{1/2},W_1$ (i.e. the set on which $|\delta_{0,1}|=1$). In particular $\pi\rho^2\leq \frac{1}{2}\|W\|_\infty^2$, i.e. $\rho\leq \frac{\|W\|_\infty}{\sqrt{2\pi}}$.
	
	Then,
	\[
	\mathcal{C}_{0,1}^-=\int_0^{1/2} (\delta_{0,1}\conv \theta)_{W_t} \d B_t 
	+\int_0^{1/2} \frac{ \langle (\delta_{0,1}\conv \theta)_{W_t},  W_{1/2}-W_t\rangle }{1/2-t}\d t
	+\int_{W_{1/2}}^{0} (\delta_{0,1}\conv \theta)_u \d u, 
	\]
	so it suffices to show these three summands have sone finite exponential moments. 
	
	$\diamond$ For the first one, we have 
	\begin{align*} 
	\mathbb{E}'[\exp (\beta   \int_0^{1/2} (\delta_{0,1}\conv \theta)_{W_t} \d B_t  ) | ]	
	&\leq 
	\mathbb{E}'[\exp (\beta^2   \int_0^{1/2} |(\delta_{0,1}\conv \theta)_{W_t}|^2 \d t  ) | ]^{\frac{1}{2}}\\	
	&\leq \mathbb{E}'[\exp (\frac{\beta^2}{2}    \|\delta_{0,1}\conv \theta\|_\infty^2   ) | ]^{\frac{1}{2}}.
	\end{align*}
	Using Hardy--Littlewood rearrangement inequality, we get 
	\begin{equation}
	\label{eq:boundeddeltatheta}
	|(\delta_{0,1}\conv \theta )_x|\leq \int_{B_{\rho} } |\theta_z| \d z= \rho\leq \|W\|_\infty. 
	\end{equation}
	Thus, 
	\[ 
	\mathbb{E}[\exp (\beta   \int_0^{1/2} (\delta_{0,1}\conv \theta)_{W_t} \d B_t  ) | ]	
	\leq \mathbb{E}[ \exp ( \frac{\beta^2}{2}  \|W\|_\infty^2  )],
	\]
	which is finite for $|\beta|$ small enough. 
	
	$\diamond$
	For the second summand, using \eqref{eq:boundeddeltatheta}, then time reversal, then  $ab\leq (a^2+b^2)/2$ and  Cauchy--Schwarz inequality,
	we get
	\begin{align*}
	\mathbb{E}[\exp (\beta
	\int_0^{1/2}   \frac{\langle (\delta_{0,1}\conv \theta)_{W_t},  W_{1/2}-W_t\rangle}{1/2-t}\d t)]
	&\leq 
	\mathbb{E}[\exp (\beta \|W\|_\infty
	\int_0^{1/2} \frac{|W_{1/2}-W_t|}{1/2-t}\d t)] \\
	&=\mathbb{E}[\exp (\beta \|W\|_\infty
	\int_0^{1/2} \frac{|W_s|}{s}\d s)]\\
	&\leq \mathbb{E}[\exp (\beta \|W\|_\infty^2) ]^\frac{1}{2}
	\mathbb{E}[\exp (\beta 
	(\int_0^{1/2} \frac{|W_s|}{s}\d s)^2)]^\frac{1}{2}.
	\end{align*}
	The first term is finite if $|\beta|$ is small enough. Furthermore, for any integer $k$, 
	\begin{align*}
	\mathbb{E}[ 
	(\int_0^{1/2} \frac{|W_s|}{s}\d s)^{k}]
	&=\int_{[0,1/2]^k} \frac{\mathbb{E}[|W_{s_1}|\dots |W_{s_k} |]}{s_1\dots s_k} \d s_1\dots \d s_k\\
	&\leq \int_{[0,1/2]^k} \frac{\mathbb{E}[|W_{s_1}|^k]^{\frac{1}{k}}   \dots \mathbb{E}[|W_{s_k}|^k]^{\frac{1}{k}}}{s_1\dots s_k} \d s_1\dots \d s_k\\
	&=(\int_{[0,1]} \frac{\mathbb{E}[|W_{s}|^k]^{\frac{1}{k}} }{s}\d s )^k \\
	&=\mathbb{E}[|W_{1}|^k] (\int_{[0,1]} \frac{1}{\sqrt{s}}\d s )^k=2^k\mathbb{E}[|W_{1}|^k],
	\end{align*}
	thus 
	\[
	\mathbb{E}[\exp (\beta 
	(\int_0^{1/2} \frac{|W_s|}{s}\d s)^2)]\leq \mathbb{E}[ \exp (2 |\beta|  |W_{1}|^2 )],
	\]
	which is finite for $|\beta|$ small enough.
	
	$\diamond$
	Finally, for the third summand, 
	\[ 
	\mathbb{E}[  \exp (\beta  \int_{W_{1/2}}^{0} (\delta_{0,1}\conv \theta)_u \d u) ]
	\leq 	\mathbb{E}[  \exp (|\beta|  |W_{1/2}| \|\delta_{0,1}\conv \theta\|_\infty ) ]
	\leq 	\mathbb{E}[  \exp (|\beta|  \|W\|_\infty^2 ) ],
	\]
	which is finite for $|\beta|$ small enough, concluding the proof.
\end{proof}

\section{Exponential moments for \texorpdfstring{$\tilde{\mathcal{A}}_W$}{Ã(W)}   }
\label{sec:expoX}    

In this section we prove that $\tilde{\mathcal{A}}_W$ admits some (two-sided) exponential moments: there exists $\beta_0>0$ such that for $\beta\in(-\beta,\beta)$, $
\mathbb{E}[ \exp(\beta \tilde{\mathcal{A}}_W)  ]<\infty$. We do not track an explicit value for $\beta_c$, but the method gives one. 
Thanks to Theorem \ref{th:cvX}, it suffices to show that the three random variables  $\mathcal{B}=\sum_{k=0}^\infty \sum_{j=1}^{2^k} \mathcal{B}_{k,j}$, $\tilde{\mathcal{C}}=\sum_{k=0}^\infty \sum_{j=1}^{2^k} \tilde{\mathcal{C}}_{k,j}$ and $\tilde{\mathcal{T}}=\sum_{k=0}^\infty \sum_{j=1}^{2^k} \tilde{\mathcal{T}}_{k,j}$ admit some finite exponential moments. Remark that we already know that $\mathcal{B}_{0,1}$ admits such moments (Theorem \ref{th:mainY}). 
Since  $\mathcal{T}_{0,1}^2\leq \|W\|_\infty^2$, the random variable $\tilde{\mathcal{T}}_{0,1}$ also admits some finite exponential moments. The fact that $\tilde{\mathcal{C}}_{0,1}$ admits some finite exponential moments is given by Lemma~\ref{le:Cexpmom}. 

Since we only consider centred random variables which admit the same exact scaling, the existence of finite exponential moments for $\mathcal{B}$,  $\mathcal{C}$ and $T$ follow from the same generic following result. 
\begin{proposition}
	\label{prop:transferExpoMom} 
	Let $\beta_0>0$ and let $(X_{k,j})_{k\in \mathbb{N}, j\in \{1,\dots, 2^k\}}$ be a collection of random variables which satisfies the properties:
	\begin{itemize}
		\item $\mathbb{E}[\exp \beta_0 |X_{0,1} |]<\infty$ and $\mathbb{E}[X_{0,1}]=0$.
		\item For all $k,j$, $X_{k,j}$ is equal in distribution to $2^{-k} X_{0,1}$. 
		\item For any fixed $k$, the random variables $(X_{k,j})_{ j\in \{1,\dots, 2^k\}}$ are globally independent.	
	\end{itemize}
	Then, the sequence of random variables $n\mapsto \sum_{k=0}^n \sum_{j=1}^{2^k} X_{k,j}$ converges almost surely as $n\to \infty$, and the limit $X= \sum_{k=0}^{+\infty} \sum_{j=1}^{2^k} X_{k,j}$ satisfies, for all $\beta<\beta_0/2$,  
	\[ 
	\mathbb{E}[\exp(\beta X )]<\infty. 
	\]
\end{proposition} 
%
To prove this we will use the following H\"older inequality for infinite product.\footnote{I believe there is a small mistake in \cite{Karakostas}: in the case the infinite product $\prod \|f_k\|_{p_k}$ is infinite, it does not hold in general that the infinite product $f$ lies in $L^1$. 
	I am using here a corrected version.}
\begin{theorem}[\cite{Karakostas}]
	\label{th:Karakostas}
	Let $(p_k)$ be a sequence of reals larger than $1$, such that ${\sum_{k\in \mathbb{N}} 1/p_k=1}$. Let $(f_k)$ be a sequence of functions with $f_k\in L^{p_k}(X,\mu)$, for some $\sigma$-finite measure space. If $\prod \|f_k\|_{p_k}$ converges to a limit in $(0,+\infty]$ and $\prod f_k$ converges almost everywhere toward $f$, then 
	$\|f\|_{L^1}\leq \prod_i \|f_k\|_{L^{p_k}}$. In particular, if the right-hand side is finite, then $f\in L^1(X,\mu)$.
\end{theorem}
We will also use the following: 
\begin{lemma} 
	\label{le:beta2}
	Let $X$ be a random variable with zero mean and let $\beta_0>0$ and $\beta_1 \in(0,\beta_0)$. Assume  $\mathbb{E}[\exp(\beta_0 |X|)]<\infty$. Then, there exists a constant $C$ such that for all $\beta'\in[-\beta_1,\beta_1]$, 
	\[ 
	\mathbb{E}[\exp(\beta' X)]\leq \exp( C {\beta'}^2).
	\]
\end{lemma} 
\begin{proof} 
	Expand the exponential into a series, which is absolutely convergent, locally uniformly below the radius of convergence hence uniformly on $[-\beta_1,\beta_1]$. We get, for some finite positive constant $C$ which depends on $\beta_1$ and on the distribution of $X$ but not on $\beta'$, 
	\[ 
	\mathbb{E}[\exp(\beta' X)]= 1+ 0+ \sum_{k=2}^\infty \frac{\beta'^k}{k!} \mathbb{E}[X^k]\leq 1+ C \beta'^2 \leq \exp(C \beta'^2).
	\]
\end{proof} 

\begin{proof}[Proof of Proposition \ref{prop:transferExpoMom} ]
	Let $\beta_0,\beta$ as in the statement of the proposition.
	Let $q<2 $ be such that $\frac{\beta q }{q-1}<\beta_0$, which exists because $2\beta<\beta_0$ and $\frac{q}{q-1}\underset{q\to 2}\longrightarrow 2$. For $k$ a non-negative integer, let $p_k=q^{k+1}/(q-1)\in(1,\infty)$. Remark $\sum_{k=0}^\infty 1/p_k=1$. Let $f_k =\sum_{j=1}^{2^k} X_{k,j}$. 
	
	Using the independence between the $X_{k,j}$ and then the scaling relation $X_{k,j}\overset{(d)}=2^{-k}X_{0,1}$, we get $\mathbb{E}[f_k^2]=\sum_{j=1}^{2^k} \mathbb{E}[X_{k,j}^2]=2^{-k} \mathbb{E}[X_{0,1}^2]$.
	Since this is summable over $k$ and $\mathbb{E}[f_k]=0$, we deduce that $n\mapsto \sum_{k=0}^n f_k$ is a martingale which is bounded in $L^2$, hence convergent, almost surely and in $L^2$. In particular, for all $\beta \in \mathbb{R}$, $\exp(\beta \sum_{k=0}^{n}  f_k )$ converges almost surely toward $\exp(\beta X)$. Furthermore, by Jensen inequality, for all $\beta$, 
	\[
	\mathbb{E}[\exp( \beta p_k f_k ) ]\geq \exp(  \beta p_k \mathbb{E}[f_k] )=1.
	\]
	Hence, the sequence $ n \mapsto \prod_{k=1}^n \| \exp( \beta  f_k)\|_{p_k}$ is non-decreasing, thus converges in $(1,+\infty]$. We can thus apply the H\"older inequality for infinite product (Theorem \ref{th:Karakostas} ), and we deduce 
	\begin{equation}
		\label{eq:tempX1}
	\mathbb{E}[ \exp( \beta X)]\leq  \prod_{k=0}^\infty \mathbb{E}[ \exp(\beta \frac{q^{k+1}}{q-1} f_k   ) ]^{ \frac{q-1}{q^{k+1}} }=\Big(  \prod_{k=0}^\infty \mathbb{E}[ \exp(\frac{\beta q}{q-1}  q^{k } f_k   ) ]^{q^{-k} } \Big)^{ \frac{q-1}{q} }. 
	\end{equation}
	Using again the independence and scaling properties, it holds, for any $k$ and $\beta'$, 
	\begin{equation}
		\label{eq:tempX2} 
	\mathbb{E}[\exp( \beta' f_k  )]= 
	\mathbb{E}[\exp( \beta'  2^{-k}  X_{0,1}  )]^{2^k}.
	\end{equation}
	Combining \eqref{eq:tempX1} and \eqref{eq:tempX2}, we obtain 
	\begin{equation}
	\label{eq:tempX3}
	\mathbb{E}[ \exp( \beta X)]
	\leq  \Big( \prod_{k=0}^\infty \mathbb{E}[ \exp(\frac{\beta q }{q-1} q^k 2^{-k}   X_{0,1}  ) ]^{ q^{-k} 2^k  }  \Big)^{ \frac{q-1}{q } } 
	\end{equation}
For all $k\geq 0$, 	$\frac{\beta q }{q-1} q^k 2^{-k}\leq {\beta q}{q-1}<\beta_0$ (recall $q<2$). Applying  Lemma \ref{le:beta2} with $\beta_1= {\beta q}{q-1}$, we deduce that there exists a constant $C$ which depends on $\beta$, $q$, and on the distribution of $X_{0,1}$, but not on $k$, such that for all $k\geq 0$, 
\begin{equation}
	\label{eq:tempX4}
\mathbb{E}[ \exp(\frac{\beta q }{q-1} q^k 2^{-k}   X_{0,1}  ) ]
\leq \exp(C   q^{2k} 2^{-2k}      )).
\end{equation}
Combining  \eqref{eq:tempX3} and \eqref{eq:tempX4}, we obtain 
\[
	\mathbb{E}[ \exp( \beta X)]
\leq  \Big( \prod_{k=0}^\infty  \exp(C   q^{2k} 2^{-2k}      )^{ q^{-k} 2^k  }  \Big)^{ \frac{q-1}{q } } 
=   \exp(  \tfrac{C(q-1)}{q }   \sum_{k=0}^\infty  q^{k} 2^{-k} ) ,
\]
	which is finite since $q<2$.
\end{proof}

\section{Estimation of the average \texorpdfstring{$\mathbb{E}[\mathcal{A}^\epsilon_W]$}{E[Ã(W) ]}   } 
\label{sec:average}
\subsection{Existence of an asymptotic expansion}

\begin{lemma} 
	\label{le:expand}
	There exist constants $C,C'$ such that, for $\epsilon=2^{-n/2}, n\in \mathbb{N}$,
	as $\epsilon\to 0$, \[\mathbb{E}[\mathcal{A}^{\epsilon}_W]=  
	C \log(\epsilon^{-1})+ C' +O(\epsilon^{\frac{1}{2}}).\]
	The constant $C$ does not depend on $\varphi$.
\end{lemma}
\begin{proof} 
	By scaling, for $\epsilon=2^{-n/2}$, \[
	\mathbb{E}[\mathcal{A}^{\epsilon}_{n,j}  ]=\epsilon^2\mathbb{E}[\mathcal{A}^{1}_{0,1}],\quad  
	\mathbb{E}[\mathcal{T}_{k,j}  ]=2^{-k}\mathbb{E}[\mathcal{T}_{0,1}],\quad  
	\mathbb{E}[\mathcal{C}_{k,j} ]=2^{-k}\mathbb{E}[\mathcal{C}_{0,1}],\]\[  
	\mathbb{E}[\mathcal{C}^{\epsilon}_{k,j}-\mathcal{C}_{k,j} ]=2^{-k}\mathbb{E}[\mathcal{C}^{2^{(k-n)/2}}_{0,1}\hspace{-0.1cm}-\mathcal{C}_{0,1}], \quad 
	\mathbb{E}[\mathcal{T}^{\epsilon}_{k,j}-\mathcal{T}_{k,j} ]=2^{-k}\mathbb{E}[\mathcal{T}^{2^{(k-n)/2}}_{0,1}\hspace{-0.1cm}-\mathcal{T}_{0,1}].
	\]
	Using the decomposition \eqref{eq:decompoX} with $m=n$, the fact that $\mathbb{E}[\mathcal{B}^\epsilon_{k,j}]=0$, and the scaling relations above, we obtain 
	\begin{align*}
		\mathbb{E}[\mathcal{A}^\epsilon_W] 
&=
 \sum_{j=1}^{2^n}  \mathbb{E}[\mathcal{A}^\epsilon_{m,j}]
+ \sum_{k=0}^{n-1}  \sum_{j=1}^{2^k} \big( -\mathbb{E}[\mathcal  T_{k,j}]
- \mathbb{E}[\mathcal  T_{k,j}^\epsilon-\mathcal T_{k,j}]
+ 2\mathbb{E}[\mathcal{C}_{k,j}]
+ 2  \mathbb{E}[\mathcal{C}^\epsilon_{k,j}-\mathcal{C}_{k,j}]\big) \\
&= \mathbb{E}[\mathcal{A}^{1}_{0,1}]
+ \sum_{k=0}^{n-1} \big( - \mathbb{E}[\mathcal  T_{0,1}]
-  \mathbb{E}[\mathcal  T_{0,1}^{2^{(k-n)/2}}-\mathcal T_{0,1}]
+ 2   \mathbb{E}[\mathcal{C}_{0,1}]
+ 2   \mathbb{E}[\mathcal{C}^{2^{(k-n)/2}}_{0,1}-\mathcal{C}_{0,1}]\big)\\
&=  \mathbb{E}[\mathcal{A}^{1}_{0,1}]+n (  2   \mathbb{E}[\mathcal{C}_{0,1}]-  \mathbb{E}[\mathcal  T_{0,1}] )
+\sum_{l=1}^{n}(2   \mathbb{E}[\mathcal{C}^{2^{-l/2}}_{0,1}-\mathcal{C}_{0,1}] 
- \mathbb{E}[\mathcal  T_{0,1}^{2^{-l/2}}-\mathcal T_{0,1}]
)
\end{align*} 
By Lemma \ref{le:Tcv} and Lemma \ref{le:Zcv}, there is a universal constant $C''$ such that for all $l\geq 0$, 
\[
|\mathbb{E}[\mathcal{T}^{2^{-l/2}}_{0,1}-\mathcal{T}_{0,1}]| \leq C'' 2^{-l/2} 
\quad \text{and} \quad 
|\mathbb{E}[\mathcal{C}^{2^{-l/2}}_{0,1}-\mathcal{C}_{0,1}]| \leq C'' 2^{-l/4}.  
\]
In particular, 
\[C'\coloneqq 
\mathbb{E}[\mathcal{A}^{1}_{0,1}]
+ 
\sum_{l=1}^\infty (2   \mathbb{E}[\mathcal{C}^{2^{-l/2}}_{0,1}-\mathcal{C}_{0,1}] 
- \mathbb{E}[\mathcal  T_{0,1}^{2^{-l/2}}-\mathcal T_{0,1}]
)
\] 
is well-defined, as the sum is absolutely convergent. 
Furthermore, 
	\begin{align*}
\big|	\mathbb{E}[\mathcal{A}^\epsilon_W] - n (  2   \mathbb{E}[\mathcal{C}_{0,1}]-  \mathbb{E}[\mathcal  T_{0,1}] ) - C' \big|
	&=
\big| \sum_{l=n+1}^{\infty}(2   \mathbb{E}[\mathcal{C}^{2^{-l/2}}_{0,1}-\mathcal{C}_{0,1}] 
- \mathbb{E}[\mathcal  T_{0,1}^{2^{-l/2}}-\mathcal T_{0,1}] ) \big|\\
&\leq C'' \sum_{l=n+1}^{\infty}  2^{-l/4+1}+2^{-l/2}\leq C''' 2^{-n/4}=C'''\epsilon^\frac12,
\end{align*} 
which concludes the proof. 
\end{proof} 
Remark that when we consider a Brownian motion with duration $t\neq 1$, we can easily see how these constants $C,C'$ depends on $t$ by a scaling argument: 
\[ 
\mathbb{E}_t[\mathcal{A}^{\epsilon}_W]=  
t \mathbb{E}_1[\mathcal{A}^{\epsilon/\sqrt{t}}_W]=
tC \log(\epsilon^{-1}\sqrt{t})+ C' +O(\epsilon^\frac12)
= tC \log(\epsilon^{-1})+ (t C'+ \frac{t \log(t) C}{2} ) +O(\epsilon^\frac12).
\]

We will now prove $C=\frac{1}{2\pi}$. This is expected from \cite[Theorem 1]{Werner3}, but the regularisation method used in \cite{Werner3} is not the same as ours, and it doesn't seem that there is a way to prove the coefficients coming from both regularisations should be equal, although it is intuitively clear. The intuitive argument goes as follows. The points at distance $\geq \max(1,C_\phi) \epsilon$ to $\operatorname{Range}(W)$ do not `feel' any of the two regularisations, and there contribution to $\mathcal{A}^\epsilon_W$ will be the same. As for the points at distance $\leq \max(1,C_\phi) \epsilon$ to $\operatorname{Range}(W)$, there total contribution, for both regularisation, is of order $1$. The difference in the regularisation methods should thus only affect the $O(1)$ term in the asymptotic expansion of $\mathbb{E}[\mathcal{A}^\epsilon_W]$ (or in fact, of $\mathcal{A}^\epsilon_W$). 

\subsection{Computation of the Constant $C$}
\label{sec:average2} 
In order to compute $C$ ($=1/2\pi$, \emph{a posteriori}), we consider, for $W$ a Brownian motion with duration $1$ and started from $0$,  
\[ 
\n_W^{\epsilon,\circ}(z)= \int_0^1 \theta^\epsilon_{W_t-z} \circ \d W_t. 
\]
Compared with $\n_W^\epsilon$, this does not include the `Boundary term' $\int_{W_1}^{W_0}
\theta^\epsilon_{v-z}  \d v$. When $\epsilon\to 0$, the difference between $\n_W^\epsilon$ and 
$\n_W^{\epsilon,\circ}$ converges to the angle $\widehat{W_0 z W_1}$, divided by $2\pi$. This angle has very poor integrability properties with respect to $z$, in particular 
$\int_{\mathbb{R}^2} (\widehat{W_0 z W_1})^2 \d z=+\infty$,
so we must introduce a cut-off to limit infrared divergences. We therefore fix a large radius $C_\epsilon$, which is such that there is only a small probability for $W$ to exit $B_{C_\epsilon}$: in the large probability event that this does not happen and for $z$ outside this ball, we will use $\n_W^{\epsilon,\circ}(z)=0$ rather than comparing $\n_W^{\epsilon}(z)$ with $\n_W^{\epsilon,\circ}(z)\simeq |z|^{-1}$ . Inside this ball, we compare $\n^\epsilon_{W}(z)$ with $\n^{\epsilon,\circ}_{W}(z)$ and estimate the later. To keep track of the order of the various error terms, let us say already that we will ultimately choose $C_\epsilon=8 \sqrt{\log(\log \epsilon^{-1})}$, although we will not use this until the very end of the computation.
\begin{lemma}
	\label{le:momBulk}
	Let $C_\epsilon$ be a function of $\epsilon$ such that $C_\epsilon\geq 4\sqrt{2}$. As $\epsilon \to 0$,
	\[\mathbb{E}\Big[
	\int_{B_{C_\epsilon}}  \n^\epsilon_{W}(z)^2 \d z \Big]  =
	\frac{1}{2\pi}\log(\epsilon^{-1})  +O\big(C_\epsilon^2 + e^{-\frac{C_\epsilon^2}{8} }\log(\epsilon^{-1} ) +C_\epsilon\sqrt{\log(\epsilon^{-1})} \big),
	\]
\end{lemma}
\begin{proof}
	Let $\bar{I}$ be the left-hand side, and let
	\[
	I=\mathbb{E}\Big[
	\int_{B_{C_\epsilon}}  \n^{\epsilon,\circ}_{W}(z)^2 \d z \Big]. \]
Again, the Stratonovich integral defining $\n^{\epsilon,\circ}_{W}(z)$ is equal to the identical It\^o integral, because $\div \theta=0$. By It\^o isometry, 
	\begin{equation}
	\label{eq:temp:expectation}
			 \mathbb{E}[\n^{\epsilon,\circ}_W(z)^2]=\int_{\mathbb{R}^2} \int_0^1  p_t(0,y) |\varphi^\epsilon\conv \theta(z-y)|^2 \d y \d t.
	\end{equation}

	Splitting the integral over $y\in \mathbb{R}^2$ into three parts, we get $I=J+R_1+R_2$ with
	\begin{align*}
	J&\coloneqq \epsilon^{-2}\int_{B_{C_\epsilon}} \int_{B_\frac{C_\epsilon}{2}} \int_0^1 p_t(0,y)  |(\varphi^\epsilon \conv \theta)( z- y)|^2 \d t\d y\d z,\\
	R_1&\coloneqq  \int_{B_{C_\epsilon}} \int_{B_{2C_\epsilon}\setminus B_\frac{C_\epsilon}{2} } \int_0^1 p_t(0,y)  |(\varphi^\epsilon \conv \theta)( z- y)|^2 \d t\d y\d z,\\
	\shortintertext{and} R_2&\coloneqq   \int_{B_{C_\epsilon} } \int_{\mathbb{R}^2\setminus B_{2C_\epsilon}}\int_0^1 p_t(0,y)  |(\varphi^\epsilon \conv \theta)( z- y)|^2 \d t\d y\d z.
	\end{align*}
	To estimate $J$, we swap the two integrals in its definition, then we use the change of variable $z'=y-z$, and then the inclusions
	\begin{equation}
	\label{eq:inc}
	\forall y\in B_\frac{C_\epsilon}{2} , \qquad B_\frac{ C_\epsilon }{2}
	\subset B_{C_\epsilon}(y) \subset B_{\frac{3}{2} C_\epsilon}.
	\end{equation}
We get 
\begin{equation} \label{temp:Jbetw} \int_{B_\frac{C_\epsilon}{2}}  \int_0^1 p_t(0,y) \d t \d y  \int_{B_\frac{C_\epsilon}{2}}|(\varphi^\epsilon \conv \theta)( z')|^2 \d z'
	\leq J\leq      \int_{B_\frac{3C_\epsilon}{2}}    |(\varphi^\epsilon \conv \theta)(z')|^2 \d z',
\end{equation}
where for the second inequality we used $\int_{\mathbb{R}^2} p_t(0,z)\d z=1$. Let $k\in \{1,3\}$. 
Using the definition of $\varphi^\epsilon$ and the scaling property $\theta (\lambda z)=\lambda^{-1}\theta(z)$, we deduce $\varphi^\epsilon\conv \theta(z)=\epsilon^{-1} (\varphi \conv \theta) (\epsilon^{-1} z)$, hence 
\[ 
\int_{B_\frac{kC_\epsilon}{2}}   |(\varphi \conv \theta)(\epsilon^{-1}z')|^2 (\epsilon^{-2}\d z')=
\int_{B_\frac{k\epsilon^{-1}C_\epsilon}{2}}   |(\varphi \conv \theta)(w)|^2 \d w.
\]
Using the fact that $\varphi\conv \theta$ is locally integrable in the vicinity of $0$,
that $|(\varphi \conv \theta)(w)|= 1/(2\pi|w|)+O(1/|w|^2)$ as $|w|\to \infty$ (recall $\varphi$ is compactly supported), and that $C_\epsilon$ is bounded below, we get 
\begin{equation}
\label{eq:temp:thel} 
\int_{B_\frac{k\epsilon^{-1}C_\epsilon}{2}}   |(\varphi \conv \theta)(w)|^2 \d w
\underset{\epsilon\to \infty}{=} \frac{1}{2\pi} \int_1^{ \frac{k\epsilon^{-1}C_\epsilon}{2}} \frac{\d r}{r}+O(1)= \frac{\ln(\epsilon^{-1})+\ln(C_\epsilon ) }{2\pi}+O(1).
\end{equation}
	For $y\geq \sqrt{2t}$, it holds $\partial_t p_t(0,y)\geq 0$, and it follows from the condition $C_\epsilon\geq 4\sqrt{2}$ that for all
	$y\in \mathbb{R}^2\setminus B_{C_\epsilon/2}$ and all $t\in [0,1]$, $p_t(0,y)\leq p_1(0,y)$. Thus,
	\begin{equation}
		\label{eq:temp:expb} 
	\int_{\mathbb{R}^2 \setminus B_\frac{C_\epsilon}{2}} \int_0^1 p_t(0,y) \d t\d y\leq \int_{\mathbb{R}^2 \setminus B_\frac{C_\epsilon}{2}}  p_1(0,y) \d y = e^{-\frac{C_\epsilon^2}{8}}.
	\end{equation}
	Inserting \eqref{eq:temp:thel} and \eqref{eq:temp:expb} inside \eqref{temp:Jbetw}, we obtain
	\[  (1- e^{-\frac{C_\epsilon^2}{8}}) ( \frac{\log( \epsilon^{-1})+\log(C_\epsilon) }{2\pi}   +O(1)) \leq J\leq  \frac{\log( \epsilon^{-1})+\log(C_\epsilon) }{2\pi}   +O(1).\]
	We now focus on the estimation of the residual terms $R_1$ and $R_2$. Let $C<\infty$ be such that for all $y \in \mathbb{R}^2$, $|(\varphi \conv \theta)(y)|^2\leq C|y|^{-2}$. Then (recall $\varphi^\epsilon\conv \theta(z)=\epsilon^{-1} (\varphi \conv \theta) (\epsilon^{-1} z)$), 
	$|(\varphi^\epsilon \conv \theta)(y)|^2\leq C|y|^{-2}$, and therefore 
	\[R_1\leq \int_{B_{2C_\epsilon} \setminus B_{C_\epsilon/2}} \frac{C}{  |y|^2} \d y
	\leq \pi (2 C_\epsilon)^2 \frac{C}{  (C_\epsilon/2)^2  }  =
	16 \pi C.\]
	As for $R_2$, we notice that for all $y\in \mathbb{R}^2\setminus B_{2C_\epsilon}$ and $z\in B_{C_\epsilon}$, it holds $|y-z|\geq \frac{|y|}{2}$, thus $p_t(y,z)\leq p_t(0, \frac{y}{2})$. Hence,
	\begin{align*}
	R_2&\leq \int_{B_{C_\epsilon}} \int_{\mathbb{R}^2\setminus B_{2C_\epsilon}} \int_0^1 p_t(0,\frac{y}{2})\d t  \frac{C}{|y|^2 }  \d y\d z=(\pi   C_\epsilon^2) 2\pi  \int_0^1 \int_{2C_\epsilon}^\infty    p_t(0,\frac{r}{2})\frac{C}{r^2} r \d r\d t=o(1).
	\end{align*}
	Altogether, we obtain
	\begin{align}
	I& = J+R_1+R_2= \frac{\log(\epsilon^{-1})+ \log(C_\epsilon )}{2\pi}+O(e^{-\frac{C_\epsilon^2}{8} }\log(\epsilon^{-1} C_\epsilon)  )+O(1)\nonumber\\
	&=\frac{1}{2\pi}\log(\epsilon^{-1} )+O(\log(C_\epsilon)+ e^{-\frac{C_\epsilon^2}{8} }\log(\epsilon^{-1} )  ).
	\label{eq:temp:prebulk}
	\end{align}
	
	In order to conclude, we compare $I$ with $\bar{I}$, for which we apply the generic identity
	\begin{equation}
	\label{eq:trigo}
	\big|\| a\|^2- \| b\|^2 \big|=  \big|\| a-b\|^2+2\langle a-b, b\rangle  \big|    \leq \| a-b\|^2
	+ 2 \| a-b\| \| b\|
	\end{equation}
	with $a=\n^\epsilon_{W}$, $b=\n^{\epsilon,\circ}_{W}$, and with the Hilbert norm $\|f\|^2=\int \mathbb{E}_0[f(z)^2]\d z$ on $L^2(\Omega\times B_{C_\epsilon} )$. Since  $|\n^\epsilon_{W}(z)-\n^{\epsilon,\circ}_W(z) |\leq 1$
	for all $z$, we have $\|a-b\|\leq \sqrt{\pi} C_\epsilon$. By \eqref{eq:temp:prebulk}, $\|b\|^2=O( \log(\epsilon^{-1}C_\epsilon) )$, and using \eqref{eq:temp:prebulk} a second time we obtain
	\begin{align}
	\int_{B_{C_\epsilon}} \mathbb{E}_0[\n^\epsilon_{W}(z)^2]\d z &= \|b\|^2+O(\|a-b\|^2+\|a-b\|\|b\|)\nonumber\\
	&  =\frac{1}{2\pi}\log(\epsilon^{-1})  +O(C_\epsilon^2 + e^{-\frac{C_\epsilon^2}{8} }\log(\epsilon^{-1} ) +C_\epsilon\sqrt{\log(\epsilon^{-1})} ),
	\label{eq:temp:bulk}
	\end{align}
	where we simplified some smaller order terms such as $O(C_\epsilon\sqrt{\log(C_\epsilon)})$ for the last inequality.
\end{proof}

\begin{lemma}
	\label{le:momTail}
	Let $C_\epsilon$ be a function of $\epsilon$ such that $C_\epsilon\geq 4\sqrt{2}$. Then, as $\epsilon  \to 0$,
	\[  \mathbb{E}_0[ \int_{\mathbb{R}^2\setminus B_{C_\epsilon}} \n^\epsilon_{W}(z)^2]\d z=O\Big( \frac{e^{-\frac{17 }{512}C_\epsilon^2 }}{C_\epsilon}\log(\epsilon^{-1})\Big)+o(1). \]
\end{lemma}
\begin{proof}
	For $|z|\geq C_\epsilon$, let $\tau_z$ be the first time when $|W|=\frac{|z|}{2}$, and notice that, for topological reasons, \[\tau_z\geq 1\implies \forall y\in \mathbb{R}^2\setminus B_{|z|/2},\ \n_{W}(y)=0.\]
	For $\epsilon$ sufficiently small (which we now assume), $C_\epsilon\geq 2 \epsilon \sup\{ |y|, y\in \Supp(\varphi) \}$. Then, for all $z\in \mathbb{R}^2\setminus B_{C_\epsilon} $,
	\[\tau_z\geq 1\implies  \n^\epsilon_{W}(z)=0.\]
	This follows from the fact that, under the condition $\tau_z\geq 1$, the vector field  $\phi^\epsilon\conv \theta$ is harmonic on the simply connected set $B_{\|W\|}$, and the Stratonovich integral of a harmonic vector field along a closed loop is equal to $0$.
		
	On the event $\tau_z<1$, let $\mathcal{W}:t\in [0,1-\tau_z] \mapsto W_{\tau_z+t}$. Conditionally on
	$(\tau_z, (W_t)_{t\leq \tau_z})$, $\mathcal{W}$ is a Brownian motion of duration $\tau_z$ started from $W_{\tau_z}$. 
	The integer $\n^\epsilon_{W}(z)-\n^\epsilon_{\mathcal{W}}(z)$ is equal (up to sign) to the winding around $z$ of the triangle with vertices $0$, $W_{\tau_z}$, and $W_1$, which is at most $1$. Thus, $\n^\epsilon_{W}(0)$ and $\n^\epsilon_{\mathcal{W}}(0)$ differs from each other from at most $2$.
	Using $(a+b)^2\leq 2(a^2+b^2)$, we deduce
	\begin{align}
		\mathbb{E}[ \n^\epsilon_{W}(z)^2  ]
		&\leq \mathbb{P}(\tau_z\leq 1)
		\sup_{\substack{ y \in \partial B_{|z|/2}\\ t\in [0,1]   } }\mathbb{E}_{t,y}[ (\n^\epsilon_{\mathcal{W}}(z)+R^\epsilon)^2  ],&& |R^\epsilon|\leq 2 \nonumber\\
		&\leq 2\mathbb{P}(\tau_z\leq 1)
		( \sup_{\substack{ y \in \partial B_{|z|/2}\\ t\in [0,1]   } } \mathbb{E}_{t,y}[ (\n^\epsilon_{\mathcal{W}}(z))^2  ]+4)\nonumber \\
		&= 2\mathbb{P}_0(\tau_z\leq 1)
		( \sup_{\substack{ y \in \partial B_{|z|/2}\\ t\in [0,1]   } } \mathbb{E}_{t,0}[ (\n^\epsilon_{\mathcal{W}}(z-y))^2  ]+4)\label{eq:temp:abc1}
	\end{align}
	For $z\in\mathbb{R}^2\setminus B_{4\sqrt{2}}$, $y\in\partial B_{|z|/2}  $, and $w\in B_1$, we have $|z-y|-|w|\geq |z|/2-1\geq |z|/4\geq \sqrt{2}$, hence for all and $s\in[0,1]$, $p_s(z-y,w)\leq p_s(0,|z|/4)\leq p_1(0,|z|/4)$.
	From Equation \eqref{eq:temp:expectation} and using $|\varphi^\epsilon\conv \theta(y)|^2\leq C/(\epsilon^2 +|y|^2)$, for a constant $C<\infty$ which does not depend on $\epsilon\in(0,1] $ nor on $y\in \mathbb{R}^2$, we have, uniformly over $z: |z|\geq 4\sqrt{2}$ and $\epsilon\in(0,1]$, 
	\begin{align}
		\mathbb{E}_{t,0}[ \n^\epsilon_{W}(z-y)^2  ]
		&\leq \int_0^t \int_{\mathbb{R}^2} p_s(z-y,w)\frac{C}{\epsilon^2+|w|^2} \d w \d s \nonumber\\
		&\leq  p_1(0,\tfrac{z}{4}) \d s \int_{B_1}  \frac{C}{\epsilon^2+|w|^2} \d w +\int_0^t \int_{\mathbb{R}^2\setminus B_1} p_s(z-y,w) C \d w \d s\nonumber\\
		&\leq C p_1(0,\tfrac{z}{4}) (\log(\epsilon^{-1})+\frac{\log(2)}{2}) +C t\nonumber\\
		&= O( e^{-\frac{|z|^2 }{32}} \log(\epsilon^{-1}))+O(1). \label{eq:temp:abc2}
	\end{align}
	
	Besides, since $\{(x^1,x^2): x^1\leq c, x^2\leq c\}\subset B_{\sqrt{2}c}$ for all $c$, using the reflection principle for the $1$-dimensional Brownian motion, we have
	\begin{align}
		\mathbb{P}_{0}(\tau_z\leq 1)&\leq 2 \mathbb{P}_0( \exists s\in [0,1]: |W^1_s|\geq \frac{|z|}{4})\nonumber\\
		&\leq 4 \mathbb{P}_0( \exists s\in [0,1]: W^1_s\geq \frac{|z|}{4})\nonumber\\
		&=4 \mathbb{P}_0( |W^1_1|\geq \frac{|z|}{4})= 8 \mathbb{P}_0( W^1_1\geq \frac{|z|}{4})= \frac{8}{\sqrt{2\pi} } \int_{\frac{|z|}{4}}^\infty e^{- \frac{r^2}{32}} \d r \underset{|z|\to \infty}{=}O \Big( \frac{e^{-\frac{|z|^2}{512} }}{|z|}\Big). \label{eq:probaTau}
	\end{align}
	Together with \eqref{eq:temp:abc1} and \eqref{eq:temp:abc2}, we get
	\[ \mathbb{E}[ \n^\epsilon_{W}(z)^2]=O \Big(\frac{e^{-\frac{17 |z|^2}{512} }}{|z|}\log(\epsilon^{-1})+\frac{e^{-\frac{|z|^2}{512} }}{|z|}\Big).\]
	By integrating over $z$, we deduce that
	\[
	\int_{\mathbb{R}^2\setminus B_{C_\epsilon}}
	\mathbb{E}[ \n^\epsilon_{W}(z)^2] \d z=O\Big( \frac{e^{-\frac{17 }{512}C_\epsilon^2 }}{C_\epsilon}\log(\epsilon^{-1})\Big)+o(1).
	\]
\end{proof}
\begin{corollary}
	\label{coro:mom}
	There exists a constant $C$, which depends on the mollifier $\varphi$, such that, as $\epsilon\to 0$,
	\[
	\mathbb{E}\Big[ \int_{\mathbb{R}^2} \n^\epsilon_{W}(z)^2\d z \Big]= \frac{\log(\epsilon^{-1})}{2\pi}+ C+ O(\epsilon).
	\]
\end{corollary}
\begin{proof}
Applying Lemma \ref{le:momBulk} and Lemma \ref{le:momTail} with $C_\epsilon=8 \sqrt{\log(\log(\epsilon^{-1}))}$, we get 
	\[
\mathbb{E}_x \Big[ \int_{\mathbb{R}^2} \n^\epsilon_{W}(z)^2\d z \Big]= \frac{\log(\epsilon^{-1})}{2\pi}+ O\big(\log(\log(\epsilon^{-1}))\sqrt{\log(\epsilon^{-1})} \big).
\]
Thus, the constant $C$ in Lemma \ref{le:expand} is equal to $1/(2\pi)$, which concludes the proof.  
%
\end{proof}

\appendix 

\section{Derivation below rough path integrals with parameter} 
\label{sec:RPderivation}
\begin{lemma}
	\label{le:derivBelowRP} 
	Let $V: \mathbb{R}^a\times \mathbb{R}^b \to\mathcal{L}(\mathbb{R}^a, \mathbb{R}^c)$, where $\mathcal{L}(\mathbb{R}^a, \mathbb{R}^c)$ is the space of linear maps from $\mathbb{R}^a$ to $\mathbb{R}^b$, and assume that $V\in \mathcal{C}^3$. Let $\mathbf{X}=(X,\mathbb{X})$ be a $\mathscr{C}^\alpha([0,T], \mathbb{R}^a)$-rough path, with $\alpha>1/3$. For $x\in \mathbb{R}^a$ and $u\in \mathbb{R}^b$, let $\nabla^{\mathbb{R}^b}_h V(x,u)\coloneqq (\nabla_h( V(x,\cdot)))(u) $. 
	
	Then, for all $u,h\in \mathbb{R}^c$, the rough path integrals 
	\[ 
	\int_0^T V(X,u) \d \mathbf{X} \qquad \text{and} \qquad 	\int_0^T  \nabla^{\mathbb{R}^b}_h  V(X,u) \d \mathbf{X}
	\]
	varie continuously $b$. Furthermore, the map $u\mapsto \int_0^T V(X,u) \d \mathbf{X}$ is $\mathcal{C}^1$, and 
	\[ 
	\int_0^T  \nabla^{\mathbb{R}^b}_h  V(X,u) \d \mathbf{X}=  \nabla_h 	\int_0^T   V(X,u) \d \mathbf{X}.
	\]
\end{lemma}
\begin{proof}
	Fix $u\in \mathbb{R}^b$.  Eventually replacing $V$ with a different linear map which agrees with $V$ in a neighbourhood of $\Range(X)\times \{u\}$, which does not change the value of the integrals  involved, we can assume without loss of generality that the third order derivatives of $V$ are globally bounded and uniformly continuous. 
	It follows from \cite[Theorem 4.4, Equation (4.16)]{FrizHairer}\footnote{This bound can certainly be attributed to Terry J. Lyons, as  \cite{FrizHairer} acknowledges.}
	that there exists a finite constant $C_\mathbf{X}$, which depends on the rough path $\mathbf{X}$, such that 
	that for $F:\mathbb{R}^a \to \mathcal{L}(\mathbb{R}^a, \mathbb{R}^c)$, 
	\[ 
	| \int_0^T  \nabla^{\mathbb{R}^b}_h  F(X,u) \d \mathbf{X}|\leq C_\mathbf{X} \|F\|_{\mathcal{C}^2_b}.
	\]
	In particular, for any $\epsilon>0$ and $h\in \mathbb{R}^b$, 
	\begin{align*}
		| \int_0^T  ( \frac{V(X,u+\epsilon h)-V(X,u)}{\epsilon} & - \nabla^{\mathbb{R}^b}_h V(X,u)) \d \mathbf{X}|\\
		& \leq C_\mathbf{X} \| x\mapsto \frac{V(x,u+\epsilon h)-V(x,u)}{\epsilon}  - \nabla^{\mathbb{R}^b}_h V(x,u)  \|_{\mathcal{C}^2_b}\\
		& \leq C_\mathbf{X} \max_{ i,j } \|  \frac{\partial_{ij} V(x,u+\epsilon h)-\partial_{ij} V(x,u)}{\epsilon}  -  \nabla^{\mathbb{R}^b}_h \partial_{ij} V(x,u) \|_{\infty}\\
		&\leq C_\mathbf{X} \sup_{\substack{  i,j \\ \epsilon'\in [0,\epsilon]}} 
		\| \nabla^{\mathbb{R}^b}_h \partial_{ij} V(x,u+\epsilon'h) - \nabla^{\mathbb{R}^b}_h \partial_{ij} V(x,u)\|_\infty \\
		&\leq b \ C_\mathbf{X} \max_{i,j,k} \omega_{\epsilon |h|}(\partial_{ijk} V),
	\end{align*}
	where $\omega_{\epsilon'}(f)$ is the continuity modulus of $f$, i.e. 
	\[ 
	\omega_{\epsilon'}(f)=\sup \{ |f(x)-f(y)|: |x-y|\leq \epsilon'\}.
	\] 
	Since $\partial_{ijk} V$ is uniformly continuous, $\omega_{\epsilon'}(\partial_{ijk} V)\underset{\epsilon' \to 0}\longrightarrow 0$, and thus 
	\begin{align*}
		\frac{	\int_0^T  V(X,u+\epsilon h) \d \mathbf{X}  -\int_0^T V(X,u) \d \mathbf{X} 
		}{\epsilon}&=
		\int_0^T  \frac{V(X,u+\epsilon h)-V(X,u)}{\epsilon} \d \mathbf{X}\\
		&	\underset{\epsilon \to 0}\longrightarrow \int  \nabla^{\mathbb{R}^b}_h V(X,u)  \d \mathbf{X},
	\end{align*}
	which proves that $u\mapsto \int_0^T V(X,u) \d \mathbf{X}$ admits a directional derivative in the direction $h$, equal to $\int  \nabla^{\mathbb{R}^b}_h V(X,u)  \d \mathbf{X}$. As for the continuity of these directional derivatives, it is obtained by employing again the same bound: 
	\begin{align*}
		|	\int  \nabla^{\mathbb{R}^b}_h V(X,u)  \d \mathbf{X}-\int  \nabla^{\mathbb{R}^b}_h V(X,v)  \d \mathbf{X}|
		&\leq C_\mathbf{X} \|x\mapsto 	\nabla^{\mathbb{R}^b}_hV(x,u)-\nabla^{\mathbb{R}^b}_hV(x,v)\|_{\mathcal{C}^2_b}\\
		&	\leq b C_\mathbf{X} \max_{i,j,k} \omega_{|u-v|} (\partial_{i,j,k} V)\underset{v\to u}\longrightarrow 0.
	\end{align*}
\end{proof}

\section{Proof of the convergence \ref{eq:conv} } 
\label{sec:appendixTechnic}
Here we prove the technical estimation \ref{eq:conv}. For this, we start with an approximation lemma.
\begin{lemma}
	\label{le:mollifok} 
	Let $r\in(2,\infty)$. Then, $ ( \varphi^\epsilon\conv\varphi^{\epsilon'} -\delta )\conv G\in L^r(\mathbb{R}^2)$ for all $\epsilon,{\epsilon'}>0$, and 
	\[ \| ( \varphi^\epsilon\conv\varphi^{\epsilon'} -\delta )\conv G\|_{L^{r}(\mathbb{R}^2) }\underset{\epsilon,\epsilon'\to 0}\longrightarrow 0.\]
\end{lemma}
\begin{proof}
	We treat separately the integrals over $B_1$ and over $\mathbb{R}^2\setminus B_1$, we start with the one over $B_1$. 
	
	For an arbitrary $r\in (1,\infty)$, 
	\[ \| (( \varphi^\epsilon\conv\varphi^{\epsilon'} -\delta )\conv G) \mathbbm{1}_{B_1} \|_{L^r}
	\leq 
	\| ( \varphi^\epsilon\conv\varphi^{\epsilon'} -\delta )\conv (G \mathbbm{1}_{B_1}) \|_{L^r}
	+\| \varphi^\epsilon\conv\varphi^{\epsilon'}  \conv (G \mathbbm{1}_{B_1}) - ( \varphi^\epsilon\conv\varphi^{\epsilon'}\conv G) \mathbbm{1}_{B_1} \|_{L^r}.
	\]
	Since $G \mathbbm{1}_{B_1}\in L^r(\mathbb{R}^2)$ for all $r\in(1,\infty)$, the first term goes to $0$ by general property of convolutions with mollifier. As for the second term, remark the considered function 
	\[\varphi^\epsilon\conv\varphi^{\epsilon'}  \conv (G \mathbbm{1}_{B_1}) - ( \varphi^\epsilon\conv\varphi^{\epsilon'}\conv G) \mathbbm{1}_{B_1}\] is non-zero only on $B_{1+\epsilon+\epsilon'}\setminus B_{1-\epsilon-\epsilon'}$. On this set, it is bounded in absolute value by $2M_{\epsilon,\epsilon'}$, where $M_{\epsilon,\epsilon'}$ is the supremum of $G$ on $B_{1+2\epsilon+2\epsilon'}\setminus B_{1-2\epsilon-2\epsilon'}$. The function $(\epsilon,\epsilon')\to M_{\epsilon,\epsilon'}$ is non-decreasing in both variables, and finite for $\epsilon,\epsilon'\leq\frac{1}{8}$. Thus, there is a finite $M$ such that for all  $\epsilon,\epsilon'\leq\frac{1}{8}$, $M_{\epsilon,\epsilon'}\leq M$  and then 
	\[ 
	\| \varphi^\epsilon\conv\varphi^{\epsilon'}  \conv (G \mathbbm{1}_{B_1}) - ( \varphi^\epsilon\conv\varphi^{\epsilon'}\conv G) \mathbbm{1}_{B_1} \|_{L^r}\leq 2\pi \cdot 2(\epsilon+\epsilon')\cdot (1+\epsilon+\epsilon') \cdot  2M, \] 
	which goes to $0$ as $\epsilon\to 0$. 
	Hence, 
	\[\| (( \varphi^\epsilon\conv\varphi^{\epsilon'} -\delta )\conv G) \mathbbm{1}_{B_1} \|_{L^r}  \underset{\epsilon\to 0}\longrightarrow 0.
	\]
	We now consider the part of the integral on $\mathbb{R}^2\setminus B_1$, which is why we assumed $r>2$ in the lemma: since $\varphi\otimes \varphi$ is positive with integral $1$ over $B_{K_{\varphi}}^2$, for any measurable function $g$, \[(\int_{B_{K_{\varphi}}^2} g(v,w) \varphi(v) \varphi(w)\d v \d w)^r\leq \int_{B_{K_{\varphi}}^2} |g(v,w)|^r \varphi(v) \varphi(w)\d v \d w.\]
	Thus,
	\begin{align*} 
		\| (( \varphi^\epsilon\conv\varphi^{\epsilon'} -\delta )\conv G) \mathbbm{1}_{\mathbb{R}^2\setminus B_1} \|^r_{L^r}
		&= 
		\int_{\mathbb{R}^2\setminus B_1} \Big( \int_{B_{K_\varphi}^2} \varphi(v)\varphi(w) (   G(x)-G(x-\epsilon v-\epsilon'w ) ) \d v\d w \Big)^r \d x\\
		&\leq  
		\int_{\mathbb{R}^2\setminus B_1}  \int_{B_{K_\varphi}^2} \varphi(v)  \varphi(w)|   G(x)-G(x-\epsilon v-\epsilon'w )|^r \d v \d w \d x\\
		& =
		\int_{\mathbb{R}^2\setminus B_1} \int_{B_{K_\varphi}^2} \frac{\varphi(v) \varphi(w)}{(2\pi)^r} \big| \log|1+ \frac{\epsilon v+\epsilon' w}{|x|} |  \big|^r \d v\d w \d x\\
		&\leq 
		\int_{\mathbb{R}^2\setminus B_1} \int_{B_{K_\varphi}^2} \frac{\varphi(v)\varphi(w)}{(2\pi)^r} K_\varphi^r (\epsilon+\epsilon')^r |x|^{-r}  \d v\d w \d x\\
		&= \frac{K_\varphi^r (\epsilon+\epsilon')^r}{(2\pi)^{r-1} } 
		\int_1^\infty   \rho^{1-r}  \d\rho \underset{\epsilon,\epsilon' \to 0}\longrightarrow 0, \qquad \text{provided $r>2$,}
	\end{align*} 
which concludes the proof.
\end{proof}

\begin{lemma} \label{le:ap2} 
	\begin{equation}
		\text{For $i\in \{1,2\}$,}
		\int_0^{\bar{T}'} \big( \int_0^{\bar{T}}( (\varphi^\epsilon\conv\varphi^{\epsilon'} -\delta )\conv G)(W_s-W'_t)
		\mathrm{d} W^i_s \big)  \mathrm{d} W'^i_t\  \overset{L^2(\mathbb{P})}{\underset{\epsilon,\epsilon'\to 0}\longrightarrow} 0.
	\end{equation}
\end{lemma} 
\begin{proof}
	By symmetry we can assume $i=1$. 	We split the integrals into four parts, at $T'$ and $T$. Since the integral is invariant by reparametrisation, we can assume without loosing generality that $\bar{T}=T+1$ and that $W_t$ is parameterised by unit speed from $T$ to $\bar{T}=T+1$ (recall on this time interval, $W$ is not a Brownian motion by just a straight line segment from $W_T$ to $W_0$), i.e. for $t\in [T,T+1]$,  $W_t=W_T+ (t-T) (W_0-W_{T})$, and similarly for $W'$. In particular, for such a $t$, $\d W_t= (W_0-W_{T}) \d t$.
	Set $C=\mathbb{E}[(W_1-W_0)^{4}]^{\frac{1}{2}}<\infty$.
	\begin{itemize} 
		\item  Using \begin{align*} \mathbb{E}[ f(W,W')^2  (W_0-W_{T})^2 (W'_0-W'_{T'})^2   ]
			& \leq \mathbb{E}[ f(W,W')^4]^\frac12 
		\mathbb{E}[(W_0-W_{T})^4  (W'_0-W'_{T'})^4  ]^\frac12\\
		&= \mathbb{E}[ f(W,W')^4]^\frac12 C^2 T T',
		\end{align*}
		we get
		\begin{align*}
			&	\mathbb{E}\big[ \big( 
			\int_{T'}^{\bar{T}'} \big( \int_{T}^{\bar{T}} (\varphi^\epsilon\conv\varphi^{\epsilon'} -\delta )\conv G (W_s-W'_t)
			\mathrm{d} W^1_s \big)  \d {{W'}}\phantom{\hspace{-0.1cm})}^1_t\big)^2 \big]\\
			=& 
			\mathbb{E}\big[   (	\int_0^1\int_0^1  (\varphi^\epsilon\conv\varphi^{\epsilon'} -\delta )\conv G (s W_T- t W'_{T'} ) (W^1_T-W^1_0)(W'^1_{T'}-W'^1_0) \d s \d t )^2 \big] 
			\\
			\leq & 
			\mathbb{E}\big[   	\int_0^1\int_0^1 \big( (\varphi^\epsilon\conv\varphi^{\epsilon'} -\delta )\conv G (s W_T- t W'_{T'} )(W^1_T-W^1_0)(W'^1_{T'}-W'^1_0) \big)^2 \d s \d t  \big] \\
			\leq 	&	
			C^2 T T'  \big( 
			\int_0^1 \int_0^1 \mathbb{E}\big[   ( \varphi^\epsilon\conv\varphi^{\epsilon'} -\delta )\conv G (s W_T -t  W'_{T'})^4 \big] \d s \d t\big)^\frac12\\
			=& C^2 T T' \big(\int_0^1 \int_0^1 \int_{(\mathbb{R}^2)^2} p_T(0,z) p_{T'}(0,w)
			( \varphi^\epsilon\conv\varphi^{\epsilon'} -\delta )\conv G(s z-tw )^4 \d z \d w \d s \d t  \big)^\frac12\\
			=& C^2 T T' \big(\int_0^1 \int_0^1 \int_{(\mathbb{R}^2)^2} p_T(0,x/s) p_{T'}(0,y/t)
			( \varphi^\epsilon\conv\varphi^{\epsilon'} -\delta )\conv G(x-y )^4 \d x \d y \frac{\d s}{s^2} \frac{\d t}{t^2}  \big)^\frac12.
		\end{align*}
		Ley $p\in(1,2)$ and $q\in(1,\infty)$ such that $2p{-1}+q^{-1}=2$, for example $p=q=\frac{3}{2}$. Then, using the convolution inequality 
		\begin{equation}
			\label{eq:convfgh}
		\int_{(\mathbb{R}^2)^2} f(x) g(x-y) h(y)  \d x \d y\leq \| f\|_{L^p}\| g\|_{L^p}\| f\|_{L^q}, 
		\end{equation}
		and the equality 
		\[ \| x\mapsto p_T(0,x/s)\|_{L^p}=
		\big( \int_{\mathbb{R}^2}  p_T(0,x/s)^p \d x 
		\big)^{\frac1p} = \big( \int_{\mathbb{R}^2}  p_T(0,z)^p s^2 \d z \big)^\frac1p =s^{\frac2p} C_{p,T}, 
		\]
		we deduce 
		\begin{align*} 
			&	\int_0^1 \int_0^1 \int_{(\mathbb{R}^2)^2} p_T(0,x/s) p_{T'}(0,y/t)
			( \varphi^\epsilon\conv\varphi^{\epsilon'} -\delta )\conv G(x-y )^4 \d x \d y \frac{\d s}{s^2} \frac{\d t}{t^2}  \\ 
			&\leq C_{p,T}C_{p,T'} \int_0^1 \int_0^1 s^{\frac2p-2} t^{\frac2p-2} \| (( \varphi^\epsilon\conv\varphi^{\epsilon'} -\delta )\conv G) \|_{L^{4r} }^4 \d s \d t\\ 
			&\leq C'\| (( \varphi^\epsilon\conv\varphi^{\epsilon'} -\delta )\conv G)  \|_{L^{4q} }^4 \underset{\epsilon,\epsilon' \to 0}\longrightarrow 0, \qquad \text{by Lemma \ref{le:mollifok}.}
		\end{align*}

		\item The other pieces are estimated rather similarly. 
		\begin{align*}
			&	\mathbb{E}\big[ \big( 
			\int_{T'}^{\bar{T'}} \big( \int_0^{T} (\varphi^\epsilon\conv\varphi^{\epsilon'} -\delta )\conv G (W_s-W'_t)
			\mathrm{d} W^1_s \big)  \d W'^1_t\big)^2 \big]\\
			=& 
			\int_0^1\int_0^T \mathbb{E}\big[   ( \varphi^\epsilon\conv\varphi^{\epsilon'} -\delta )\conv G (W_s- t W'_{T'} )^2 (W'^1_{T'} )^2 \big] \d s \d t
			\\
			\leq & 
			C T'
			\int_0^T \sup_a \big( \int_{\mathbb{R}^2}  p_t(a,z)   ( \varphi^\epsilon\conv\varphi^{\epsilon'} -\delta )\conv G (z )^4 \d z\big)^\frac12 \d s,   
		\end{align*}
		which goes to $0$ as $\epsilon,\epsilon' \to 0$ by a simpler version of the same argument above. Remark applying the convolution inequality \eqref{eq:convfgh} immediately gets rid of the supremum over $a$. 
		Swapping the integrals, this also deals with the term \[ 
		\mathbb{E}\big[ \big( 
		\int_{0}^{T'} \big( \int_T^{\bar{T}} (\varphi^\epsilon\conv\varphi^{\epsilon'} -\delta )\conv G (W_s-W'_t)
		\mathrm{d} W^1_s \big)  \mathrm{d} {{W'}}\phantom{\hspace{-0.1cm})}^1_t\big)^2 \big].\]
		\item The last term 
		\begin{multline*}
			\mathbb{E}\big[ \big( 
			\int_{0}^{T'} \big( \int_0^{T} (\varphi^\epsilon\conv\varphi^{\epsilon'} -\delta )\conv G (W_s-W'_t)
			\mathrm{d} W^1_s \big)  \mathrm{d} {{W'}}\phantom{\hspace{-0.1cm})}^1_t\big)^2 \big]
			\\=   \int_{0}^{T'} \int_0^T \mathbb{E}\big[ 
			(\varphi^\epsilon\conv\varphi^{\epsilon'} -\delta )\conv G (W_s-W'_t)^2 \big] \d s \d t
		\end{multline*}
		is also treated this way.\vspace{-0.6cm}
	\end{itemize}
\end{proof}

\section*{Acknowledgments}
I am deeply grateful to Pierre Perruchaud who first introduced me to differential geometry and to  Thierry Lévy who first introduced me to the Yang--Mills field. This paper would not have been possible without the knowledge they shared with me. 
 
\section*{Fundings}  
This work was supported by the EPSRC grants no
EP/R014604/1 and EP/W006227/1 .
I would like to thank the Isaac Newton Institute for Mathematical Sciences,
Cambridge, for the support and hospitality during the programme Self-interacting Processes
where work on this paper was undertaken.

\bibliographystyle{plain}

\begin{thebibliography}{10}
	
	\bibitem{AlbeverioKusuoka}
	Sergio Albeverio and Shigeo Kusuoka.
	\newblock A basic estimate for two-dimensional stochastic holonomy along
	{B}rownian bridges.
	\newblock {\em J. Funct. Anal.}, 127(1):132--154, 1995.
	
	\bibitem{Banchoff}
	Thomas~F. Banchoff and William~F. Pohl.
	\newblock A generalization of the isoperimetric inequality.
	\newblock {\em J. Differential Geometry}, 6:175--192, 1971/72.
	
	\bibitem{chandra}
	Ajay Chandra and Ilya Chevyrev.
	\newblock Gauge field marginal of an abelian {H}iggs model, 2022.
	
	\bibitem{Dynkin}
	E.~B. Dynkin.
	\newblock Regularized self-intersection local times of planar {B}rownian
	motion.
	\newblock {\em Ann. Probab.}, 16(1):58--74, 1988.
	
	\bibitem{magneticFK}
	L{\'a}szl{\'o} Erd{\H o}s.
	\newblock Magnetic {L}ieb-{T}hirring inequalities.
	\newblock {\em Comm. Math. Phys.}, 170(3):629--668, 1995.
	
	\bibitem{FrizHairer}
	Peter~K. Friz and Martin Hairer.
	\newblock {\em A course on rough paths}.
	\newblock Universitext. Springer, Cham, second edition, [2020] \copyright 2020.
	\newblock With an introduction to regularity structures.
	
	\bibitem{Geman}
	Donald Geman, Joseph Horowitz, and Jay Rosen.
	\newblock A local time analysis of intersections of {B}rownian paths in the
	plane.
	\newblock {\em Ann. Probab.}, 12(1):86--107, 1984.
	
	\bibitem{Hutton}
	James~E. Hutton and Paul~I. Nelson.
	\newblock Interchanging the order of differentiation and stochastic
	integration.
	\newblock {\em Stochastic Process. Appl.}, 18(2):371--377, 1984.
	
	\bibitem{Karakostas}
	George~L. Karakostas.
	\newblock An extension of {H}\"older's inequality and some results on infinite
	products.
	\newblock {\em Indian J. Math.}, 50(2):303--307, 2008.
	
	\bibitem{LeGall}
	Jean-Fran\c{c}ois Le~Gall.
	\newblock Exponential moments for the renormalized self-intersection local time
	of planar {B}rownian motion.
	\newblock In {\em S\'{e}minaire de {P}robabilit\'{e}s, {XXVIII}}, volume 1583
	of {\em Lecture Notes in Math.}, pages 172--180. Springer, Berlin, 1994.
	
	\bibitem{Thierry}
	Thierry L\'evy.
	\newblock Yang-{M}ills measure on compact surfaces.
	\newblock {\em Mem. Amer. Math. Soc.}, 166(790):xiv+122, 2003.
	
	\bibitem{MansuyYor}
	Roger Mansuy and Marc Yor.
	\newblock {\em Aspects of {B}rownian motion}.
	\newblock Universitext. Springer-Verlag, Berlin, 2008.
	
	\bibitem{Nualart}
	David Nualart.
	\newblock Stochastic anticipating calculus.
	\newblock In {\em Probability towards 2000 ({N}ew {Y}ork, 1995)}, volume 128 of
	{\em Lect. Notes Stat.}, pages 249--262. Springer, New York, 1998.
	
	\bibitem{zeta}
	Pierre Perruchaud and Isao Sauzedde.
	\newblock Loop soup representation of zeta-regularised determinants and
	equivariant symanzik identities, 2024.
	\newblock to appear in Annals of probability.
	
	\bibitem{Rosenberg}
	Steven Rosenberg.
	\newblock {\em The {L}aplacian on a {R}iemannian manifold}, volume~31 of {\em
		London Mathematical Society Student Texts}.
	\newblock Cambridge University Press, Cambridge, 1997.
	\newblock An introduction to analysis on manifolds.
	
	\bibitem{LAWA}
	Isao Sauzedde.
	\newblock {Lévy area without approximation}.
	\newblock {\em Annales de l'Institut Henri Poincaré, Probabilités et
		Statistiques}, 58(4):2165 -- 2200, 2022.
	
	\bibitem{Symanzik}
	Kurt Symanzik.
	\newblock Euclidean quantum field theory. {I}. {E}quations for a scalar model.
	\newblock {\em J. Mathematical Phys.}, 7:510--525, 1966.
	
	\bibitem{StochasticFubini}
	Mark Veraar.
	\newblock The stochastic {F}ubini theorem revisited.
	\newblock {\em Stochastics}, 84(4):543--551, 2012.
	
	\bibitem{Werner3}
	Wendelin Werner.
	\newblock Rate of explosion of the {Amperean} area of the planar {Brownian}
	loop.
	\newblock In {\em S\'eminaire de Probabilit\'es XXVIII}, pages 153--163.
	Berlin: Springer, 1994.
	
	\bibitem{Werner2}
	Wendelin Werner.
	\newblock Formule de {G}reen, lacet brownien plan et aire de {L}\'{e}vy.
	\newblock {\em Stochastic Process. Appl.}, 57(2):225--245, 1995.
	
	\bibitem{Zakai}
	Moshe Zakai.
	\newblock Some moment inequalities for stochastic integrals and for solutions
	of stochastic differential equations.
	\newblock {\em Israel J. Math.}, 5:170--176, 1967.
	
\end{thebibliography}

%
%
%

\end{document}